\nonstopmode

\documentclass[]{svmult}
\usepackage{amsmath}          
\usepackage{verbatim}% for comment enviroment
\usepackage{amssymb}
\usepackage{graphicx}
\usepackage{epsfig}
\allowdisplaybreaks
\usepackage[all]{xy}% xymatrix for CD's
\newdir{ >}{{}*!/-10pt/@{>}}% tail for mono arrow
\newdir^{ (}{{}*!/-5pt/@^{(}}% tail for hookarrow
\newdir_{ (}{{}*!/-5pt/@_{(}}% tail for hookarrow
\def\cgaps#1{}
\def\Cgaps#1{}
\def\AMStag#1{}
\def\AMSunderset#1\to#2{\underset{#1}{#2}}
\def\AMSoverset#1\to#2{\overset{#1}{#2}}
\def\undersetbrace#1\to#2{\underbrace{#2}_{#1}}
\def\oversetbrace#1\to#2{\overbrace{#2}^{#1}}
\def\3{\ss}
\def\spcheck{\check{}}

\def\therosteritem#1{(#1)}
\def\Cal{\mathcal}
\newenvironment{proclaim}[1]{\par\medskip\noindent{\bf #1.}\it\newline}{\par\smallskip}
\newenvironment{demo}[1]{\par\smallskip\noindent{\bf #1.}}{\par\smallskip}
\newcommand{\nmb}[2]{%
\if0#1{#2\setcounter{section}{0}}%
\else\if1#1{#2\addtocounter{section}{1}}%
\else\if:#1{#2}%
\else\if.#1{#2}%
\else\if!#1{(#2)}%
\else\if|#1{(#2)}%
\else{(#2)}%
\fi\fi\fi\fi\fi\fi}
\def\East#1#2{-\raisebox{0.1pt}{$\mkern-16mu\frac{\;\;#1\;}{\;\;#2\;}\mkern-16mu$}\to}
\def\West#1#2{\gets\raisebox{0.1pt}{$\mkern-16mu\frac{\;#1\;\;}{\;#2\;\;}\mkern-16mu$}-}
\def\cit#1#2{\ifx#1!\cite{2}\else#2\fi} %for citing references 
\def\idx{}               % for producing index, invoked by kms-book.sty 
\def\ign#1{}             %=ignore, invisible entry for the index only 
\def\o{\circ} 
\def\X{\mathfrak X} 
\def\al{\alpha} 
\def\be{\beta} 
\def\ga{\gamma} 
\def\de{\delta} 
\def\ep{\varepsilon} 
\def\ze{\zeta} 
\def\et{\eta} 
\def\th{\theta} 
 
\def\ka{\kappa} 
\def\la{\lambda} 
\def\rh{\rho} 
\def\si{\sigma} 
\def\ta{\tau} 
\def\ph{\varphi} 
\def\ch{\chi} 
\def\ps{\psi} 
\def\om{\omega} 
\def\Ga{\Gamma}

\def\Om{\Omega} 
\def\i{^{-1}} 
\def\x{\times} 
\def\g{{\mathfrak g}} 
\def\p{\partial}
\def\ad{\operatorname{ad}} 
\def\Ad{\operatorname{Ad}} 
\def\L{\mathcal{L}}
\def\on{\operatorname}
\def\today{\ifcase\month\or 
 January\or February\or March\or April\or May\or June\or 
 July\or August\or September\or October\or November\or December\fi 
 \space\number\day, \number\year} 
\title*{Some Geometric Evolution Equations
\\
Arising as Geodesic Equations on Groups of Diffeomorphisms 
Including the Hamiltonian Approach 
} 
\titlerunning{Geometric Evolution Equations}
\author{Peter W. Michor 
\thanks{Supported by `Fonds zur
F\"orderung der wissenschaftlichen                    
Forschung', Projekt P~17108. 
Work partly done at the Program for Evolutionary Dynamics, Harvard
University}
} 
\institute{
Fakult\"at f\"ur Mathematik, Universit\"at Wien, 
Nordbergstrasse 15, A-1090 Wien, Austria; and  
Erwin Schr\"odinger International Institute of Mathematical Physics,   
Boltzmanngasse 9, A-1090 Wien, Austria. 
\\
peter.michor@esi.ac.at
}

\begin{document}

\maketitle


\section*{Table of contents } 

\noindent Introduction
\par
\noindent 1. A general setting and a motivating example 
\par 
\noindent 2. Weak symplectic manifolds 
\par 
\noindent 3. Right invariant weak Riemannian metrics on Lie groups 
\par 
\noindent 4. The Hamiltonian approach
\par 
\noindent 5. Vanishing geodesic distance on groups of diffeomorphisms
\par 
\noindent 6. The regular Lie group of rapidly decreasing diffeomorphisms 
\par 
\noindent 7. The diffeomorphism group of $S^1$ or $\mathbb R$, and Burgers' hierarchy 
\par 
\noindent 8. The Virasoro-Bott group and the Korteweg-de Vries hierarchy 
\par 
\noindent Appendix A. Smooth Calculus beyond Banach spaces
\par 
\noindent Appendix B. Regular infinite dimensional Lie groups
\par
\noindent References
\par

\section*{Introduction}

This is the extended version of a lecture course given at the University of
Vienna in the spring term 2005. Many thanks to the audience of this course
for many keen questions. The main aim of this course was to understand the
papers \cite{10} and \cite{11}.

The purpose of this review article is to give a complete account of existence and
uniqueness of the solutions of the members of higher order of the
hierarchies of Burgers' equation and the Korteweg-de Vries equation,
including their derivation and all the necessary background. 
We do this both on the circle, and on the real line in the setting of
rapidly decreasing functions. 
These are all
geodesic equations of infinite dimensional regular Lie groups, namely the
diffeomorphism group of the line or the circle and the corresponding
Virasoro group. 

Let us describe the content: Appendix A is a short description of
convenient calculus in infinite dimensions (beyond Banach spaces) where
everything is based on smooth curves: A mapping is $C^\infty$ if it maps
smooth curves to smooth curves. It is a theorem that smooth curves in a 
space of smooth functions are just smooth functions of one variable more;
this is the basic assumption of variational calculus. 
Appendix B gives a short account of infinite dimensional regular Lie
groups.
Here regularity means that a smooth curve in the Lie algebra can be
integrated to a smooth curve in the group whose right (or left) logarithmic
derivative equals the given curve. No infinite dimensional Lie group is known
which is not regular.
Section 1, as a motivating example, computes the geodesics and the
curvature of the most naive Riemannian metric on the space of embeddings of
the real line to itself and shows that this can be converted into Burgers'
equation. 
Section 2 treats Hamiltonian mechanics on infinite dimensional weak
symplectic manifolds. Here `weak' means that the symplectic 2-form is
injective as a mapping from the tangent bundle to the cotangent bundle. 
Section 3 computes geodesics and curvatures of right invariant Riemannian
metrics on regular Lie groups as done by Arnold \cite{1}. 
Section 4 redoes this in the symplectic
approach and computes the associated momentum mappings and conserved
quantities. 
Section 5 shows that the geodesic distance vanishes on any full
diffeomorphis group for the right invariant metric coming from the
$L^2$-metric on the Lie algebra of vector fields for a given Riemannian
metric on a manifold. In particular, Burgers' equation is the geodesic
equation of such a metric.
Section 6 treats the group of diffeomorphisms of the real line which
decrease rapidly to the identity as a regular Lie group. This will be
important for Burgers' equation as geodesic equation on this group, and
also for the KdV equation.
Here we also give a short presentation of Sobolev spaces on the real line
and of the scale of $HC^n$-spaces for which we were
able to give simple proofs of the results which we need later. 
Section 7 treats geodesic equations on the diffeomorphism groups of the real
line or $S^1$ which leads to Burgers' hierarchy. We solve these
equations starting at certain higher order, following \cite{11}. 
Section 8 does this for the Virasoro groups on the real line or $S^1$. For
the solution of the higher order equations we follow \cite{10}.

Note that in this paper we concentrate on in the smooth ($=C^\infty$) aspect.
We also do not treat complete integrability for Burgers' and KdV equation,
although we prepared almost all of the necessary background.

\section*{\nmb0{1}. A general setting and a motivating example  
}       
 
\subsection*{\nmb.{1.1}. The principal bundle of embeddings}
 
Let $M$ and $N$ be smooth finite dimensional manifolds, connected and  
second countable without boundary, such that $\dim M\leq\dim N$. 
Then the space $\on{Emb}(M,N)$ of all embeddings (immersions which are 
homeomorphisms on their images) from $M$ into $N$ is an 
open submanifold of $C^\infty(M,N)$ which is stable under the right 
action of the diffeo\-mor\-phism group of $M$. 
Here $C^\infty(M,N)$ is a smooth manifold modeled on spaces of 
sections with compact support $\Ga_c(f^*TN)$. In particular 
the tangent space at 
$f$ is canonically iso\-mor\-phic to the space of vector fields along $f$ 
with compact support in $M$. If $f$ and $g$ differ on a non-compact 
set then they belong to different connected components of 
$C^\infty(M,N)$. See \cite{23} and \cite{27}.
Then $\on{Emb}(M,N)$ is the 
total space of a smooth principal fiber bundle with structure 
group the diffeo\-mor\-phism group of $M$; the base is called $B(M,N)$, it 
is a Hausdorff smooth manifold modeled on nuclear (LF)-spaces. 
It can be thought of as the "nonlinear Grassmannian" or "differentiable
Chow variety" of all
submanifolds of $N$ which are of type $M$. This result is based 
on an idea implicitly contained in \cite{39}, it
was fully proved in \cite{6} for compact $M$ 
and for general $M$ in \cite{26}. See also \cite{27}, section 13 and
\cite{23}.
If we take a Hilbert space $H$ instead of $N$, then  $B(M,H)$ is 
the classifying space for $\on{Diff}(M)$ if $M$ is compact, and the 
classifying bundle $\on{Emb}(M,H)$ carries also a universal 
connection. This is shown in \cite{28}. 
 
\subsection*{\nmb.{1.2} } 
If $(N,g)$ is a Riemannian manifold then on the manifold $\on{Emb}(M,N)$  
there is a naturally induced weak Riemannian metric 
given, for  $s_1,s_2\in\Ga_c(f^*TN)$ and $\ph\in \on{Emb}(M,N)$, by
\begin{displaymath}
G_\phi(s_1,s_2)=\int_Mg(s_1,s_2)\operatorname{vol}(\phi^*g), \quad
\phi \in \on{Emb}(M,N),  
\end{displaymath}
where $\operatorname{vol}(g)$ denotes the volume form on $N$ induced
by the Riemannian metric $g$ and $\operatorname{vol}(\phi^*g)$ the
volume form on $M$ induced by the pull back metric $\phi^*g$. The
covariant derivative and curvature of the Levi-Civita connection
induced by $G$ were investigated in
\cite{5} if $N=\mathbb R^{\dim M+1}$ (endowed with the standard
inner product) and in \cite{18} for the general case. 
In \cite{30} it was shown that the geodesic distance (topological metric) 
on the base manifold
$B(M,N)=\on{Emb}(M,N)/\on{Diff}(M)$ induced by this Riemannian metric
vanishes.
 
This weak Riemannian  
metric is invariant under the action of the diffeo\-mor\-phism group  
$\on{Diff}(M)$ by composition from the right and hence it induces a  
Rie\-mannian metric on the base manifold $B(M,N)$. 
 
\subsection*{\nmb.{1.3}. Example }
Let us consider the special case $M=N=\mathbb R$, that is, the space
$\on{Emb}(\mathbb R,\mathbb R)$  of all embeddings of the real line into
itself, which contains the   diffeomorphism group $\on{Diff}(\mathbb R)$ as
an open subset.  The case $M=N=S^1$ is treated in a similar fashion
and the results of this paper are also valid in this situation, where 
$\on{Emb}(S^1,S^1)=\on{Diff}(S^1)$. 
For our purposes, we may restrict attention to
the space of orientation-preserving embeddings, denoted by
$\on{Emb}^{+}(\mathbb R,\mathbb R)$.
The weak Riemannian metric has thus the expression 
\begin{displaymath}
G_f(h,k) = \int_\mathbb R h(x)k(x)|f'(x)|\,dx,\quad  
f\in\on{Emb}(\mathbb R,\mathbb R), \quad h,k\in C^\infty_c(\mathbb R,\mathbb R). 
\end{displaymath}
We shall compute the geodesic equation for this metric by variational  
calculus. The energy of a curve $f$ of embeddings is
\begin{displaymath}
E(f) = \tfrac12\int_a^b G_f(f_t,f_t) dt  
     = \tfrac12\int_a^b \!\!\int_{\mathbb R} f_t^2f_x\, dx dt. 
\end{displaymath}
If we assume that $f(x,t,s)$ is a smooth function and that the
variations are with fixed endpoints, then the derivative with
respect to $s$ of the energy is  
\begin{align*} 
\p s|_0 E(f(\quad,s)) 
&= \p s|_0 \tfrac12\int_a^b \!\!\int_{\mathbb R}
         f_t^2f_x\, dx dt
\\&
= \tfrac12\int_a^b \!\!\int_{\mathbb R}
      (2f_tf_{ts}f_x+f_t^2f_{xs})dx dt
\\&
=  -\tfrac12\int_a^b \!\!\int_{\mathbb R}  
       (2f_{tt}f_{s}f_x+2f_{t}f_{s}f_{tx}+2f_tf_{tx}f_s)dx dt
\\&
= -\int_a^b \!\!\int_{\mathbb R}  
     \left(f_{tt}+2\frac{f_{t}f_{tx}}{f_x}\right)f_{s}f_xdx dt,
\\
\end{align*}so that the geodesic equation with its initial data is:  
\begin{align} 
f_{tt}&=-2\frac{f_{t}f_{tx}}{f_x},\quad  
     f(\quad,0)\in \on{Emb}^+(\mathbb R,\mathbb R), \quad 
     f_t(\quad,0)\in C^\infty_c(\mathbb R,\mathbb R)\tag{\nmb:{1}}
\\&
=: \Ga_f(f_t,f_t), 
\notag\end{align}where the Christoffel symbol  
$\Ga:\on{Emb}(\mathbb R,\mathbb R)\x C^\infty_c(\mathbb R,\mathbb R)\x  
C^\infty_c(\mathbb R,\mathbb R)\to C^\infty_c(\mathbb R,\mathbb R)$ 
is given by symmetrisation:
\begin{equation}
\Ga_f(h,k) := -\frac{hk_x+h_xk}{f_x}= -\frac{(hk)_x}{f_x}\,.\tag{\nmb:{2}} 
\end{equation}
For vector fields $X,Y$ on $\on{Emb}(\mathbb R,\mathbb R)$ the  
covariant derivative is given by the expression  
$\nabla^{\on{Emb}}_XY= dY(X)-\Ga(X,Y)$. 
The Riemannian curvature  
$R(X,Y)Z = (\nabla_X\nabla_Y - \nabla_Y\nabla_X -\nabla_{[X,Y]})Z$ 
is then determined in terms of the Christoffel form by
\begin{align*}
R&(X,Y)Z = (\nabla_X\nabla_Y - \nabla_Y\nabla_X -\nabla_{[X,Y]})Z
\\&
= \nabla_X(dZ(Y)-\Ga(Y,Z)) - \nabla_Y(dZ(X)-\Ga(X,Z)) 
\\&\qquad\qquad     - dZ([X,Y]) + \Ga([X,Y],Z)
\\&
= d^2Z(X,Y) + dZ(dY(X)) - \Ga(X,dZ(Y)) 
\\&
\quad- d\Ga(X)(Y,Z) - \Ga(dY(X),Z) - \Ga(Y,dZ(X)) + \Ga(X,\Ga(Y,Z))
\\&
\quad- d^2Z(Y,X) - dZ(dX(Y)) + \Ga(Y,dZ(X)) 
\\&
\quad+ d\Ga(Y)(X,Z) + \Ga(dX(Y),Z) + \Ga(X,dZ(Y)) - \Ga(Y,\Ga(X,Z))
\\&
\quad- dZ(dY(X)-dX(Y)) + \Ga(dY(X)-dX(Y),Z)
\\&
= - d\Ga(X)(Y,Z) + \Ga(X,\Ga(Y,Z)) + d\Ga(Y)(X,Z) - \Ga(Y,\Ga(X,Z)
\\
\end{align*}so that
\begin{align*} 
&R_f(h,k)\ell =
\\&
= -d\Ga(f)(h)(k,\ell) + d\Ga(f)(k)(h,\ell) + 
     \Ga_f(h,\Ga_f(k,\ell)) - \Ga_f(k,\Ga_f(h,\ell))
\\&
= -\frac{h_x(k\ell)_x}{f_x^2} +\frac{k_x(h\ell)_x}{f_x^2} 
     +\frac{\left( h\frac{(k\ell)_x}{f_x} \right)_x}{f_x} 
     -\frac{\left( k\frac{(h\ell)_x}{f_x} \right)_x}{f_x}
\tag{\nmb:{3}}\\&
= \frac1{f_x^3}\Bigl( 
     f_{xx}h_xk\ell - f_{xx}hk_x\ell  
     +f_xhk_{xx}\ell - f_xh_{xx}k\ell 
     +2f_xhk_x\ell_x -2f_xh_xk\ell_x\Bigr). 
\end{align*}
Now let us consider the trivialisation of $T\on{Emb}(\mathbb R,\mathbb R)$ by  
right translation (this is most useful for $\on{Diff}(\mathbb R)$). The
derivative of the inversion $\operatorname{Inv}:g\mapsto  g\i$ is given by
$$
T_g(\operatorname{Inv})h = - T(g\i)\o h \o g\i = 
-\frac{h\o g\i}{g_x\o g\i}
$$ 
for $g \in \on{Emb}(\mathbb R, \mathbb R),\, h \in C^\infty_c(\mathbb R,\mathbb R)$.
Defining  
\begin{displaymath}
u:= f_t\o f\i,\quad\text{ or, in more detail, }\quad
u(t,x)=f_t(t,f(t,\quad)\i(x)) ,
\end{displaymath}
we have
\begin{align*} 
u_x &= (f_t \o f\i)_x = (f_{tx}\o f\i)\frac1{f_x\o f\i}
     = \frac{f_{tx}}{f_x} \o f\i,
\\
u_t &= (f_t \o f\i)_t = f_{tt}\o f\i + (f_{tx}\o f\i)(f\i)_t 
\\&
= f_{tt}\o f\i - (f_{tx}\o f\i)\frac{1}{f_x\ f\i}(f_t\ f\i)
\end{align*}which, by {\nmb|{1}} and the first equation becomes
\begin{displaymath}
u_t = f_{tt}\o f\i - \left( \frac{f_{tx}f_t}{f_x} \right)\o f\i
     = -3\left(\frac{f_{tx}f_t}{f_x}\right)\o f\i 
     = -3 u_xu.
\end{displaymath}
The geodesic equation on $\on{Emb}(\mathbb R, \mathbb R)$ in right
trivialization, that is, in Eulerian formulation, is hence  
\begin{equation}
u_t = -3 u_xu \,,\tag{\nmb:{4}} 
\end{equation}
which is just Burgers' equation. 

Finally let us solve Burgers' equation and also describe its universal
completion, see see \cite{9}, \cite{A}, or \cite{KhM}.

In $\mathbb R^2$ with coordinates
$(x,y)$ consider the vector field
$Y(x,y)=(3y,0)=3y\p_x$ with differential equation $\dot x=3y, \dot y=0$.
It has the complete flow $\on{Fl}^Y_t(x,y)=(x+3ty,y)$. 

Let now $t\mapsto u(t,x)$ be a curve of functions on $\mathbb R$. 
We ask when the graph of $u$ can
be reparametrized in such a way that it becomes a solution curve of the
push forward vector field $Y_*:f\mapsto Y\o f$ on the space of embeddings
$\on{Emb}(\mathbb R,\mathbb R^{2})$. 
Thus consider a time dependent reparametrization $z\mapsto x(t,z)$, 
i.e., $x\in C^\infty(\mathbb R^{2},\mathbb R)$.  
The curve $t\mapsto (x(t,z),u(t,x(z,t)))$ in $\mathbb R^{2}$ is an 
integral curve of $Y$ if and only if 
\begin{align*} 
&\begin{pmatrix} 3u\o x \\
          0 \end{pmatrix}
=\p_t\begin{pmatrix} x \\
          u\o x \end{pmatrix}
=\begin{pmatrix} x_t \\
          u_t\o x + (u_x\o x)\cdot x_t
         \end{pmatrix} 
\\&
\Longleftrightarrow 
{\begin{cases} x_t = 3u\o x \\  
         0 = (u_t + 3u u_x)\o x
           \end{cases}} 
\end{align*}This implies that the 
graph of $u(t,\cdot)$, namely the curve $t\mapsto (x\mapsto (x,u(t,x)))$, 
may be parameterized as a solution curve 
of the vector field $Y_*$ on the space of embeddings 
$\on{Emb}(\mathbb R,\mathbb R^{2})$
starting at $x\mapsto (x,u(0,x))$ 
if and only if $u$ is a solution of the partial differential equation 
$u_t+ 3uu_x=0$. 
The  parameterization  $z\mapsto x(z,t)$ is then given by  
$x_t(z,t)=3u(x(t,z))$ with $x(0,z) = z\in \mathbb R$.
\begin{figure} \begin{center}
\epsfig{width=8cm,file=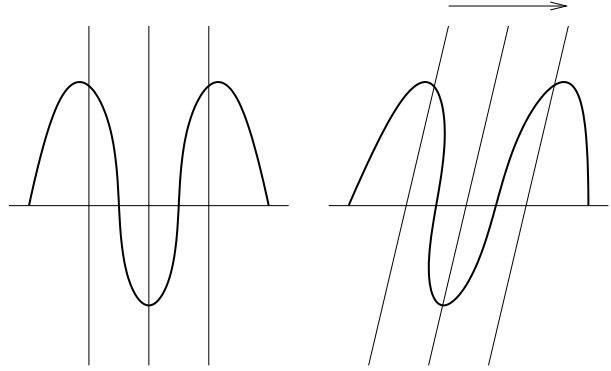}
\end{center} 
\caption{The characteristic flow of
the inviscid Burgers' equation tilts the plane.}
\end{figure}

This has a simple physical meaning. Consider
freely flying particles in $\mathbb R$, and trace a trajectory $x(t)$ of one
of the particles. Denote the velocity of a particle at the position
$x$ at the moment $t$ by $u(t,x)$, or
rather, by $3u(t,x):=\dot x(t)$. 
Due to the absence of interaction, the Newton
equation of any particle is $\ddot x(t)=0$. 

Let us illustrate this:
The flow of the vector 
field $Y=3u\p_x$ is tilting the plane to the right with constant speed. The 
illustration shows how a graph of an honest function is moved through a 
shock (when the derivatives become infinite) towards the graph of a
multivalued function; each piece of it is still a local solution.

\section*{\nmb0{2}. Weak symplectic manifolds } 

\subsection*{\nmb.{2.1}. Review } For a finite dimensional 
symplectic manifold $(M,\om)$ we have the following exact sequence of
Lie algebras:
\begin{displaymath}
0\to  H^0(M)\to C^\infty(M,\mathbb R)  
     \East{{\on{grad}^\om}}{}\X(M,\om) \East{}{}H^1(M)\to 0.
\end{displaymath}
Here $H^*(M)$ is the real De Rham cohomology of $M$, the space
$C^\infty(M,\mathbb R)$ is equipped with the Poisson bracket $\{\quad,\quad\}$,
$\X(M,\om)$ consists of all vector fields $\xi$ with $\L_\xi\om=0$
(the locally Hamiltonian vector fields), which is a Lie algebra
for the Lie bracket. Furthermore, 
$\on{grad}^\om f$ is the Hamiltonian vector field
for $f\in C^\infty(M,\mathbb R)$ given by $i(\on{grad}^\om f)\om = df$ and 
$\ga(\xi) = [i_\xi\om]$. The spaces $H^0(M)$ and $H^1(M)$ are
equipped with the zero bracket.

Consider a symplectic right action $r:M\x G\to  M$ of a connected
Lie group $G$ on $M$; we use the notation $r(x,g)=r^g(x)=r_x(g)=x.g$.
By $\ze_X(x)=T_e(r_x)X$ we get
a mapping $\ze:\g \to  \X(M,\om)$ which sends each element
$X$ of the Lie algebra $\g$ of $G$ to the fundamental vector
field $X$. This is a Lie algebra homomorphism (for right actions!).  
\begin{displaymath}
\cgaps{1.5;0.8;0.8;1.5}\xymatrix{
H^0(M) \ar[r]^{i} & 
     C^\infty(M,\mathbb R)  \ar[rr]^{\on{grad}^\om} & & 
     \X(M,\om) \ar[r]^{\ga} & H^1(M)
\\
 & & \g \ar[lu]^{j} \ar[ru]_{\ze} & &
}\end{displaymath}
A linear lift $j:\g\to  C^\infty(M,\mathbb R)$ of $\ze$  with 
$\on{grad}^\om\o j=\ze$ exists if
and only if $\ga\o \ze=0$ in $H^1(M)$. This lift $j$ may be
changed to a Lie algebra homomorphism if and only if the
$2$-cocycle $\bar \jmath:\g\x\g\to  H^0(M)$, given by 
$(i\o\bar \jmath)(X,Y) = \{j(X),j(Y)\} - j([X,Y])$, vanishes
in the Lie algebra cohomology $H^2(\g,H^0(M))$, 
for if $\bar\jmath = \de\al$ then $j-i\o\al$
is a Lie algebra homomorphism.

If $j:\g\to C^\infty(M,\mathbb R)$ is a Lie algebra homomorphism, we may
associate the {\it moment mapping} 
$\mu:M\to \g'=L(\g,\mathbb R)$ to it, which is 
given by
$\mu(x)(X) = \ch(X)(x)$ for $x\in M$ and $X\in \g$.
It is $G$-equivariant for a suitably chosen (in general affine)
action of $G$ on $\g'$.

\subsection*{\nmb.{2.2}} We now want to carry over to infinite 
dimensional manifolds 
the procedure of subsection \nmb!{2.1}. First we need the appropriate notions in 
infinite dimensions. So let $M$ be a manifold, which in general is 
infinite dimensional.

A $2$-form  $\om\in\Om^2(M)$ is called a \idx{\it weak symplectic 
structure} on $M$ if it is closed $(d\om=0)$ and if 
its associated vector bundle homomorphism
$ \om\spcheck: TM \to T^*M$ is injective. 

A $2$-form  $\om\in\Om^2(M)$ is called a \idx{\it strong 
symplectic structure} on $M$ if it is closed $(d\om=0)$ and if 
its associated vector bundle homomorphism
$ \om\spcheck: TM \to T^*M$ is invertible with smooth 
inverse. In this case, the vector bundle $TM$ has reflexive fibers 
$T_xM$: Let $i:T_xM\to (T_xM)''$ be the canonical mapping onto the 
bidual. Skew symmetry of $\om$ is equivalent to the fact that the 
transposed $(\om\spcheck)^t=(\om\spcheck)^*\o i:T_xM\to (T_xM)'$ satisfies 
$(\om\spcheck)^t=-\om\spcheck$. Thus, $i=-((\om\spcheck)\i)^*\o  \om\spcheck$  is an 
isomorphism.

\subsection*{\nmb.{2.3}}
Every cotangent bundle $T^*M$, viewed as a manifold, carries a 
canonical weak symplectic structure $\om_M\in\Om^2(T^*M)$, 
which is defined as follows. Let 
$\pi_M^*:T^*M\to M$ be the projection. Then the \idx{\it Liouville 
form} $\th_M\in \Om^1(T^*M)$ is given by 
$\th_M(X)=\langle \pi_{T^*M}(X),T(\pi_M^*)(X)\rangle $ for  
$X\in T(T^*M)$,
where $\langle \quad,\quad\rangle$ denotes the duality pairing 
$T^*M\x_M TM\to \mathbb R$. Then the symplectic structure on $T^*M$ is 
given by $\om_M = - d\th_M$, which of course in a local chart looks 
like 
$\om_E((v,v'),(w,w'))=\langle w',v\rangle_E-\langle v',w\rangle_E$.
The associated mapping $ \om\spcheck: T_{(0,0)}(E\x E')=E\x E' 
\to E'\x E''$ is given by $(v,v')\mapsto (-v',i_E(v))$,
where $i_E:E\to E''$ is the embedding into the bidual. 
So the canonical symplectic structure on $T^*M$ is strong if and only if 
all model spaces of the manifold $M$ are reflexive.

\subsection*{\nmb.{2.4}} Let $M$ be a weak symplectic manifold.
The first thing to note is that
the Hamiltonian mapping 
$\on{grad}^\om:C^\infty(M,\mathbb R) \to  \X(M,\om)$
does not make sense in general, since 
$ \om\spcheck: TM \to  T^*M$ is not invertible.
Namely, $\on{grad}^\om f = (\om\spcheck)\i\o df$ is defined only for 
those $f\in C^\infty(M,\mathbb R)$ with $df(x)$ in the image of 
$\om\spcheck$  for all 
$x\in M$. A similar difficulty arises for the
definition of the Poisson bracket on $C^\infty(M,\mathbb R)$.

\subsection*{Definition } For a weak symplectic manifold $(M,\om)$ 
let $T_x^\om M$ denote the real linear subspace 
$T_x^\om M= \om\spcheck_ x(T_x M) \subset T_x^*M = L(T_x M,\mathbb R)$,
and let us call it the {\it smooth cotangent space} with respect to 
the symplectic structure $\om$
of $M$ at $x$ in view of the embedding of test functions
into distributions. 
These vector spaces fit together to form a subbundle of $T^*M$ which 
is isomorphic to the tangent bundle $TM$ via 
$ \om\spcheck: TM\to T^\om M \subseteq T^*M$. It is in general not a 
splitting subbundle.

\subsection*{\nmb.{2.5}. Definition } 
For a weak symplectic vector space 
$(E,\om)$ let 
\begin{displaymath}
C^\infty_\om(E,\mathbb R) \subset C^\infty(E,\mathbb R)
\end{displaymath}
denote the linear subspace consisting of all smooth functions
$f:E \to  \mathbb R$ such that each iterated derivative 
$d^kf(x)\in L^k_{\text{sym}}(E;\mathbb R)$ has the
property that 
\begin{displaymath}
d^kf(x)(\quad,y_{2},\dots,y_{k}) \in {E}^\om
\end{displaymath}
is actually in the smooth dual 
$E^\om \subset {E}'$ for all $x,y_{2},\dots,y_{k} \in E$, and
that the mapping 
\begin{align*}
\prod^kE &\to  E
\\
(x,y_{2},\dots,y_{k})&\mapsto (\om\spcheck)\i(df(x)(\quad,y_{2},\dots,y_{k}))
\end{align*}is smooth. By the symmetry of higher derivatives, this is then 
true for all entries of $d^kf(x)$, for all $x$.  

\begin{proclaim}{\nmb.{2.6}. Lemma} For $f \in C^\infty(E,\mathbb R)$ the
following assertions are equivalent:
\begin{enumerate}
\item[(\nmb:{1})] $df:E \to  E{}'$ factors to a
        smooth mapping $E \to  E^\om$.
\item[(\nmb:{2})] $f$ has a smooth $\om$-gradient 
        $\on{grad}^\om f\in \X(E) = C^\infty(E,E)$ which 
        satisfies  
        $df(x)y = \om(\on{grad}^\om f(x),y)$.
\item[(\nmb:{3})] $f \in C^\infty_\om(E,\mathbb R)$.
\end{enumerate}
\end{proclaim}
\begin{demo}{Proof} Clearly, \nmb|{3} $\Rightarrow$ \nmb|{2}
$\Leftrightarrow$ \nmb|{1}. 
We have to show that
\nmb|{2} $\Rightarrow$ \nmb|{3}.
\newline
Suppose that $f:E \to  \mathbb R$ is smooth and
$df(x)y = \om(\on{grad}^\om f(x),y)$. Then 
\begin{align*}
d^kf(x)(y_{1},\dots,y_{k}) &= d^kf(x)(y_{2},\dots,y_{k},y_1)  
\\&
= (d^{k-1}(df))(x)(y_{2},\dots,y_{k})(y_1)  
\\&
= \om\bigl(d^{k-1}(\on{grad}^\om f)(x)(y_{2},\dots,y_{k}),y_1\bigr).\qed
\end{align*}\end{demo}

\subsection*{\nmb.{2.7}} For a weak symplectic manifold  
$(M,\om)$
let 
\begin{displaymath}
C^\infty_\om (M,\mathbb R) \subset C^\infty(M,\mathbb R)
\end{displaymath}
denote the linear subspace consisting of all smooth functions
$f:M \to  \mathbb R$ such that the differential $df:M\to T^*M$ factors 
to a smooth mapping $M\to T^\om M$. In view of lemma \nmb!{2.6} 
these are exactly those smooth functions on $M$ which admit a smooth 
$\om$-gradient $\on{grad}^\om f\in \X(M)$. Also the condition 
\nmb!{2.6.1} translates to a local differential condition
describing the functions in $C^\infty_\om(M,\mathbb R)$.

\begin{proclaim}{\nmb.{2.8}. Theorem} 
The Hamiltonian mapping 
$\on{grad}^\om :C^\infty_\om (M,\mathbb R) \to  \X(M,\om)$, which is given by
\begin{displaymath}
i_{\on{grad}^\om f}\om=df\quad\text{ or }\quad 
     \on{grad}^\om f:= (\om\spcheck)\i\o df
\end{displaymath}
is well defined. Also the Poisson bracket 
\begin{align*}
\{\quad,\quad\}&: C^\infty_\om (M,\mathbb R) \x
C^\infty_\om (M,\mathbb R) \to C^\infty_\om (M,\mathbb R)
\\
\{f,g\}&:= i_{\on{grad}^\om f}i_{\on{grad}^\om g}\om = 
     \om(\on{grad}^\om g,\on{grad}^\om f) = 
\\&
= dg(\on{grad}^\om f) = (\on{grad}^\om f)(g) 
\end{align*}is well defined and gives a Lie algebra structure to
the space $C^\infty_\om (M,\mathbb R)$, which also fulfills
\begin{displaymath}
\{f,gh\}=\{f,g\}h+g\{f,h\}.
\end{displaymath}
We have the following long exact sequence of Lie algebras
and Lie algebra homomorphisms:
\begin{displaymath}
0\to  H^0(M)\to C^\infty_\om (M,\mathbb R)  
     \East{{\on{grad}^\om}}{}\X(M,\om) \East{\ga}{}H^1_\om(M) \to 0,
\end{displaymath}
where $H^0(M)$ is the space of locally constant functions, and 
\begin{displaymath}
H^1_\om(M)=\frac{\{\ph\in C^\infty(M\gets T^\om M):d\ph=0\}}
     {\{df:f\in C^\infty_\om(M,\mathbb R)\}}
\end{displaymath}
is the first symplectic cohomology space of $(M,\om)$, a linear 
subspace of the De~Rham cohomology space $H^1(M)$.
\end{proclaim}

\begin{demo}{Proof} It is clear from lemma \nmb!{2.6}, that the
Hamiltonian mapping $\on{grad}^\om$ is well defined 
and has values in $\X(M,\om)$, 
since by \cite{23},~34.18.6 we have
\begin{displaymath}
\L_{\on{grad}^\om f}\om = i_{\on{grad}^\om f}d\om + 
     di_{\on{grad}^\om f}\om =ddf =0.
\end{displaymath}
By \cite{23},~34.18.7, the space $\X(M,\om)$ is a Lie 
subalgebra of $\X(M)$.
The Poisson
bracket is well defined as a mapping
$\{\quad,\quad\}: C^\infty_\om (M,\mathbb R) \x
C^\infty_\om (M,\mathbb R) \to  
C^\infty(M,\mathbb R) $;
it only remains to check that it has values in the subspace
$C^\infty_\om (M,\mathbb R)$.  

This is a local question, so we may assume that $M$ is an open subset 
of a con\-ven\-ient vector space equipped with a (non-constant) weak 
symplectic structure. 
So let $f$, $g\in C^\infty_\om (M,\mathbb R)$, then 
$\{f,g\}(x) = dg(x)(\on{grad}^\om f(x))$, and we have
\begin{align*} &d(\{f,g\})(x)y 
= d(dg(\quad)y)(x).\on{grad}^\om f(x) + dg(x)(d(\on{grad}^\om f)(x)y)
\\&
 = d(\om(\on{grad}^\om g(\quad),y)(x).\on{grad}^\om f(x)
     + \om\Bigl(\on{grad}^\om g(x),d(\on{grad}^\om f)(x)y\Bigr)
\\&
 = \om\Bigl(d(\on{grad}^\om g)(x)(\on{grad}^\om f(x)) 
    - d(\on{grad}^\om f)(x)(\on{grad}^\om g(x)),y\Bigr), 
\end{align*}since $\on{grad}^\om f\in\X(M,\om)$ and for any
$X\in\X(M,\om)$ the condition $\L_X\om=0$ implies
$\om(dX(x)y_1,y_2) = -\om(y_1,dX(x)y_2)$. So 
\nmb!{2.6.2} is satisfied, and thus 
$\{f,g\}\in C^\infty_\om (M,\mathbb R)$.

If $X\in\X(M,\om)$ then $di_X\om=\L_X\om = 0$, so 
$[i_X\om]\in H^1(M)$ is well defined, and by $i_X\om= \om\spcheck\ o X$ 
we even have $\ga(X):=[i_X\om]\in H^1_\om(M)$, so $\ga$ is well 
defined. 

Now we show that the sequence is exact.
Obviously, it is exact at $H^0(M)$ and at 
$C^\infty_\om(M,\mathbb R)$, since the kernel of $\on{grad}^\om$ consists of 
the locally constant functions. 
If $\ga(X)=0$ then $ \om\spcheck\o X =i_X\om=df$ 
for $f\in C^\infty_\om(M,\mathbb R)$, and clearly $X=\on{grad}^\om f$.
Now let us suppose that $\ph\in \Ga(T^\om M)\subset\Om^1(M)$ with 
$d\ph=0$. Then $X:=(\om\spcheck)\i\o\ph\in\X(M)$ is well defined and 
$\L_X\om=di_X\om = d\ph = 0$, so $X\in\X(M,\om)$ and $\ga(X)=[\ph]$.

Moreover, $H^1_\om(M)$ is a linear subspace of $H^1(M)$ since 
for $\ph\in \Ga(T^\om M)\subset\Om^1(M)$ with 
$\ph=df$ for $f\in C^\infty(M,\mathbb R)$ the vector field 
$X:=(\om\spcheck)\i\o\ph\in\X(M)$ is well defined, and since $ \om\spcheck\ o X 
=\ph=df$ by \nmb!{2.6.1} we have 
$f\in C^\infty_\om(M,\mathbb R)$ with $X=\on{grad}^\om f$. 

The mapping $\on{grad}^\om$ maps the Poisson bracket into the Lie bracket, 
since by \cite{23},~34.18 we have
\begin{align*}
i_{\on{grad}^\om \{f,g\}}\om &= d\{f,g\} 
= d\L_{\on{grad}^\om f}g = \L_{\on{grad}^\om f}dg =
\\&
= \L_{\on{grad}^\om f}i_{\on{grad}^\om g}\om - 
     i_{\on{grad}^\om g}\L_{\on{grad}^\om f}\om 
\\&
= [\L_{\on{grad}^\om f},i_{\on{grad}^\om g}]\om = 
     i_{[\on{grad}^\om f,\on{grad}^\om g]}\om.
\end{align*}Let us now check the properties of the Poisson bracket. By 
definition, it is skew symmetric, and we have
\begin{align*}
\{\{f,g\},h\} &= \L_{\on{grad}^\om \{f,g\}}h 
     = \L_{[\on{grad}^\om f,\on{grad}^\om g]}h 
     = [\L_{\on{grad}^\om f},\L_{\on{grad}^\om g}]h = 
\\&
= \L_{\on{grad}^\om f}\L_{\on{grad}^\om g}h - \L_{\on{grad}^\om g}\L_{\on{grad}^\om f}h 
=\{f,\{g,h\}\} - \{g,\{f,h\}\}
\\
\{f,gh\} &= \L_{\on{grad}^\om f}(gh) 
     = (\L_{\on{grad}^\om f}g)h + g\L_{\on{grad}^\om f}h =
\\&
= \{f,g\}h + g\{f,h\}.
\end{align*}Finally, it remains to show that all mappings in the sequence are Lie 
algebra homomorphisms, where we put the zero bracket on both 
cohomology spaces. For locally constant functions we have 
$\{c_1,c_2\}=\L_{\on{grad}^\om c_1}c_2 = 0$. We have already checked that 
$\on{grad}^\om$ is a Lie algebra homomorphism. For $X,Y\in\X(M,\om)$ 
\begin{displaymath}
i_{[X,Y]}\om = [\L_X,i_Y]\om = \L_Xi_Y\om + 0 = d i_Xi_Y\om + 
     i_X\L_Y\om = d i_Xi_Y\om
\end{displaymath}
is exact.
\qed \end{demo}

\subsection*{\nmb.{2.9}. Weakly symplectic group actions
} 
Let us suppose that an infinite dimensional regular Lie group $G$ with Lie
algebra $\g$ acts from the
right on a weak 
symplectic manifold $(M,\om)$ by 
$r:M\times G\to M$ in a way which respects $\om$, so that each 
transformation $r^g$ is a symplectomorphism. 
This is called a \idx{\it symplectic group action}. 
We shall use the notation $r(x,g)=r^g(x)=r_x(g)$.
Let us list some immediate consequences:

{\it \therosteritem{\nmb:{1}} The space $C^\infty_\om(M)^G$ of 
$G$-invariant smooth functions with $\om$-gradients 
is a Lie subalgebra for the Poisson 
bracket,} since for each $g\in G$ and $f,h\in C^\infty(M)^G$
we have $(r^g)^*\{f,h\} = \{(r^g)^*f,(r^g)^*h\}=\{f,h\}$. 

{\it\therosteritem{\nmb:{2}} 
For $x\in M$ the pullback of $\om$ to
the orbit $x.G$ is a 2-form, invariant under 
the action of $G$ on the orbit.}
In the finite dimensional case the orbit is an initial submanifold. In our
case this has to be checked directly in each example. In any case we have
something like a tangent bundle $T_x(x.G)=T(r_x)\g$.
If $i:x.G\to M$ is the embedding of the orbit then $r^g\o i = 
i\o r^g$, so that $i^*\om = i^*(r^g)^*\om= 
(r^g)^*i^*\om$ holds for each $g\in G$ and thus $i^*\om$ is 
invariant. 

\therosteritem{\nmb:{3}} 
The 
fundamental vector field mapping $\ze:\mathfrak g\to \mathfrak X(M,\om)$, 
given by $\ze_X(x)= T_e(r_x)X$ for $X\in \mathfrak g$ and 
$x\in M$, is a homomorphism of Lie algebras, where 
$\mathfrak g$ is the Lie algebra of $G$ (for a left action we get an 
anti homomorphism of Lie algebras). 
Moreover, $\ze$ takes values in $\X(M,\om)$. 
Let us consider again the exact sequence of Lie algebra homomorphisms 
from \nmb!{2.8}:
\begin{displaymath}
\cgaps{.4;.8;.9;.8;.4}\xymatrix{
0  \ar[r] & H^0(M) \ar[r]^{\al} & 
     C^\infty_\om(M) \ar[r]^{\on{grad}^\om} & \X(M,\om) \ar[r]^{\ga} & 
     H^1_\om(M) \ar[r] & 0
\\
 &  &  & \mathfrak g \ar@{.>}[lu]^{j} \ar[u]^{\ze} &  &
}\end{displaymath}
One can lift $\ze$ to a 
linear mapping $j:\mathfrak g\to C^\infty(M)$ if and only if 
$\gamma\circ\ze=0$. In this case the action of $G$ is called a 
\idx{\it Hamiltonian group action}, and the linear mapping 
$j:\mathfrak g\to C^\infty(M)$ is called a 
\idx{\it generalized Hamiltonian function} for the group action. 
It is unique up to addition of a mapping 
$\al\o\ta$ for $\ta:\mathfrak g\to H^0(M)$. 

{\it \therosteritem{\nmb:{4}} 
If $H^1_\om(M)=0$ then any symplectic 
action on $(M,\om)$ is a Hamiltonian action. 
But if $\ga\o\ze\ne 0$ we can replace $\mathfrak g$ by its Lie 
subalgebra $\on{ker}(\gamma\circ\ze)\subset \mathfrak g$ and consider the 
corresponding Lie subgroup $G$ which then admits a Hamiltonian 
action.}

{\it \therosteritem{\nmb:{5}} 
If the Lie algebra $\g$ is equal to its commutator subalgebra 
$[\g,\g]$, 
the linear span of all $[X,Y]$ for $X,Y\in \g$ 
(true for all full diffeomorphism groups), then any 
infinitesimal symplectic action $\ze:\g\to \X(M,\om)$ is a 
Hamiltonian action,} since then any $Z\in \g$ can be written as 
$Z=\sum_i[X_i,Y_i]$ so that 
$\ze_Z=\sum[\ze_{X_i},\ze_{Y_i}]\in \on{im}(\on{grad}^\om)$ since 
$\ga:\X(M,\om)\to H^1(M)$ is a homomorphism into the zero Lie 
bracket. 

\therosteritem{\nmb:{6}} If $j:\g\to (C^\infty_\om(M),\{\quad,\quad\})$
happens to be not a homomorphism of Lie algebras then 
$c(X,Y)=\{j(X),j(Y)\}-j([X,Y])$ lies in $H^0(M)$, and indeed 
$c:\g\x\g\to H^0(M)$ is a cocycle for the Lie algebra cohomology:
$c([X,Y],Z)+c([Y,Z],X)+c([Z,X],Y)=0$. If $c$ is a coboundary, i.e.,
$c(X,Y)=-b([X,Y])$, then $j+\al\o b$ is a Lie algebra homomorphism. 
If the cocycle $c$ is non-trivial we can use the central extension 
$H^0(M)\x_c \g$ with bracket $[(a,X),(b,Y)]=(c(X,Y),[X,Y])$ in the diagram
\begin{displaymath}
\cgaps{.4;.8;.9;.8;.4}\xymatrix{
0  \ar[r] & H^0(M) \ar[r]^{\al} & 
     C^\infty_\om(M) \ar[r]^{\on{grad}^\om} & \X(M,\om) \ar[r]^{\ga} & 
     H^1_\om(M) \ar[r] & 0
\\
 &  & H^1(M)\x_c\g \ar[r]^{\on{pr}_2} \ar[u]^{\bar\jmath} & 
                 \mathfrak g  \ar[u]^{\ze} &  &
}\end{displaymath}
where $\bar\jmath(a,X)=j(X)+\al(a)$. Then $\bar\jmath$ is a homomorphism of
Lie algebras.

\begin{proclaim}{\nmb.{2.10}. Momentum mapping} 
For an infinitesimal symplectic action, 
i.e.\ a homomorphism $\ze:\mathfrak g\to \X(M,\om)$ of Lie algebras, 
we can find a linear lift $j:\mathfrak g\to C^\infty_\om(M)$ if and only if 
there exists a mapping 
\begin{displaymath}
J\in C^\infty_\om(M,\g^*):=\{f\in C^\infty(M,\g^*):\langle f(\quad),X
\rangle\in C^\infty_\om(M)\text{  for all }X\in\g\}
\end{displaymath}
such that
\begin{displaymath}
\on{grad}^\om(\langle J,X \rangle)=\ze_X\quad\text{ for all }X\in \g.
\end{displaymath}
\end{proclaim}

The mapping $J\in C^\infty_\om(M,\g^*)$ is called the \idx{\it momentum mapping} 
for the infinitesimal action $\ze:\g\to \X(M,\om)$.
%This holds even for a Poisson manifold $(M,P)$ 
%and an infinitesimal action of a Lie algebra $\ze:\g\to \X(M,P)$ by 
%Poisson morphisms. 
Let us note again the relations between the generalized 
Hamiltonian $j$ and the momentum mapping $J$:
\begin{gather}
J:M\to \mathfrak g^*,\quad j:\g\to C^\infty_\om(M),\quad \ze:\g\to \X(M,\om)
\notag\\
\langle J, X\rangle = j(X)\in C^\infty_\om(M), 
\quad \on{grad}^\om(j(X)) =\ze(X), \quad 
X\in \mathfrak g,
\tag{\nmb:{1}}\\
i_{\ze(X)}\om = dj(X) = d\langle J,X \rangle,
\notag\end{gather}
where $\langle \; ,\; \rangle:\g^*\x \g\to \mathbb R$ is the duality pairing. 

\subsection*{\nmb.{2.11}. Basic properties of the momentum mapping }
Let $r:M\x G\to M$ be a Hamiltonian right action of an infinite dimensional
regular Lie group $G$
on a weak symplectic manifold $M$, let $j:\g\to C^\infty_\om(M)$ be a 
generalized Hamiltonian and let $J\in C^\infty_\om(M,\g^*)$ be the associated 
momentum mapping. 

{\it \therosteritem{\nmb:{1}} For $x\in M$, the transposed mapping of the 
linear mapping $dJ(x):T_xM\to \g^*$ is} 
\begin{displaymath}
dJ(x)^\top:\g\to T_x^*M, \qquad dJ(x)^\top = \check\om_x\o \ze,
\end{displaymath}
since for $\xi\in T_xM$ and $X\in \g$ we have 
\begin{displaymath}
\langle dJ(\xi),X \rangle = \langle i_{\xi}dJ, X \rangle
=  i_{\xi}d\langle J, X \rangle 
=  i_{\xi}i_{\ze_X}\om = \langle \check\om_x(\ze_X(x)),\xi \rangle.
\end{displaymath}

{\it \therosteritem{\nmb:{2}} The closure of the image $dJ(T_xM)$ of 
$dJ(x):T_xM\to \g*$ is the annihilator $\g_x^\o$ of the isotropy Lie 
algeba $\g_x:=\{X\in \g:\ze_X(x)=0\}$ in $\g^*$,} since the 
annihilator of the image is 
the kernel of the transposed mapping, 
\begin{displaymath}
\on{im}(dJ(x))^\o=\ker(dJ(x)^\top) = \ker(\check\om_x\o \ze) = 
\ker(\on{ev}_x\o\ze) = \g_x.
\end{displaymath}

{\it \therosteritem{\nmb:{3}}
The kernel of $dJ(x)$ is the symplectic orthogonal 
\begin{displaymath}
(T(r_x)\g)^{\bot,\om}=(T_x(x.G))^{\bot,\om}\subseteq T_xM,
\end{displaymath}
}
since for the annihilator of the kernel we have
\begin{multline*}
\ker(dJ(x))^\o = \overline{\on{im}(dJ(x)^\top)}
=\overline{\on{im}(\check\om_x\o\ze)} =
\\
= \overline{\{\check\om_x(\ze_X(x)):X\in \g\}} = 
     \overline{\check\om_x(T_x(x.G))}.
\end{multline*}
{\it \therosteritem{\nmb:{4}}
If $G$ is connected, $x\in M$ is a fixed point for the $G$-action if 
and only if $x$ is a critical point of $J$, i.e. $dJ(x)=0$.}

{\it \therosteritem{\nmb:{5}} (Emmy Noether's theorem)
Let $h\in C^\infty_\om(M)$ be a Hamiltonian function which is invariant 
under the Hamiltonian $G$ action. Then 
$dJ(\on{grad}^\om(h))=0$.
Thus the momentum mapping 
$J:M\to\g^*$ is constant on each trajectory (if it exists) of the Hamiltonian 
vector field $\on{grad}^\om(h)$.}
Namely,
\begin{multline*}
\langle dJ(\on{grad}^\om(h)),X \rangle =
d\langle J,X \rangle(\on{grad}^\om(h)) =
dj(X)(\on{grad}^\om(h)) =
\\
= \{h,j(X)\} = -dh(\on{grad}^\om j(X)) =
dh(\ze_X)=0.
\end{multline*}
E.\ Noether's theorem admits the following generalization.

\begin{proclaim}{\nmb.{2.12}. Theorem}
Let $G_1$ and $G_2$ be two regular Lie groups 
which act by Hamiltonian actions $r_1$ and $r_2$ on the weakly symplectic 
manifold $(M,\om)$, with momentum mappings $J_1$ and $J_2$, 
respectively. We assume that $J_2$ is $G_1$-invariant, i.e.\ $J_2$ is 
constant along all $G_1$-orbits, and that $G_2$ is connected.

Then $J_1$ is constant on the $G_2$-orbits and the two actions 
commute. 
\end{proclaim}

\begin{demo}{Proof} Let $\ze^i:\g_i\to \X(M,\om)$ be the two infinitesimal 
actions. Then for $X_1\in \g_1$ and $X_2\in \g_2$ we have
\begin{align*}
\L_{\ze^2_{X_2}}\langle J_1,X_1 \rangle &
     = i_{\ze^2_{X_2}}d\langle J_1,X_1 \rangle 
     = i_{\ze^2_{X_2}}i_{\ze^1_{X_1}}\om 
     = \{\langle J_2,X_2 \rangle,\langle J_1,X_1\rangle\} 
\\&
= -\{\langle J_1,X_1\rangle,\langle J_2,X_2 \rangle\}
     = -i_{\ze^1_{X_1}}d\langle J_2,X_2 \rangle
     = -\L_{\ze^1_{X_1}}\langle J_2,X_2 \rangle = 0
\end{align*}since $J_2$ is constant along each $G_1$-orbit. Since $G_2$ is 
assumed to be connected, $J_1$ is also 
constant along each $G_2$-orbit. We also saw that each Poisson bracket 
$\{\langle J_2,X_2 \rangle,\langle J_1,X_1\rangle\}$ vanishes; by 
$\on{grad}^\om{\langle J_i,X_i \rangle}=\ze^i_{X_i}$ we conclude that 
$[\ze^1_{X_1},\ze^2_{X_2}]=0$ for all $X_i\in \g_i$ which implies the 
result if also $G_1$ is connected. In the general case we can argue 
as follows: 
\begin{align*}
(r_1^{g_1}&)^* \ze^2_{X_2} =
     (r_1^{g_1})^* \on{grad}^\om{\langle J_2,X_2 \rangle} =
     (r_1^{g_1})^* (\check\om\i d{\langle J_2,X_2 \rangle}) 
\\&
=(((r_1^{g_1})^*\om)\spcheck) \i d{\langle (r_1^{g_1})^*J_2,X_2 \rangle} =
     (\check\om\i d{\langle J_2,X_2 \rangle} =
     \on{grad}^\om{\langle J_2,X_2 \rangle} = \ze^2_{X_2}.
\end{align*}Thus $r_1^{g_1}$ commutes with each $r_2^{\exp(t X_2)}$  and thus 
with each $r_2^{g_2}$, since $G_2$ is connected. 
\qed\end{demo}

\section*{\nmb0{3}. Right invariant weak Riemannian metrics on Lie groups  
} 
 
\subsection*{\nmb.{3.1}. Notation on Lie groups}
Let $G$ be a Lie group which may be infinite dimensional, but then is
supposed to be regular, with Lie
algebra $\g$. See appendix \nmb!{B} for more information. 
Let $\mu:G\x G\to G$ be the multiplication, let $\mu_x$ be left  
translation and $\mu^y$ be right translation, 
given by $\mu_x(y)=\mu^y(x)=xy=\mu(x,y)$. 

Let $L,R:\g\to \X(G)$ be the left 
and right invariant vector field mappings, given by 
$L_X(g)=T_e(\mu_g).X$ and $R_X=T_e(\mu^g).X$, respectively. 
They are related by $L_X(g)=R_{\Ad(g)X}(g)$.
Their flows are given by 
\begin{displaymath}
\on{Fl}^{L_X}_t(g)= g.\exp(tX)=\mu^{\exp(tX)}(g),\quad
\on{Fl}^{R_X}_t(g)= \exp(tX).g=\mu_{\exp(tX)}(g).\end{displaymath}

We also need the right  
Maurer-Cartan form $\ka=\ka^r\in\Om^1(G,\g)$, given by  
$\ka_x(\xi):=T_x(\mu^{x\i})\cdot \xi$. It satisfies the right
Maurer-Cartan equation $d\ka-\tfrac12[\ka,\ka]_\wedge=0$, where
$[\;,\;]_\wedge$ denotes the wedge product of $\g$-valued forms on
$G$ induced by the Lie bracket. Note that
$\tfrac12[\ka,\ka]_\wedge (\xi,\et) = [\ka(\xi),\ka(\et)]$.
The (exterior) derivative of the function $\Ad:G\to GL(\g)$ can be 
expressed by
\begin{displaymath}
d\Ad = \Ad.(\ad\o\ka^l) = (\ad\o \ka^r).\Ad,\end{displaymath}
since we have
$d\Ad(T\mu_g.X) = \frac d{dt}|_0 \Ad(g.\exp(tX))
= \Ad(g).\ad(\ka^l(T\mu_g.X))$.

\subsection*{\nmb.{3.2}. Geodesics of a right invariant metric on a Lie  
group } 
 
Let $\ga=\langle \;,\;\rangle:\g\x\g\to\mathbb R$ be a positive  
definite bounded (weak) inner product. Then  
\begin{equation}
\ga_x(\xi,\et)=\langle T(\mu^{x\i})\cdot\xi,
     T(\mu^{x\i})\cdot\et)\rangle =  
     \langle \ka(\xi),\ka(\et)\rangle  
\tag{\nmb:{1}}
\end{equation}
is a right invariant (weak) Riemannian metric on $G$, and any
(weak) right invariant bounded Riemannian metric is of this form, for
suitable  $\langle \;,\;\rangle$. 
 
Let $g:[a,b]\to G$ be a smooth curve.  
The velocity field of $g$, viewed in the right trivializations, 
coincides  with the right logarithmic derivative
\begin{displaymath}
\de^r(g)=T(\mu^{g\i})\cdot \partial_t g =  
\ka(\partial_t g) = (g^*\ka)(\partial_t),
\text{ where }
\partial_t=\frac{\partial}{\partial t}. 
\end{displaymath}
The energy of the curve $g(t)$ is given by  
\begin{displaymath}
E(g) = \tfrac12\int_a^bG_g(g',g')dt = \tfrac12\int_a^b 
     \langle (g^*\ka)(\partial_t),(g^*\ka)(\partial_t)\rangle\, dt. 
\end{displaymath}
For a variation $g(s,t)$ with fixed endpoints we have then, using the  
right Maurer-Cartan equation and integration by parts, 
\begin{align*} 
&\partial_sE(g) = \tfrac12\int_a^b2 
     \langle \partial_s(g^*\ka)(\partial_t),\,
                            (g^*\ka)(\partial_t)\rangle\, dt
\\&
= \int_a^b \langle \partial_t(g^*\ka)(\partial_s) - 
         d(g^*\ka)(\partial_t,\partial_s),\,
                (g^*\ka)(\partial_t)\rangle\,dt
\\&
= \int_a^b \left(-\langle (g^*\ka)(\partial_s),\, 
     \partial_t(g^*\ka)(\partial_t)\rangle - 
     \langle [(g^*\ka)(\partial_t),(g^*\ka)(\partial_s)],\, 
     (g^*\ka)(\partial_t)\rangle\right)\, dt
\\&
= -\int_a^b \langle (g^*\ka)(\partial_s),\, 
     \partial_t(g^*\ka)(\partial_t) + 
     \ad((g^*\ka)(\partial_t))^{\top}
           ((g^*\ka)(\partial_t))\rangle\, dt
\\
\end{align*}
where $\ad((g^*\ka)(\partial_t))^{\top}:\g\to\g$ is the adjoint of  
$\ad((g^*\ka)(\partial_t))$ with respect to the inner product  
$\langle \;,\; \rangle$. In infinite dimensions one also has to  
check the existence of this adjoint. 
In terms of the right logarithmic derivative $u:[a,b]\to \g$ of  
$g:[a,b]\to G$, given by  
$u(t):= g^*\ka(\partial_t) = T_{g(t)}(\mu^{g(t)\i})\cdot g'(t)$, 
the geodesic equation has the expression: 
\begin{equation}
\boxed{\quad
u_t = - \ad(u)^{\top}u\quad} \tag{\nmb:{2}}
\end{equation}
This is, of course, just the Euler-Poincar\'e equation for right
invariant systems using the Lagrangian given by the kinetic energy
(see \cite{25}, section 13). 
 
\subsection*{\nmb.{3.3}. The covariant derivative } 
Our next aim is to derive the Riemannian curvature and for that we  
develop the basis-free version of Cartan's method of moving frames 
in this setting, which also works in infinite dimensions. 
The right trivialization, or framing, $(\pi_G,\ka):TG\to G\x\g$  
induces the isomorphism $R:C^\infty(G,\g)\to \X(G)$, given by  
$R(X)(x):= R_X(x):=T_e(\mu^x)\cdot X(x)$, for $X\in C^\infty(G,\g)$ and
$x\in G$. Here $\X(G):=\Ga(TG)$ denote the Lie algebra of all vector 
fields. For the Lie bracket and the Riemannian metric we have 
\begin{gather} 
[R_X,R_Y] = R(-[X,Y]_\g + dY\cdot R_X - dX\cdot R_Y),
\tag{\nmb:{1}}\\
R\i[R_X,R_Y] = -[X,Y]_\g + R_X(Y) - R_Y(X),
\notag\\
\ga_x(R_X(x),R_Y(x)) = \ga( X(x),Y(x))\,,\, x\in G. 
\notag\end{gather}
In the sequel we shall compute in $C^\infty(G,\g)$ instead of 
$\X(G)$. In particular, we shall use the convention
\begin{displaymath}
\nabla_XY := R\i(\nabla_{R_X}R_Y)\quad\text{ for }X,Y\in C^\infty(G,\g).
\end{displaymath}
to express the Levi-Civita covariant derivative. 
 
\begin{proclaim}{Lemma} 
Assume that for all $\xi\in\g$ the adjoint $\ad(\xi)^\top$  
with respect to the inner  
product $\langle \;,\;\rangle$ exists and that  
$\xi\mapsto \ad(\xi)^\top$ is bounded. 
Then the Levi-Civita covariant derivative of the metric  
\nmb!{3.2.1} exists and is given for any $X,Y \in C^\infty(G,\g)$  in
terms of the isomorphism $R$ by
\begin{equation}
\nabla_XY= dY.R_X + \tfrac12\ad(X)^\top Y 
     +\tfrac12\ad(Y)^\top X - \tfrac12\ad(X)Y.
\tag{\nmb:{2}}\end{equation}
\end{proclaim} 
 
\begin{demo}{Proof}  
Easy computations show that this formula satisfies the axioms of a
covariant derivative, that relative to it the Riemannian
metric is covariantly constant, since 
\begin{displaymath}
R_X\ga( Y,Z) = \ga( dY.R_X,Z) + \ga( Y,dZ.R_X)
     = \ga( \nabla_XY,Z) + \ga( Y,\nabla_XZ),
\end{displaymath}
and that it is torsion free, since 
\begin{displaymath}
\nabla_XY-\nabla_YX + [X,Y]_\g - dY.R_X + dX.R_Y = 0.\qed
\end{displaymath}
\end{demo}
\smallskip 
 
For $\xi \in \g$ define $\al(\xi):\g\to\g$ by
$\al(\xi)\eta:=\ad(\eta)^\top \xi$. With this notation, the previous
lemma states that for all $X\in C^\infty(G,\g)$ the covariant
derivative of the Levi-Civita connection has the expression
\begin{equation}
\nabla_X = R_X + \tfrac12\ad(X)^\top +\tfrac12\al(X) - 
\tfrac12\ad(X).
\tag{\nmb:{3}}\end{equation}

\subsection*{\nmb.{3.4}. The curvature } 
First note that we have the following relations: 
\begin{alignat}2 
[R_X,\ad(Y)] &= \ad(R_X(Y)),&\quad  
[R_X,\al(Y)] &= \al(R_X(Y)),
\tag{\nmb:{1}}\\
[R_X,\ad(Y)^\top] &= \ad(R_X(Y))^\top,&\quad 
[\ad(X)^\top,\ad(Y)^\top] &= -\ad([X,Y]_\g)^\top.  
\notag\end{alignat}
The Riemannian curvature is then computed by  
\begin{align*} 
&\mathcal{R}(X,Y) =  
[\nabla_X,\nabla_Y]-\nabla_{-[X,Y]_\g+R_X(Y)-R_Y(X)}
\\&
= [R_X+\tfrac12\ad(X)^\top+\tfrac12\al(X)-\tfrac12\ad(X), 
     R_Y+\tfrac12\ad(Y)^\top+\tfrac12\al(Y)-\tfrac12\ad(Y)]
\\&\quad
- R_{-[X,Y]_\g + R_X(Y) - R_Y(X)} 
     -\tfrac12\ad(-[X,Y]_\g + R_X(Y) - R_Y(X))^\top
\\&\quad
-\tfrac12\al(-[X,Y]_\g + R_X(Y) - R_Y(X)) 
     +\tfrac12\ad(-[X,Y]_\g + R_X(Y) - R_Y(X)) 
\\&
= -\tfrac14[\ad(X)^\top+\ad(X),\ad(Y)^\top+\ad(Y)]
\tag{\nmb:{2}}\\&
\quad+\tfrac14[\ad(X)^\top-\ad(X),\al(Y)]  
      +\tfrac14[\al(X),\ad(Y)^\top-\ad(Y)] 
\\&
\quad+\tfrac14[\al(X),\al(Y)] +\tfrac12\al([X,Y]_\g).
\end{align*}
If we plug in all definitions and use 4 times 
the Jacobi identity we get the following expression 
\begin{align} 
&\ga( 4\mathcal{R}(X,Y)Z,U) = 
+ 2\ga( [X,Y],[Z,U] )                   
- \ga( [Y,Z],[X,U] )                   
+ \ga( [X,Z],[Y,U] )
\notag\\&
- \ga( Z,[U,[X,Y]] )  
+ \ga( U,[Z,[X,Y]] )  
- \ga( Y,[X,[U,Z]] )                   
- \ga( X,[Y,[Z,U]] )
\notag\\&
+ \ga( \ad(X)^\top Z,\ad(Y)^\top U )   
+ \ga( \ad(X)^\top Z,\ad(U)^\top Y )   
\notag\\&
+ \ga( \ad(Z)^\top X,\ad(Y)^\top U )   
- \ga( \ad(U)^\top X,\ad(Y)^\top Z ) 
\tag{\nmb:{3}}\\&
- \ga( \ad(Y)^\top Z,\ad(X)^\top U ) 
- \ga( \ad(Z)^\top Y,\ad(X)^\top U ) 
\notag\\&
- \ga( \ad(U)^\top X,\ad(Z)^\top Y )   
+ \ga( \ad(U)^\top Y,\ad(Z)^\top X ).
\notag\end{align}
This yields the following expression which is useful for computing the
sectional curvature:
\begin{align}
4&\ga(\mathcal{R}(X,Y)X,Y) =  3\ga(\on{ad}(X)Y,\on{ad}(X)Y)                   
- 2\ga(\on{ad}(Y)^\top X,\on{ad}(X)Y)  
\notag\\&\quad
- 2\ga( \on{ad}(X)^\top Y,\on{ad}(Y)X )                   
+ 4\ga( \on{ad}(X)^\top  X,\on{ad}(Y)^\top  Y )   
\tag{\nmb:{4}}\\&\quad
- \ga( \on{ad}(X)^\top  Y+\on{ad}(Y)^\top  X,\on{ad}(X)^\top 
Y+\on{ad}(Y)^\top X ).
\notag\end{align}
\subsection*{\nmb.{3.5}. Jacobi fields, I } 
We compute first the Jacobi equation directly via variations of 
geodesics.  
So let $g:\mathbb R^2\to G$ be smooth, $t\mapsto g(t,s)$ a geodesic for  
each $s$. Let again $u=\ka(\partial_t g) = (g^*\ka)(\partial_t)$ be  
the velocity field along the geodesic in right trivialization which  
satisfies the geodesic equation  
$u_t=-\ad(u)^\top u$. 
Then $y:= \ka(\partial_s g)=(g^*\ka)(\partial_s)$ is the Jacobi field  
corresponding to this variation, written in the right trivialization. 
 From the right Maurer-Cartan equation we then have: 
\begin{align*} 
y_t &= \partial_t (g^*\ka)(\partial_s)  
     = d(g^*\ka)(\partial_t,\partial_s)  
     + \partial_s (g^*\ka)(\partial_t) + 0
\\&
= [(g^*\ka)(\partial_t),(g^*\ka)(\partial_s)]_{\mathfrak g} + u_s
\\&
= [u,y] +u_s. 
\end{align*}Using the geodesic equation, the definition of $\al$, and
the fourth relation in \nmb!{3.4.1}, this identity implies
\begin{align*} 
u_{st} &= u_{ts} = \partial_s u_t = -\partial_s (\ad(u)^\top u)  
     = -\ad(u_s)^\top u - \ad(u)^\top u_s
\\&
= -\ad(y_t+[y,u])^\top u - \ad(u)^\top (y_t+[y,u])
\\&
= -\al(u)y_t -\ad([y,u])^\top u  
     - \ad(u)^\top y_t - \ad(u)^\top ([y,u])
\\&
= - \ad(u)^\top y_t -\al(u)y_t +[\ad(y)^\top,\ad(u)^\top]u  
     - \ad(u)^\top\ad(y)u\,.
\\
\end{align*}
Finally we get the Jacobi equation as 
\begin{align} 
y_{tt} &= [u_t,y] + [u,y_t] + u_{st}
\notag\\&
= \ad(y)\ad(u)^\top u + \ad(u)y_t - \ad(u)^\top y_t 
\notag\\&
\qquad-\al(u)y_t +[\ad(y)^\top,\ad(u)^\top]u - 
\ad(u)^\top\ad(y)u\,,
\notag\\
y_{tt}&= [\ad(y)^\top+\ad(y),\ad(u)^\top]u   
     - \ad(u)^\top y_t -\al(u)y_t + \ad(u)y_t\,.
\tag{\nmb:{1}}\end{align}
\subsection*{\nmb.{3.6}. Jacobi fields, II } 
Let $y$ be a Jacobi field along a geodesic 
$g$ with right trivialized  
velocity field $u$. Then $y$ should satisfy the analogue of the 
finite dimensional Jacobi equation 
\begin{displaymath}
\nabla_{\partial_t}\nabla_{\partial_t}y + \mathcal{R}(y,u)u = 0 
\end{displaymath}
We want to show that this leads to same equation 
as \nmb!{3.5.1}. 
First note that from \nmb!{3.3.2} we have 
\begin{displaymath}
\nabla_{\partial_t}y  
     = y_t +\tfrac12\ad(u)^\top y + \tfrac12\al(u)y - 
\tfrac12\ad(u)y 
\end{displaymath}
so that, using $u_t=-\ad(u)^\top u$, we get: 
\begin{align*} 
\nabla_{\partial_t}\nabla_{\partial_t}y  
&= \nabla_{\partial_t}\Bigl(y_t +\tfrac12\ad(u)^\top y + 
\tfrac12\al(u)y -  
\tfrac12\ad(u)y\Bigr)
\\
%\allowdisplaybreak
&= y_{tt} +\tfrac12\ad(u_t)^\top y +\tfrac12\ad(u)^\top y_t  
     + \tfrac12\al(u_t)y 
\\&
\quad+ \tfrac12\al(u)y_t - \tfrac12\ad(u_t)y  
     -\tfrac12\ad(u)y_t
\\&
\quad+\tfrac12\ad(u)^\top\Bigl(y_t+\tfrac12\ad(u)^\top y  
     + \tfrac12\al(u)y - \tfrac12\ad(u)y\Bigr) 
\\&
\quad+\tfrac12\al(u)\Bigl(y_t+\tfrac12\ad(u)^\top y 
     + \tfrac12\al(u)y - \tfrac12\ad(u)y\Bigr) 
\\&
\quad-\tfrac12\ad(u)\Bigl(y_t+\tfrac12\ad(u)^\top y 
     + \tfrac12\al(u)y - \tfrac12\ad(u)y\Bigr) 
\\
%\allowdisplaybreak
&= y_{tt}  + \ad(u)^\top y_t + \al(u)y_t - \ad(u)y_t
\\&
\quad-\tfrac12\al(y)\ad(u)^\top u  
      - \tfrac12\ad(y)^\top\ad(u)^\top u  
      - \tfrac12\ad(y)\ad(u)^\top u 
\\&
\quad+\tfrac12\ad(u)^\top\Bigl(\tfrac12\al(y)u  
     + \tfrac12\ad(y)^\top u + \tfrac12\ad(y)u\Bigr) 
\\&
\quad+\tfrac12\al(u)\Bigl(\tfrac12\al(y)u 
     + \tfrac12\ad(y)^\top u + \tfrac12\ad(y)u\Bigr) 
\\&
\quad-\tfrac12\ad(u)\Bigl(\tfrac12\al(y)u 
     + \tfrac12\ad(y)^\top u + \tfrac12\ad(y)u\Bigr)\,. 
\\
\end{align*}In the second line of the last expression we use  
\begin{displaymath}
-\tfrac12\al(y)\ad(u)^\top u= -\tfrac14\al(y)\ad(u)^\top u  
     -\tfrac14\al(y)\al(u)u 
\end{displaymath}
and similar forms for the other two terms to get: 
\begin{align*} 
\nabla_{\partial_t}\nabla_{\partial_t}y  
&= y_{tt}  + \ad(u)^\top y_t + \al(u)y_t - \ad(u)y_t
\\&
\quad+\tfrac14[\ad(u)^\top,\al(y)]u  
     + \tfrac14[\ad(u)^\top,\ad(y)^\top]u  
     + \tfrac14[\ad(u)^\top,\ad(y)]u 
\\&
\quad+\tfrac14[\al(u),\al(y)]u 
     + \tfrac14[\al(u),\ad(y)^\top]u  
     + \tfrac14[\al(u),\ad(y)]u 
\\&
\quad-\tfrac14[\ad(u),\al(y)]u 
     - \tfrac14[\ad(u),\ad(y)^\top + \ad(y)]u, 
\\
\end{align*}where in the last line we also used $\ad(u)u=0$. 
We now compute the curvature term using \nmb!{3.4.2}: 
\begin{align*} 
\mathcal{R}(y,u)u 
&= -\tfrac14[\ad(y)^\top+\ad(y),\ad(u)^\top+\ad(u)]u
\\&
\quad+\tfrac14[\ad(y)^\top-\ad(y),\al(u)]u  
      +\tfrac14[\al(y),\ad(u)^\top-\ad(u)]u 
\\&
\quad+\tfrac14[\al(y),\al(u)] +\tfrac12\al([y,u])u
\\&
= -\tfrac14[\ad(y)^\top+\ad(y),\ad(u)^\top]u 
      -\tfrac14[\ad(y)^\top+\ad(y),\ad(u)]u
\\&
\quad+\tfrac14[\ad(y)^\top,\al(u)]u  
      -\tfrac14[\ad(y),\al(u)]u  
      +\tfrac14[\al(y),\ad(u)^\top-\ad(u)]u 
\\&
\quad+\tfrac14[\al(y),\al(u)]u +\tfrac12\ad(u)^\top\ad(y)u\,.
\\
\end{align*}Summing up we get 
\begin{align*} 
\nabla_{\partial_t}\nabla_{\partial_t}y + \mathcal{R}(y,u)u  
&= y_{tt}  + \ad(u)^\top y_t + \al(u)y_t - \ad(u)y_t 
\\&
\quad-\tfrac12[\ad(y)^\top+\ad(y),\ad(u)^\top]u 
\\&
\quad+ \tfrac12[\al(u),\ad(y)]u +\tfrac12\ad(u)^\top\ad(y)u\,. 
\end{align*}Finally we need the following computation using \nmb!{3.4.1}: 
\begin{align*} 
\tfrac12[\al(u),\ad(y)]u  
&= \tfrac12\al(u)[y,u] - \tfrac12\ad(y)\al(u)u
\\&
= \tfrac12\ad([y,u])^\top u - \tfrac12\ad(y)\ad(u)^\top u
\\&
= -\tfrac12[\ad(y)^\top,\ad(u)^\top] u -
       \tfrac12\ad(y)\ad(u)^\top u\,. 
\end{align*}Inserting we get the desired result: 
\begin{align*} 
\nabla_{\partial_t}\nabla_{\partial_t}y + \mathcal{R}(y,u)u  
&= y_{tt}  + \ad(u)^\top y_t + \al(u)y_t - \ad(u)y_t 
\\&
\quad-[\ad(y)^\top+\ad(y),\ad(u)^\top]u. 
\end{align*}
\subsection*{\nmb.{3.7}. The weak symplectic structure on 
the space of Jacobi  
fields } 
Let us assume now that the geodesic equation in $\g$ 
\begin{displaymath}
u_t=-\ad(u)^\top u  
\end{displaymath}
admits a unique solution for some time 
interval, depending smoothly on the  
choice of the initial value $u(0)$. 
Furthermore we assume that $G$ is  
a regular Lie group \nmb!{B.9} 
so that each smooth curve $u$ in $\g$  
is the right logarithmic derivative of a smooth curve  
$g$ in $G$ which depends smoothly on $u$, 
so that $u=(g^*\ka)(\partial_t)$.
Furthermore we have to assume 
that the Jacobi equation along $u$ 
admits a unique solution for some time,
depending smoothly on the initial values $y(0)$ and $y_t(0)$.
These are non-trivial assumptions: 
in \nmb!{A.4} there are examples of
ordinary linear differential equations 
`with constant coefficients' which
violate existence or uniqueness. 
These assumptions have to be 
checked in the special situations.
Then the space $\mathcal{J}_u$ of all Jacobi fields 
along the geodesic $g$ described by $u$
is isomorphic to the space $\g\x \g$ of all initial data.
 
There is the well known symplectic structure on the space $\mathcal{J}_u$  
of all Jacobi fields along a fixed geodesic with velocity field $u$, see
e.g.\ \cite{20}, II, p.70.  
It is given by the following expression which is constant in time $t$: 
\begin{align*} 
\om(y,z) :&= \langle y,\nabla_{\partial_t}z \rangle  
     - \langle \nabla_{\partial_t}y,z \rangle 
\\&
= \langle y, z_t + \tfrac12\ad(u)^\top z 
          + \tfrac12\al(u)z - \tfrac12\ad(u)z \rangle 
\\&
\qquad-\langle y_t + \tfrac12\ad(u)^\top y  
          + \tfrac12\al(u)y - \tfrac12\ad(u)y, z  \rangle
\\&
= \langle y,z_t \rangle - \langle y_t, z \rangle  
     + \langle [u,y],z \rangle - \langle y,[u,z] \rangle  
     - \langle [y,z],u \rangle 
\\&
= \langle y,z_t-\ad(u)z+\tfrac12\al(u)z \rangle  
     -\langle y_t-\ad(u)y+\tfrac12\al(u)y, z \rangle.  
\end{align*}It is worth while to check directly from the Jacobi field equation 
\nmb!{3.5.1} that $\om(y,z)$ is indeed constant in $t$. 
Clearly $\om$ is a weak symplectic structure on the relevant
vector space $\mathcal{J}_u\cong \g\x \g$, i.e., $\om$ gives an
injective (but in general not surjective) linear
mapping $\mathcal{J}_u\to \mathcal{J}_u^*$. This is seen most
easily by writing 
\begin{displaymath}
\om(y,z)=\langle y, z_t - \Ga_g(u,z)\rangle|_{t=0} 
        - \langle y_t-\Ga_g(u,y),z\rangle|_{t=0}
\end{displaymath}
which is induced from the standard symplectic structure on 
$\g\x \g^*$ by applying first the automorphism 
$(a,b)\mapsto (a,b-\Ga_g(u,a))$  to $\g\x\g$ and then by injecting 
the second factor $\g$ into its dual $\g^*$.

For regular (infinite dimensional) Lie groups variations of geodesics
exist, but there is no general theorem stating that they are uniquely
determined by $y(0)$ and $y_t(0)$. For concrete regular Lie groups, this
needs to be shown directly.   

\section*{\nmb0{4}. The Hamiltonian approach}

\subsection*{\nmb.{4.1}. The symplectic form on $T^*G$ and $G\x\g^*$ }
For an (infinite dimensional regular) Lie group $G$ with Lie algebra $\g$,
elements in the cotangent bundle $\pi:(T^*G,\om_{G})\to G$ are said to be in 
{\it material} or {\it Lagrangian representation}.
The cotangent bundle $T^*G$ has two trivializations, the left one 
\begin{align*}
&(\pi_G,\ka^l):T^*G\to G\x \g^*,
\\&
T^*_gG\ni \al_g\mapsto (g,T_e(\mu_g)^*\al_g=T_g^*(\mu_{g\i})\al_g),
%\\&
%T_g(\mu_{g\i})*\al \gets (g,\al)\in G\x\g^*
\end{align*}
also called the {\it body coordinate chart}, and the right one, 
\begin{align}
&(\pi_G,\ka^r):T^*G\to G\x \g^*,
\notag\\& 
T^*G\ni \al_g\mapsto (g,T_e(\mu^g)^*\al_g=T^*_g(\mu^g)\al_g),
\tag{\nmb:{1}}\\&
T_g(\mu^{g\i})^*\al \gets (g,\al) \in G\x\g^*
\notag\end{align}
also called the {\it space} or {\it Eulerian coordinate chart}. We will use
only this from now on.
The canonical 1-form in the Eulerian chart is given by (where
$\langle\quad,\quad\rangle:\g^*\x\g\to\mathbb R$ is the duality pairing):
\begin{align}
\th_{G\x\g^*}&(\xi_g,\al,\be) 
:=(((\pi,\ka^r)\i)^*\th_{G})_{(g,\al)}(\xi_g,\al,\be)
\notag\\&
=\th_G(T_{(g,\al)}(\pi,\ka^r)\i(\xi_g,\al,\be))                                                 \notag\\&
= \Bigl\langle  \pi_{T^*G}(T_{(g,\al)}(\pi,\ka^r)\i(\xi_g,\al,\be)),
  T(\pi)(T_{(g,\al)}(\pi,\ka^r)\i(\xi_g,\al,\be))\Bigr\rangle
\notag\\&
= \Bigl\langle  (\pi,\ka^r)\i(\pi_G,\pi_{\g^*})(\xi_g,\al,\be),
  T(\pi\o (\pi,\ka^r)\i)(\xi_g,\al,\be))\Bigr\rangle
\notag\\&
= \Bigl\langle  (\pi,\ka^r)\i(g,\al),
  T(\on{pr}_1)(\xi_g,\al,\be))\Bigr\rangle
= \Bigl\langle T_g(\mu^{g\i})^*\al,\xi_g        \Bigr\rangle
\notag\\&
= \langle \al, T_g(\mu^{g\i})\xi_g \rangle
=\langle \al,\ka^r(\xi_g) \rangle
\tag{\nmb:{2}}\end{align}
Now it is easy to to take the exterior derivative: For $X_i\in G$, thus 
$R_{X_i}\in\X(G)$ right
invariant vector fields, and $\g^*\ni \be_i\in \X(\g^*)$ constant vector
fields, we have
\begin{align}
\th_{G\x\g^*}&(R_{X_i}(g),(\al,\be_i)) = \langle \al,X_i \rangle
\notag\\
\th_{G\x\g^*}&(R_{X_i},\be_i) = \langle \on{Id}_{\g^*},X_i \rangle
= \langle \quad,X_i \rangle
\notag\\
\om_{G\x\g^*}&((R_{X_1},\be_1),(R_{X_2},\be_2))
=-d\th_{G\x\g^*}((R_{X_1},\be_1),(R_{X_2},\be_2))
\notag\\&
= -(R_{X_1},\be_1)(\th_{G\x\g^*}(R_{X_2},\be_2))
+(R_{X_2},\be_2)(\th_{G\x\g^*}(R_{X_1},\be_1))
\notag\\&\quad
+ (\th_{G\x\g^*}([(R_{X_1},\be_1),R_{X_2},\be_2)])
\notag\\&
=- (R_{X_1},\be_1)(\langle\quad,X_2\rangle)
+ (R_{X_2},\be_2)(\langle\quad,X_1\rangle)
\notag\\&\quad
+ (\th_{G\x\g^*}(-R_{[X_1,X_2]},0_{\g^*})
\notag\\&
=-\langle\be_1,X_2\rangle
+ \langle\be_2,X_1\rangle
-\langle\quad,[X_1,X_2]\rangle
\notag\\
(\om_{G\x\g^*}&)_{(g,\al)}((T(\mu^g).X_1,\be_1),(T(\mu^g)X_2,\be_2))
\notag\\&
=\langle\be_2,X_1\rangle
-\langle\be_1,X_2\rangle
-\langle\al,[X_1,X_2]\rangle
\tag{\nmb:{3}}\end{align}

\subsection*{\nmb.{4.2}. The symplectic form on $TG$ and $G\x\g$ and the
momentum mapping }
We consider an (infinite dimensional regular) Lie group $G$ with Lie algebra $\g$
and a bounded weak inner product $\ga:\g\x\g\to \mathbb R$ with the property
the transpose of the adjoint action of $G$ on $\g$, 
\begin{displaymath}
\ga(\Ad(g)^\top X,Y)=\ga(X,\Ad(g)X),
\end{displaymath}
exists. It is then unique and a right action of $G$ on $\g$.
By differentiating it follows that then also 
the transpose of the adjoint operation of $\g$ exists:
\begin{displaymath}
\ga(\ad(X)^\top Y,Z)= \p_t|_0 \ga(\Ad(\exp(tX))^\top Y,Z) =
\ga(Y,\ad(X)Z)
\end{displaymath}
exists.

We exted $\ga$ to
a right invariant Riemannian metric, again called $\ga$ on $G$ and consider
$\ga:TG\to T^*G$. Then we pull back the canonical symplectic structure
$\om_{G}$ to $G\x \g$ in the right or Eulerian trivialization:
\begin{align}
\ga&:G\x\g \to G\x\g^*, (g,X)\mapsto (g,\ga(X))
\notag\\
(\ga^*\om)_{(g,X)}&((T(\mu^g).X_1,X,Y_1),(T(\mu^g)X_2,X,Y_2))
\notag\\&=
\om_{(g,\ga(X))}((T(\mu^g).X_1,\ga(X),\ga(Y_1)),(T(\mu^g)X_2,\ga(X),\ga(Y_2)))
\notag\\&
=\langle\ga(Y_2),X_1\rangle
-\langle\ga(Y_1),X_2\rangle
-\langle\ga(X),[X_1,X_2]\rangle
\notag\\&
=\ga(Y_2,X_1)
-\ga(Y_1,X_2)
-\ga(X,[X_1,X_2])
\tag{\nmb:{1}}\end{align}
Since $\ga$ is a weak inner product, 
$\ga^*\om$ is again a weak symplectic structure on 
$TG\cong G\x \g$. We compute the Hamiltonian vector field mapping (symplectic
gradient) for functions $f\in C^\infty_{\ga^*\om}(G\x\g)$ admitting such
gradients:
\begin{align*}
(\ga^*\om&)_{(g,X)}
  \left(\on{grad}^{\ga^*\om}(f)(g,X),(T(\mu^g)X_2,X,Y_2)\right) 
= df(T(\mu^g)X_2;X,Y_2) 
\\&
= d_1f(g,X)(T(\mu^g)X_2) + d_2f(g,X)(Y_2)
\\&
= \ga(\ka^r(\on{grad}_1^\ga(f)(g,X)),X_2) +  
   \ga(\on{grad}_2^\ga(f)(g,X),Y_2)   
\\&
=\ga(X_1,Y_2) +\ga(-Y_1-\ad(X_1)^\top X,X_2)
\quad\text{  by \thetag{\nmb|{1}}.} 
\end{align*}Thus the Hamiltonian vector field of 
$f\in C^\infty_{\ga^*\om}(G\x\g)=C^\infty_\ga(G\x \g)$ is 
\begin{align*}
&\on{grad}^{\ga^*\om}(f)(g,X) =
\tag{\nmb:{2}}\\
&\big(T(\mu^g)\on{grad}_2^\ga(f)(g,X),X,
  -\ad(\on{grad}_2^\ga(f)(g,X))^\top X-\ka^r(\on{grad}^\ga_1(f)(g,X))\big)
\end{align*}
In particular, the Hamiltonian vector field of the function 
$(g,X)\mapsto \ga(X,X)=\|X\|_\ga^2$ on $TG$ is given by: 
\begin{displaymath}
\on{grad}^{\ga^*\om}(\tfrac12\|\quad\|_\ga^2)(g,X) =
(T(\mu^g)X;X,-\ad(X)^\top X)
\tag{\nmb:{3}}
\end{displaymath}
We can now compute again the flow equation of the Hamiltonian
vector field $\on{grad}^{\ga^*\om}(\tfrac12\|\quad\|_\ga^2)$:
For $g_t(t)\in TG$ we have 
\begin{displaymath}
(\pi_G,\ka^r)(g_t(t))=(g(t),u(t))=(g(t),T(\mu^{g(t)\i})g_t(t))
\end{displaymath}
and 
\begin{equation}
\p_t(g,u)=\on{grad}^{\ga^*\om}(\tfrac12\|\quad\|_\ga^2)(g,u)
=(T(\mu^g)u,u,-\ad(u)^\top u).
\tag{\nmb:{4}}\end{equation}
which reproduces the geodesic equation from \nmb!{3.2}.

\subsection*{\nmb.{4.3}. The momentum mapping}
Under the assumptions of \nmb!{4.2},
consider the right action of $G$ on $G$ and its prolongation to a right
action of $G$ on $TG$ in the Eulerian chart.
The corresponding fundamental vector fields 
are then given by:
\begin{align}
T(\mu^g)&:TG\to TG,
\notag\\
(\pi,\ka^r)T(\mu^g)T(\mu^h)X &= (\pi,\ka^r)T(\mu^{hg})X = (h.g,X), 
\quad
(h,X)\mapsto (hg,X)
\notag\\
\ze^{G\x\g}_X(h,Y)&
%=\ze_X^{G\x\g}(\pi_G,\ka^r)(T(\mu^h)Y) 
%= T(\pi_G,\ka^r)\ze_X^{TG}T(\mu^h)Y
%\\&
%= T(\pi_G,\ka^r)\p_t|_0T(\mu^{\exp(tX)})T(\mu^h)Y
%\\&
%=\p_t|_0 (\pi_G,\ka^r)T(\mu^{h.\exp(tX)})Y
%\\&
=\p_t|_0 (h.\exp(tX),Y) %=T\ka^r(h,\al;T(\mu_h)X,0)
=(T(\mu_h)X,0_Y)\in TG\x T\g
\tag{\nmb:{1}}\end{align}
%Likewise we get
%\begin{align*}
%T^*(\mu^g)&=T(\mu^{g\i})^*:T^*G\to T^*G,\quad (h,\al)\mapsto (h.g,\al)
%\\
%\ze^{G\x\g^*}_X(h,\al)&=\p_t|_0 (h.\exp(tX),\al) %=T\ka^r(h,\al;T(\mu_h)X,0)
%=(h,\al;\Ad(h)X,0)
%\end{align*}
Consider now the diagram from \nmb!{2.1} in the case of the weak symplectic
manifold $(M=G\x\g,\ga^*\om)$:
\begin{displaymath}
\xymatrix{
H^0 \ar[r] & 
     C^\infty_{\ga^*\om}(G\x\g,\mathbb R)  \ar[rr]^{\on{grad}^{\ga^*\om}} & & 
     \X(G\x \g,\ga^*\om) \ar[r] & H^1_{\ga^*\om}
\\
 & & \g \ar[lu]^{j} \ar[ru]_{\ze} & &
}\end{displaymath}
From the formulas derived above we see that for
$j(X)(h,Y):=\ga(\Ad(h)X,Y)$ we have: 
\begin{align*}
&
\ga(\on{grad}_2^\ga(j(X))(h,Y),Z) = d_2(j(X))(h,Y)(Z) = \ga(\Ad(h)X,Z)
\\&
\on{grad}_2^\ga(j(X))(h,Y)= \Ad(h)X
\\&
\ga(\on{grad}_1^\ga(j(X))(h,Y),T(\mu^h)Z) = d(j(X))(T(\mu^h)Z,Y,0)
\\&\quad
=\ga(d\Ad(T(\mu^h)Z)(X),Y) 
=\ga(((\ad\o\ka^r)\Ad)(T(\mu^h)Z)(X),Y) 
\\&\quad
=\ga(\ad(Z)\Ad(h)X,Y)  
=-\ga([\Ad(h)X,Z],Y)  
=-\ga(Z,\ad(\Ad(h)X)^\top Y)  
\\&
\ka^r(\on{grad}_1^\ga(j(X))(h,Y)) = - \ad(\Ad(h)X)^\top Y
\end{align*}
Thus the momentum mapping is 
\begin{align}
J:G\x \g&\to \g^*,\quad  
J\in C^\infty_{\ga^*\om}(G\x \g,\g^*) =
\notag\\ &= \{f\in
C^\infty(G\x \g,\g^*): \langle f(\quad),X \rangle\in
C^\infty_{\ga^*\om}(G\x\g)\; \forall X\in \g\}
\notag\\ 
\langle J(h,Y), X \rangle  &=j(X)(h,Y)=\ga(\Ad(h)X,Y) 
= \ga(\Ad(h)^\top Y,X)
\notag\\&
=\langle \ga(\Ad(h)^\top Y),X \rangle,
\notag\\
J(h,Y) &= \ga(\Ad(h)^\top Y) \in \g*
\notag\\
\bar J:=\ga\i\o J&:G\x \g \to \g, 
\notag\\
\bar J(h,Y) &= \Ad(h)^\top Y\in \g.
\tag{\nmb:{2}}\end{align}

\therosteritem{\nmb:{3}}
Note that the momentum mapping $J:G\x \g\to \g^*$ is equivariant for the
right $G$-action and the coadjoint action, and that $\bar J:G\x\g\to \g$
is equivariant for the right action $\Ad(\quad)^\top$ on $\g$: 
\begin{align*}
\langle J(hg,Y),X\rangle &= \langle\ga(\Ad(hg)^\top Y),X\rangle 
= \ga(\Ad(g)^\top\Ad(h)^\top Y,X)
\\&
= \ga(\Ad(h)^\top Y,\Ad(g)X)
= \langle\ga(\Ad(h)^\top Y),\Ad(g)X\rangle
\\&
= \langle\Ad(g)^*\ga(\Ad(h)^\top Y),X\rangle = \langle\Ad(g)^*J(h,Y),X\rangle
\\
\bar J(hg,Y) &= \Ad(hg)^\top Y = \Ad(g)^\top\bar J(h,Y).
\end{align*}

{\it \therosteritem{\nmb:{4}} For $x\in G\x\g$, the transposed mapping of
$d\bar J(x):T_x(G\x\g)\to \g$ is} 
\begin{displaymath}
d\bar J(x)^\top:\g\to T_x^*(G\x\g), \qquad d\bar J(x)^\top = (\ga^*\om)_x\o \ze,
\end{displaymath}
since for $\xi\in T_x(G\x\g)$ and $X\in \g$ we have 
\begin{displaymath}
\ga(d\bar J(\xi),X) = d\ga(\bar J,X)(\xi)
=dj(X)(\xi)=\langle(\ga^*\om)(\ze_X),\xi\rangle.
\end{displaymath}

{\it \therosteritem{\nmb:{5}} For $x\in G\x\g$, the closure  
 $\overline{d\bar J(T_x(G\x\g))}$ of the image of
$d\bar J(x):T_x(G\x\g)\to \g$ is the $\ga$-orthogonal space $\g_x^{\bot,\ga}$ 
of the isotropy Lie
algeba $\g_x:=\{X\in \g:\ze_X(x)=0\}$ in $\g$,} since the 
annihilator of the image is 
the kernel of the transposed mapping, 
\begin{displaymath}
\on{im}(dJ(x))^\o=\ker(dJ(x)^\top) = \ker((\ga^*\om)_x\o \ze) = 
\ker(\on{ev}_x\o\ze) = \g_x.
\end{displaymath}
Attention: the orthogonal space with respect to a weak inner product need
not be a complement.

\therosteritem{\nmb:{6}} For $(h,Y)\in G\x\g$, the $G$-orbit
$(h,Y).G=G\x\{Y\}$ is a submanifold of $G\x \g$.
{\it The kernel of $d\bar J(h,Y)$ is the symplectic orthogonal space}  
\begin{displaymath}
(T_{(h,Y)}(G\x\{Y\}))^{\bot,\ga^*\om} \subset T(\mu^h)\g\x\g
\end{displaymath}
since for the annihilator of the kernel we have 
\begin{align*}
\ker(&d\bar J(h,Y))^\o = \overline{\on{im}(d\bar J(h,Y)^\top)} 
=\overline{\on{im}((\ga^*\om_{(h,Y)}\o\ze)}, \quad\text{  by \thetag{\nmb|{4}},} 
\\&
= \overline{\{(\ga^*\om)_{(h,Y)}(\ze_X(x)):X\in \g\}} 
= \overline{(\ga^*\om)_{(h,Y)}(T_{(h,Y)}(G\x\{Y\}))},
\\&
= \bigl((T_{(h,Y)}(G\x\{Y\}))^{\bot,\ga^*\om}\bigr)^\o.
\end{align*}
The last equality holds by the bipolar theorem for the usual duality
pairing.

{\it \therosteritem{\nmb:{7}} Thus, for $(h,Y)\in G\x\g$,} 
\begin{align*}
&T(\mu^h)X_1,Y_1)\in \ker(d\bar J(h,Y))
\\&\iff
(\ga^*\om)_{(h,Y)}((T(\mu^h)X_1,Y_1),(T(\mu^h)Z,0))=0
\text{ for all }Z\in \g 
\\&\iff
0=0-\ga(Y_1,Z)-\ga(Y,[X_1,Z])=-\ga(Y_1+\ad(X_1)^\top Y,Z)
\; \forall\; Z\in \g
\\&\iff
Y_1=-\ad(X_1)^\top Y.
\end{align*}

{\it \therosteritem{\nmb:{8}} (Emmy Noether's theorem)
Let $h\in C^\infty_\om(G\x \g)$ be a Hamiltonian function which is invariant 
under the right $G$-action. Then 
$d\bar J(\on{grad}^{\ga^*\om}(h))=0\in \g$ and also
$dJ(\on{grad}^{\ga^*\om}(h))=0\in\ga(\g)\subseteq\g^*$.
Thus the momentum mappings 
$\bar J:G\x \g\to \g$ and $J:G\x \g\to \ga(\g)\subset \g^*$ are 
constant on each trajectory (if it exists) 
of the Hamiltonian vector field $\on{grad}^{\ga^*\om}(h)$.}
Namely, consider the function $\ga(\bar J,X)=\langle J,X \rangle = j(X)$.   
\begin{align*}
\ga(d\bar J(\on{grad}^{\ga^*\om}(h)),X) 
&= \on{grad}^{\ga^*\om}(h)(\ga(\bar J,X)) =
\\&
= \{h,\ga(\bar J,X)\} = -\{j(X),h\} = -\ze_X(h) =0.
\\
\langle dJ(\on{grad}^{\ga^*\om}(h)),X \rangle
&= \on{grad}^{\ga^*\om}(h)(\langle J,X\rangle) =
\\&
= \{h,j(X)\} = -\{j(X),h\} = -\ze_X(h) =0.
\end{align*}

\subsection*{\nmb.{4.4}. The geodesic equation via conserved momentum}
We consider a smooth curve $t\mapsto g(t)$ in $G$ and
$(\pi_G,\ka^r)g_t(t)=(g(t),u(t))=(g(t),T(\mu^{g(t)\i})g_t(t))$ as in
\nmb!{4.2.4}. 
Applying $\bar J:G\x \g\to \g$ to it we get $\bar J(g,u)=\Ad(g)^\top u$.
We claim that {\it the curves $t\mapsto g(t)$ in $G$ for which 
$\bar J(g(t),u(t))$ is constant in $t$ are exactly the geodesics in
$(G,\ga)$}.
Namely, by \nmb!{3.1} we have
\begin{align*}
0 &= \p_t \Ad(g(t))^\top u(t) = \big((\ad\o\ka^r)(\p_t
g(t)).\Ad(g(t))\big)^\top u(t) + \Ad(g(t))^\top \p_t u(t)
\\&
=\Ad(g(t))^\top\big(\ad(u(t))^\top u(t) + u_t(t)\big)
\\&
\iff \quad u_t = - \ad(u)^\top u.
\end{align*}

\subsection*{\nmb.{4.5}. Symplectic reduction to transposed adjoint orbits }
Under the assumptions of \nmb!{4.2} we have the following:

\therosteritem{\nmb:{1}} {\it For $X\in \bar J(G\x \g)$ the inverse image 
$\bar J\i(X)\subset G\x \g$ is a manifold.}
Namely, it is the graph of a smooth mapping:
\begin{align*}
\bar J\i(X) &=\{(h,Y)\in G\x\g: \Ad(h)^\top Y= X\}
\\&
=\{(h,\Ad(h\i)^\top X): h\in G\} \West{\cong}{} G. \qed
\end{align*}

\therosteritem{\nmb:{2}} {\it At any point of $\bar J\i(X)$, 
the kernel of the
pullback of the symplectic form $\ga^*\om$ on $G\x \g$ from
\nmb!{4.2.1} equals the tangent space 
to the orbit of the isotropy group $G_X:=\{g\in G:\Ad(g)^\top X=X\}$
through that point}. 

For $(h,Y=\Ad(h\i)^\top X)\in \bar J\i(X)$ the $G_X$-orbit is 
$h.G_X\x\{Y\}$ and its tangent space at 
$(h,Y)$ is 
$T(\mu_h)\g_X\x 0$ where $\g_X=\{Z\in\g:\ad(Z)^\top X=0\}$.
The tangent space at $(h,Y)$ of $\bar J\i(X)$ is  
\begin{align*}
T_{(h,\Ad(h\i)^\top X)}&\bar J\i(X)
=\{\p_t|_0 (\exp(tZ).h,\Ad((\exp(tZ).h)\i)^\top X): Z\in \g \}
\\&
=\{(T(\mu^h)Z,-\ad(Z)^\top\Ad(h\i)^\top X):Z\in \g\}\subset T_hG\x \g.
\end{align*}
For $Z_1,Z_2\in \g$ consider the tangent vectors 
$(T(\mu^h)\Ad(h)Z_1,Y,-\ad(Z_1)X)$ and 
$(T(\mu^h)Z,Y,-\ad(Z)^\top\Ad(h\i)^\top X)$
in $T_{(h,Y)}\bar J\i(X)$.
From \nmb!{4.2.1}, we get 
\begin{align*}
&(\ga^*\om)_{(h,Y)}\big(\!(T(\mu^h)\Ad(h)Z_1,\!-\!\ad(Z_1)^\top X),
  (T(\mu^h)Z_2,\!-\!\ad(Z_2)^\top\Ad(h\i)^\top X)\!\big)
\\&
= \ga(-\ad(Z_2)^\top\Ad(h\i)^\top X,\Ad(h)Z_1) 
- \ga(-\ad(Z_1)^\top X,Z_2) 
\\&\qquad
-\ga(Y,[\Ad(h)Z_1,Z_2])
\\&
= -\ga(\Ad(h\i)^\top X,\ad(Z_2)\Ad(h)Z_1) 
+\ga(\ad(Z_1)^\top X,Z_2) 
-\\&\qquad
-\ga(\Ad(h\i)^\top X,[\Ad(h)Z_1,Z_2]) 
\\&= \ga(\ad(Z_1)^\top X,Z_2) = 0 \quad \forall Z_2\in \g \iff Z_1 \in \g_X.
\qquad\qed
\end{align*}

\therosteritem{\nmb:{3}} {\it The reduced symplectic manifold $\bar
J\i(X)/G_X$ with symplectic form induced by $\ga^*\om|_{\bar J\i(X)}$ is
symplectomorphic to the adjoint orbit $\Ad(G)^\top X\subset \g$ with
symplectic form the pullback via $\ga:\g\to \g^*$ of the Kostant Kirillov
Souriou form 
$$
\om_\al(\ad(Y_1)^*\al,\ad(Y_2)^*\al)=\langle\al,[Y_1,Y_2]\rangle
$$
which is given by} 
\begin{align*}
&\om_Z(\ad(Y_1)^\top Z,\ad(Y_2)^\top Z) 
  = \om_{\ga(Z)}(\ga\ad(Y_1)^\top Z,\ga\ad(Y_2)^\top Z)
\\&\qquad
= \om_{\ga(Z)}(\ad(Y_1)^*\ga Z,\ad(Y_2)^*\ga Z) = \langle \ga(Z),[Y_1,Y_2]\rangle
=\ga(Z,[Y_1,Y_2]),
\end{align*}
since for $Y,Z,U\in \g$ we get
\begin{multline*}
\langle \ga\ad(Y)^\top Z,U \rangle= \ga(\ad(Y)^\top Z,U) 
=\ga(Z,\ad(Y)U)=
\\
= \langle\ga(Z),\ad(Y)U\rangle 
= \langle \ad(Y)^*\ga(Z),U\rangle.
\end{multline*}
The quotient space is 
$\bar J\i(X)/G_X = \{(h.G_X,\Ad(h\i)^\top X): h\in G\}\cong \Ad(G)^\top
X\cong G/G_X$. The 2-form $\ga^*\om|_{\bar J\i(X)}$ induces a symplectic
form on the quotient by {\nmb|{2}} and it remains to check
that it agrees with the pullback of the Kirillov Kostant Souriou symplectic
form. But this is obvious from the last computation in {\nmb|{2}}
(for the special case $h=e$ if the reader insists). \qed

\therosteritem{\nmb:{4}} Reconsider the geodesic equation on the 
reduced space  $\bar J\i(X)/G_X \cong \Ad(G)^\top X$. The energy function is 
$E(\Ad(g)^\top X)= \frac12\|\Ad(g)^\top X\|^2_\ga$. 
For $Z=\Ad(g)^\top X\in\Ad(G)^\top X$ the tangent space is given by 
$T_Z(\Ad(G)^\top X)=\{\ad(Y)^\top Z:Y\in \g\}$. 
We look for the Hamiltonian vector field of $E$ in the form
$\on{grad}^\om E(Z)=\ad(H_E(Z))^\top Z$,
for a vector field $H_E$.
The differential of the
energy function is
$dE(Z)(\ad(Y)^\top Z)=\ga(Z,\ad(Y)^\top Z)=\ga([Y,Z],Z)$
which equals 
$\om_Z(\on{grad}^\om E(Z),\ad(Y)^\top Z)
=\om_Z(\ad(H_E(Z))^\top Z,\ad(Y)^\top Z)= \ga(Z,[H_E(Z),Y])$
from which we conclude that $H_E(Z)=-Z$ will do (which is defined up to
annihilator of $Z$). Thus 
$\on{grad}^\om E(Z)=-\ad(Z)^\top Z$ which leads us back to the geodesic
equation $u_t=-\ad(u)^\top u$ again.

\section*{\nmb0{5}. Vanishing $H^0$-geodesic distance on groups of
diffeomorphisms}

This section is based on \cite{30}.

\subsection*{\nmb.{5.1} The $H^0$-metric on groups of
diffeomorphisms}
Let $(N,g)$ be a smooth connected Riemannian manifold, and let
$\on{Diff}_c(N)$ be the group of all diffeomorphisms with
compact support on $N$, and let $\on{Diff}_0(N)$ be the
subgroup of those which are diffeotopic in $\on{Diff}_c(N)$ to the identity;
this is the connected component of the identity in $\on{Diff}_c(N)$,
which is a regular Lie group in the sense of \cite{32}, section~38. This is
proved in \cite{23}, section~42. The Lie algebra is $\X_c(N)$, the space of all 
smooth vector fields with compact support on $N$, with the negative of the
usual bracket of vector fields as Lie bracket. Moreover, $\on{Diff}_0(N)$
is a simple group (has no nontrivial normal subgroups), see \cite{E}, 
\cite{T}, \cite{Ma}. 
The {\it right invariant} $H^0$-metric 
on $\on{Diff}_0(N)$ is then given as follows, where $h,k:N\to TN$ are 
vector fields with compact support along $\ph$ and where 
$X=h\o\ph\i, Y=k\o\ps\i\in\X_c(N)$: 
\begin{align}
\ga^0_\ph(h,k) &= \int_N g(h,k)\on{vol}(\ph^*g)
= \int_N g(X\o\ph,Y\o\ph)\ph^*\on{vol}(g)
\notag\\&
= \int_N g(X,Y)\on{vol}(g).
\tag{\nmb:{1}}\end{align}

\begin{proclaim}{\nmb.{5.2}.  Theorem}
Geodesic distance on $\on{Diff}_0(N)$ with respect to the $H^0$-metric vanishes.
\end{proclaim}

\begin{demo}{Proof}
Let $[0,1]\ni t\mapsto\ph(t,\quad)$ be a smooth curve in
$\on{Diff}_0(N)$ between $\ph_0$ and $\ph_1$. Consider the
curve $u=\ph_t\o\ph\i$ in $\X_c(N)$, the right logarithmic
derivative. Then for the length and the energy we have:  
\begin{align}
L_{\ga^0}(\ph)&=\int_0^1  \sqrt{\int_N \|u\|^2_g\on{vol(g)}}\;dt
\tag{\nmb:{1}}\\
E_{\ga^0}(\ph)&=\int_0^1  \int_N \|u\|^2_g\on{vol(g)}\,dt
\tag{\nmb:{2}}\\
L_{\ga^0}(\ph)^2&\le E_{\ga^0}(\ph)
\tag{\nmb:{3}}\end{align}

\noindent\thetag{\nmb:{4}} Let us denote by $\on{Diff}_0(N)^{E=0}$ the set of all
diffeomorphisms $\ph\in\on{Diff}_0(N)$ with the following
property: For each $\ep>0$ there exists a smooth curve from the
identity to $\ph$ in $\on{Diff}_0(N)$ with energy $\le \ep$.

\noindent\thetag{\nmb:{5}} {\it We claim that $\on{Diff}_0(N)^{E=0}$ coincides with
the set of all diffeomorphisms which can be reached from the
identity by a smooth curve of arbitraily short $\ga^0$-length.}
This follows by {\nmb|{3}}.

\noindent\thetag{\nmb:{6}} {\it We claim that $\on{Diff}_0(N)^{E=0}$
is a normal subgroup of $\on{Diff}_0(N)$.}
Let $\ph_1\in\on{Diff}_0(N)^{E=0}$ and $\ps\in\on{Diff}_0(N)$.
For any smooth curve $t\mapsto\ph(t,\quad)$ from the identity
to $\ph_1$ with energy $E_{\ga^0}(\ph)<\ep$ we have
\begin{align*}
E_{\ga^0}&(\ps\i\o\ph\o\ps) 
  = \int_0^1 \int_N\|T\ps\i\o\ph_t\o\ps\|_g^2
  \on{vol}((\ps\i\o\ph\o\ps)^*g)
\\&
\le \sup_{x\in N}\|T_x\ps\i\|^2 \cdot 
  \int_0^1 \int_N\|\ph_t\o\ps\|_g^2 (\ph\o\ps)^*\on{vol}((\ps\i)^*g)
\\&
\le \sup_{x\in N}\|T_x\ps\i\|^2 \cdot
  \sup_{x\in N}\frac{\on{vol}((\ps\i)^*g)}{\on{vol}(g)}\cdot
  \int_0^1 \int_N\|\ph_t\o\ps\|_g^2\,
  (\ph\o\ps)^*\on{vol}(g)
\\&
\le \sup_{x\in N}\|T_x\ps\i\|^2 \cdot
  \sup_{x\in N}\frac{\on{vol}((\ps\i)^*g)}{\on{vol}(g)}\cdot
  E_{\ga^0}(\ph).
\end{align*}
Since $\ps$ is a diffeomorphism with compact support, the two
suprema are bounded.
Thus $\ps\i\o\ph_1\o\ps\in\on{Diff}_0(N)^{E=0}$.
 
\noindent\thetag{\nmb:{7}} {\it We claim that  
$\on{Diff}_0(N)^{E=0}$
is a non-trivial subgroup.} In view of the simplicity of
$\on{Diff}_0(N)$ mentioned in \nmb!{5.1} this concludes the
proof. 

It remains to find a non-trivial diffeomorphism in $\on{Diff}_0(N)^{E=0}$. 
The idea is to use compression waves. The basic case is this: take
any non-decreasing smooth function $f:{\mathbb R} \rightarrow {\mathbb R}$
such that $f(x)\equiv 0$ if $x \ll 0$ and $f(x) \equiv 1$ if $x \gg 0$.
Define 
\begin{displaymath}
\ph(t,x) = x + f(t-\la x)\end{displaymath}
where $\la < 1/\max(f')$. Note that 
\begin{displaymath}
\ph_x(t,x) = 1 - \la f'(t-\la x) > 0,\end{displaymath}
hence each map $\ph(t,\quad)$ is a diffeomorphism of $\mathbb R$ and 
we have a path in the group of diffeomorphisms of $\mathbb R$. These 
maps are not the identity outside a compact set however. In fact, 
$\ph(x)=x+1$ if $x \ll 0$ and $\ph(x)=x$ if $x \gg 0$. As 
$t \rightarrow -\infty$, the map $\ph(t,\quad)$ approaches
the identity uniformly on compact subsets, 
while as $t \rightarrow +\infty$, the map approaches
translation by 1. This path is a moving compression wave which pushes all
points forward by a distance 1 as it passes. We calculate its energy between
two times $t_0$ and $t_1$:
\begin{align*}
E_{t_0}^{t_1}(\ph) &= \int_{t_0}^{t_1} \int_{\mathbb R} 
\ph_t(t,\ph(t,\quad)^{-1}(x))^2 dx\, dt 
%\\&
= \int_{t_0}^{t_1} \int_{\mathbb R} 
\ph_t(t,y)^2 \ph_y(t,y) dy\, dt  \\
&= \int_{t_0}^{t_1} \int_{\mathbb R} 
f'(z)^2 \cdot (1-\la f'(z)) \frac{dz}{\la}\, dt \\
& \le \frac{\max {f'}^2}{\la} \cdot (t_1-t_0) \cdot
\int_{\on{supp}(f')} (1-\la f'(z)) dz 
\end{align*}
If we let $\la = 1-\ep$ and consider the specific $f$ given by the convolution
\begin{displaymath}
f(z) = \max(0,\min(1,z))\star G_\ep(z),\end{displaymath}
where $G_\ep$ is a smoothing kernel supported on $[-\ep,+\ep]$, then the
integral is bounded by $3\ep$, hence
\begin{displaymath}
E_{t_0}^{t_1}(\ph) \le (t_1-t_0)\tfrac{3\ep}{1-\ep}.\end{displaymath}

We next need to adapt this path so that it has compact support. To do this
we have to start and stop the compression wave, which we do by giving 
it variable length. Let:
\begin{displaymath}
f_\ep(z,a) = \max(0,\min(a,z)) \star (G_\ep(z)G_\ep(a)).\end{displaymath}
The starting wave can be defined by:
\begin{displaymath}
\ph_\ep(t,x) = x + f_\ep(t-\la x, g(x)),\quad \la <1, \quad 
g \text{ increasing}.\end{displaymath}
Note that the path of an individual particle $x$ hits the wave at $t=\la x-\ep$
and leaves it at $t=\la x + g(x) + \ep$, having moved forward to $x+g(x)$. 
Calculate the derivatives:
\begin{align*}
(f_\ep)_z &= I_{0 \le z \le a} \star (G_\ep(z)G_\ep(a)) \in [0,1]\\
(f_\ep)_a &= I_{0 \le a \le z} \star (G_\ep(z)G_\ep(a)) \in [0,1] \\
(\ph_\ep)_t &= (f_\ep)_z(t-\la x, g(x))\\
(\ph_\ep)_x &= 1 - \la (f_\ep)_z(t-\la x, g(x))
+(f_\ep)_a(t-\la x, g(x))\cdot g'(x) > 0.
\end{align*}
This gives us:
\begin{align*}
E_{t_0}^{t_1}(\ph) &= \int_{t_0}^{t_1} \int_{\mathbb R} (\ph_\ep)_t^2
(\ph_\ep)_x dx \, dt \\
&\le \int_{t_0}^{t_1} \int_{\mathbb R} (f_\ep)_z^2(t-\la x,g(x))\cdot
(1-\la (f_\ep)_z(t-\la x, g(x)))dx\, dt\\
&\quad+\int_{t_0}^{t_1} \int_{\mathbb R}(f_\ep)_z^2(t-\la x,g(x))\cdot
(f_\ep)_a(t-\la x, g(x)) g'(x) dx\, dt
\end{align*}
The first integral can be bounded as in the original discussion.
The second integral is also small because the support of the $z$-derivative
is $-\ep \le t-\la x \le g(x)+\ep$, while the support of the $a$-derivative
is $-\ep \le g(x) \le t-\la x + \ep$, so together $|g(x)-(t-\la x)| \le \ep$.
Now define $x_1$ and $x_2$ by $g(x_1)+\la x_1 = t+\ep$ and 
$g(x_0)+\la x_0 = t-\ep$.
Then the inner integral is bounded by
\begin{displaymath}
\int_{|g(x)+\la x -t|\le \ep} g'(x) dx = g(x_1)-g(x_0) \le 2\ep,\end{displaymath}
and the whole second term is bounded by ${2\ep}(t_1-t_0)$. Thus the
length is $O(\ep)$.

The end of the wave can be handled by playing the beginning backwards. If
the distance that a point $x$ moves when the wave passes it is to be $g(x)$,
so that the final diffeomorphism is $x \mapsto x+g(x)$, then let $b = \max(g)$ 
and use the above definition of $\ph$ while $g' > 0$. The modification
when $g' < 0$ (but $g' > -1$ in order for $x \mapsto x+g(x)$ 
to have positive derivative) is given by:
\begin{displaymath}
\ph_\ep(t,x) = x + f_\ep(t-\la x -(1-\la)(b-g(x)),g(x)).\end{displaymath}
Consider the figure showing the trajectories $\ph_\ep(t,x)$ for sample
values of $x$.
\begin{figure*} \begin{center}
\epsfig{width=8cm,file=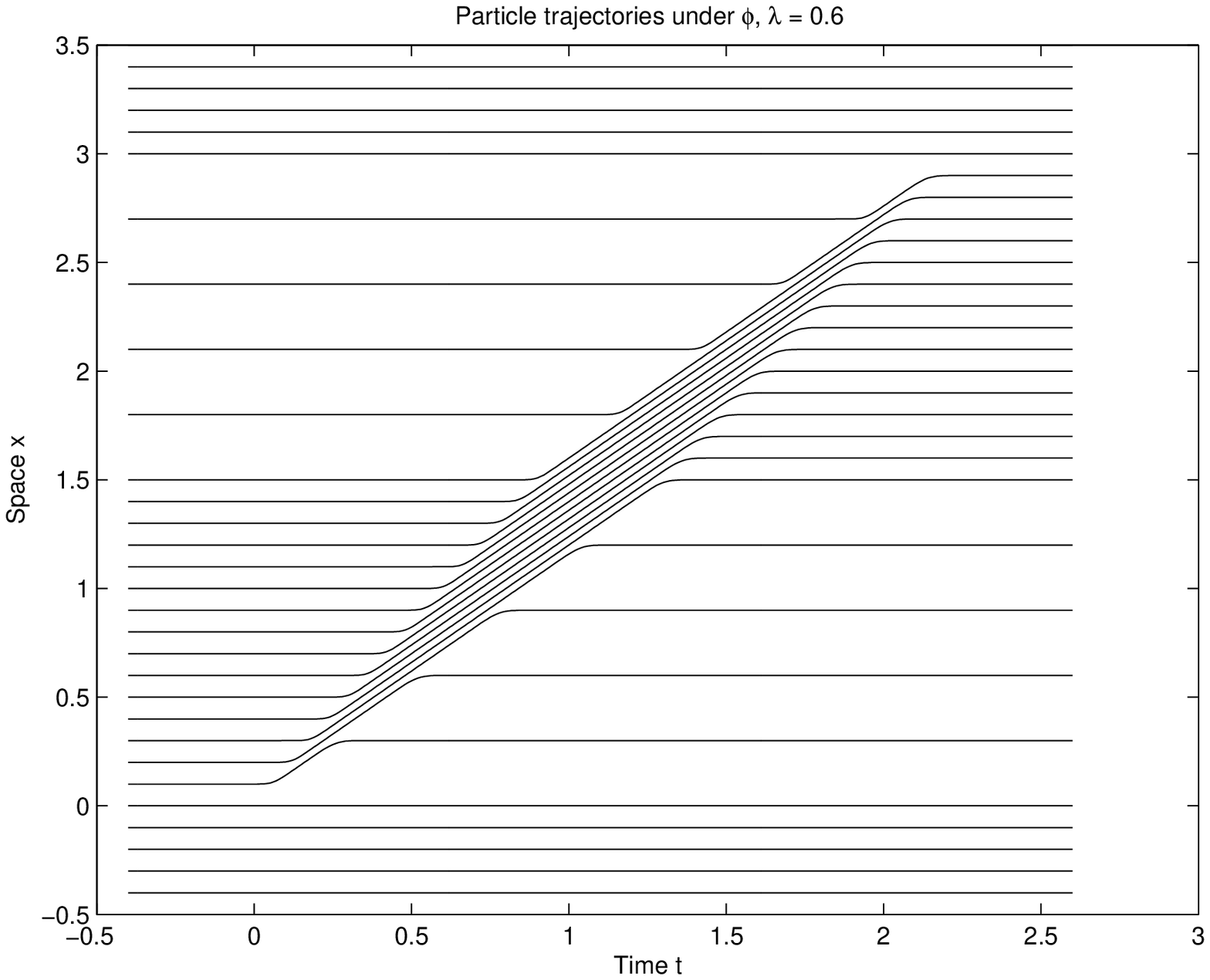}
\end{center} 
\end{figure*}

It remains to show that $\on{Diff}_0(N)^{E=0}$ is a nontrivial subgroup 
for an arbitrary Riemannian manifold. We choose a piece of a unit
speed geodesic containing no conjugate points in $N$ and Fermi coordinates 
along this geodesic; so we can assume that we are in an open set in 
$\mathbb R^m$ which is a tube around a piece of the $u^1$-axis. Now we 
use a small bump function in the slice orthogonal to the $u^1$-axis 
and multiply it with the construction from above for the coordinate $u^1$. 
Then it follows that we get a nontrivial diffeomorphism in 
$\on{Diff}_0(N)^{E=0}$ again. 
\qed\end{demo}

\subsection*{Remark} 
Theorem \nmb!{5.2} can be proved directly without
the help of the simplicity of $\on{Diff}_0(N)$. For $N=\mathbb R$ one can
use the method of \nmb!{5.2.7} in the parameter space of a
curve, and for general $N$ one can use a Morse function on $N$ to produce a
special coordinate for applying the same method.

\subsection*{\nmb.{5.3}. Geodesics and sectional curvature for $\ga^0$ on $\on{Diff}(N)$
}
According to \nmb!{3.2}, \nmb!{3.4}, or \nmb!{4.4}, for a right invariant
weak Riemannian metric $G$ on an (possibly infinite dimensional) Lie group
the geodesic equation and the curvature are given in terms of the
transposed operator (with respect to $G$, if it exists) of the Lie bracket
by the following formulas:
\begin{align*} 
u_t&=-\on{ad}(u)^*u, \quad u = \ph_t\o \ph\i
\\                                                                                      
G(\on{ad}(X)^*Y,Z)&:=G(Y,\on{ad}(X)Z)
\\
4G(R(X,Y)X,Y) &=  3G(\on{ad}(X)Y,\on{ad}(X)Y)                   
- 2G(\on{ad}(Y)^*X,\on{ad}(X)Y)  
\\&\quad
- 2G( \on{ad}(X)^*Y,\on{ad}(Y)X )                   
+ 4G( \on{ad}(X)^* X,\on{ad}(Y)^* Y )   
\\&\quad
- G( \on{ad}(X)^* Y+\on{ad}(Y)^* X,\on{ad}(X)^* Y+\on{ad}(Y)^* X ) 
\end{align*}
In our case, for $\on{Diff}_0(N)$, we have $\on{ad}(X)Y=-[X,Y]$ 
(the bracket on the Lie algebra $\X_c(N)$ of vector fields with
compact support is the negative of the usual one), and:
\begin{align*}
\ga^0(X,Y) &= \int_N g(X,Y)\on{vol}(g)
\\
\ga^0(\on{ad}(Y)^*X,Z)&= \ga^0(X,-[Y,Z]) = \int_N g(X,-\L_YZ)\on{vol}(g) 
\\&
= \int_N g\Bigl(\L_YX+(g\i\L_Yg)X+\on{div}^g(Y)X,Z\Bigl)\on{vol}(g)
\\
\on{ad}(Y)^*&= \L_Y+g\i\L_Y(g)+\on{div}^g(Y)\on{Id}_TN = \L_Y +\be(Y),
\end{align*}
where the tensor field $\be(Y)=g\i\L_Y(g)+\on{div}^g(Y)\on{Id}:TN\to TN$ 
is self adjoint with respect to $g$. 
Thus the geodesic equation is 
\begin{displaymath}
u_t = -(g\i\L_u(g))(u)-\on{div}^g(u)u=-\be(u)u, \qquad u = \ph_t\o \ph\i.
\end{displaymath}
The main part of the sectional curvature is given by:
\begin{align*} 
&4G(R(X,Y)X,Y) =
\\&
= \int_N\Bigl( 3\|[X,Y]\|_g^2
+ 2g((\L_Y+\be(Y))X,[X,Y])  + 2g((\L_X+\be(X))Y,[Y,X] )                   
\\&\qquad\qquad
+ 4g(\be(X)X,\be(Y)Y) -\|\be(X)Y+\be(Y)X\|_g^2
\Bigr)\on{vol}(g)
\\&
= \int_N\Bigl(-\|\be(X)Y-\be(Y)X+[X,Y]\|_g^2
-4g([\be(X),\be(Y)]X,Y)
\Bigr)\on{vol}(g)
\end{align*}
So sectional curvature consists of a part which is visibly non-negative, and
another part which is difficult to decompose further. 

\subsection*{\nmb.{5.4} Example: $n$-dimensional analog of Burgers' equation
}
For $(N,g)=(\mathbb R^n,\on{can})$ or $((S^1)^n,\on{can})$ we have:
\begin{align*}
(\on{ad}(X)Y)^k&=\sum_i((\p_iX^k)Y^i-X^i(\p_iY^k))
\\
%\ga^0(\on{ad}(X)Y, Z)
%     &= \int_{\mathbb R^n}\langle dX.Y-dY.X,Z \rangle dx\\ 
%&= \int_{\mathbb R^n}\sum_{i,j} \Bigl( 
%     (\p_jX^i)Y^jZ^i -(\p_jY^i)X^jZ^i \Bigr) dx\\ 
%&= \int_{\mathbb R^n}\sum_{i,j} \Bigl( 
%     (\p_jX^i)Y^jZ^i +Y^i(\p_jX^j)Z^i +Y^iX^j(\p_jZ^i)\Bigr) dx\\ 
%&= \int_{\mathbb R^n}\sum_{i,k} Y^k\Bigl( 
%     (\p_kX^i)Z^i + (\p_iX^i)Z^k + X^i(\p_iZ^k)\Bigr) dx\\ 
(\on{ad}(X)^* Z)^k &= \sum_{i} \Bigl( 
     (\p_kX^i)Z^i + (\p_iX^i)Z^k + X^i(\p_iZ^k)\Bigr), 
\end{align*}
so that the geodesic equation is given by 
\begin{displaymath}
\p_t u^k = - (\on{ad}(u)^\top u)^k  
=  -\sum_{i} \Bigl((\p_ku^i)u^i + (\p_iu^i)u^k + u^i(\p_iu^k)\Bigr), 
\end{displaymath}
the $n$-dimensional analog of Burgers' equation.

\subsection*{\nmb.{5.5}. Stronger metrics on $\on{Diff}_0(N)$}
A very small strengthening of the weak Riemannian $H^0$-metric on 
$\on{Diff}_0(N)$ makes it into a true metric. We define the stronger 
right invariant semi-Riemannian metric by the formula:
\begin{displaymath}
G^A_\ph(X\o\ph,Y\o\ph)
= \int_N (g(X,Y)+A\on{div}_g(X).\on{div}_g(Y)) \on{vol}(g).
\end{displaymath}
Then the following holds:

\begin{proclaim}{Theorem}
For any distinct diffeomorphisms $\ph_0, \ph_1$, the infimum of the
lengths of all paths from $\ph_0$ to $\ph_1$ with respect to $G^A$
is positive.
\end{proclaim}

\begin{demo}{Proof}
We may suppose that $\ph_0=\on{Id}_N$. 
If $\ph_1 \ne \on{Id}_N$, there are two
functions $\rh$ and $f$ on $N$ with compact support such that:
\begin{displaymath}
\int_N \rh(y) f(\ph_1(y)) \on{vol}(g)(y) \ne 
\int_N \rh(y) f(y) \on{vol}(g)(y). \end{displaymath}
Now consider any path $\ph(t,y)$ between $\ph_0=\on{Id}_N$ to $\ph_1$
with left logarithmic derivative
$u = T(\ph)\i\o\ph_t$ and a path in $\X_c(N)$. 
Then we have:
\begin{align*}
\int_N &\rh (f\o\ph_1) \on{vol}(g) -\int_N \rh f\on{vol}(g) 
=\int_0^1 \int_N \rh \p t f(\ph(t,\quad) \on{vol}(g) dt
\\&
= \int_0^1\int_N \rh (df.\ph_t)\on{vol}(g)\, dt
= \int_0^1\int_N \rh (df.T\ph. u)\on{vol}(g) dt
\\&
= \int_0^1\int_N (df.T\ph.(\ph u))\on{vol}(g) dt
\end{align*}
Locally, on orientable pieces of $N$, we have:  
\begin{align*}
\on{div}((f\o\ph)\rh u)&\on{vol}(g)=\Cal L_{(f\o\ph)\rh u}\on{vol}(g)
= (i_{(f\o\ph)\rh u}d+di_{(f\o\ph)\rh u})\on{vol}(g)
\\&
= d((f\o\ph)i_{\rh u}\on{vol}(g))
=d(f\o\ph)\wedge i_{\rh u}\on{vol}(g)+\rh\on{div}(u)\on{vol}(g),
\\&
=d(f\o\ph)(\rh u) \on{vol}(g)+(f\o\ph)\on{div}(\rh u)\on{vol}(g),
\qquad\text{  since }
\\
d(f\o\ph)\wedge i_{\rh u}&\on{vol}(g) 
= - i_{\rh u}(d(f\o\ph)\wedge \on{vol}(g)) + 
(i_{\rh u} d(f\o\ph))\on{vol}(g)).  
\end{align*}
Thus on $N$ we have:
\begin{align*}
0 &= \int_N \on{div}((f\o\ph)\rh u)\on{vol}(g)
\\&
= \int_N d(f\o\ph)(\rh u) \on{vol}(g)+\int_N(f\o\ph)\on{div}(\rh u)\on{vol}(g) 
\end{align*}
and hence 
\begin{align*}
0&\le\Big|\int_N \rh (f\o\ph_1) \on{vol}(g) -\int_N \rh f\on{vol}(g)\Big| 
= \Big|\int_0^1\int_N d(f\o\ph)(\ph u))\on{vol}(g) dt \Big|
\\&
= \Big|\int_0^1\int_N -(f\o\ph)\on{div}(\rh u)\on{vol}(g) dt \Big|
\\&
\le \sup |f| \cdot \int_0^1 \sqrt{\int_N C_\rh \|u\|^2 + 
C'_\rh|\on{div}(u)|^2 \on{vol}(g)}\;dt
\end{align*}
for constants $C_\rh, C'_\rh$ depending only on $\rh$. 
Clearly the right hand 
side gives a lower bound for the length of any path from $\ph_0$ to $\ph_1$.
\qed\end{demo}

\subsection*{\nmb.{5.6}. Geodesics and sectional curvature for $G^A$ on
$\on{Diff}(\mathbb R)$}
We consider the groups $\on{Diff}_c(\mathbb R)$ or $\on{Diff}(S^1)$ with Lie  
algebras $\X_c(\mathbb R)$ or $\X(S^1)$ whose Lie brackets are 
$\on{ad}(X)Y=-[X,Y]=X'Y-XY'$.
The $G^A$-metric equals the $H^1$-metric on $\X_c(\mathbb R)$, and we have: 
\begin{align*} 
G^A(X,Y) &= \int_{\mathbb R} (XY + AX'Y')dx=
\int_{\mathbb R}X(1-A\p_x^2)Y\,dx,\\ 
G^A(\on{ad}&(X)^* Y, Z)= \int_{\mathbb R} (YX'Z-YXZ'+AY'(X'Z-XZ')') dx\\ 
%&= \int_{\mathbb R} (YX'Z+Y'XZ+YX'Z-AY''(X'Z-XZ')) dx\\ 
%&= \int_{\mathbb R} Z(YX'+Y'X+YX'-AY''X'-AY'''X-AY''X') dx\\ 
&= \int_{\mathbb R}  
Z(1-\p_x^2)(1-\p_x^2)\i(2YX'+Y'X-2AY''X'-AY'''X) dx,\\ 
\on{ad}(X)^* Y &= (1-\p_x^2)\i(2YX'+Y'X-2AY''X'-AY'''X)\\ 
\on{ad}(X)^* &= (1-\p_x^2)\i(2X'+X\p_x)(1-A\p_x^2)
\end{align*}
so that the geodesic equation in 
Eulerian representation   
$u=(\p_t f)\o f\i \in \X_c(\mathbb R)$ or $\X(S^1)$ is
\begin{align*} 
\p_t u &= - \on{ad}(u)^* u
=  -(1-\p_x^2)\i(3uu'-2Au''u'-Au'''u),\text{ or }\\ 
u_t-u_{txx}&= Au_{xxx}.u+2Au_{xx}.u_x-3u_x.u,
\end{align*}
which for $A=1$ is the dispersionless version of the 
{\it Camassa-Holm equation}, see \nmb!{7.3.4}. 
Note that here geodesic distance is a well defined metric describing the
topology. 

\section*{\nmb0{6}. 
The regular Lie group 
of rapidly decreasing diffeomorphisms
} 

\begin{proclaim}{\nmb.{6.1}.  Lemma} 
For smooth functions of one variable we have: 
\begin{align*}
(f\o g)^{(p)}(x) &= p! \sum_{m\ge 0} \frac{f^{(m)}(g(x))}{m!} 
  \sum_{\substack{\al\in \mathbb N_{>0}^m\\ \al_1+\dots+\al_m =p}} \prod_{i=1}^m
  \frac{g^{(\al_i)}(x)}{\al_i!}
\\&
= \sum_{m\ge 0} f^{(m)}(g(x)) 
  \sum_{\substack{\la=(\la_n)\in \mathbb N_{\ge 0}^{\mathbb N_{>0}}\\ \sum_n\la_n=m \\
  \sum_n\la_n n= p}}
  \frac{p!}{\la!}\prod_{n>0}
  \left(\frac{g^{(n)}(x)}{n!}\right)^{\la_n}
\end{align*}
Let $f\in C^\infty(\mathbb R^k)$ and let 
$g=(g_1,\dots,g_k)\in C^\infty(\mathbb R^n,\mathbb R^k)$. Then for a 
multiindex $\ga\in \mathbb N^n$ the partial 
derivative $\partial^\ga(f\o g)(x)$ of the composition is given by 
the following formula, 
where we use multiindex-notation heavily.
\begin{align*}
&\p^\ga(f\o g)(x) = \\
&=\sum_{\be\in\mathbb N^k} 
(\p^\be f)(g(x))
\sum_{\substack{\la=(\la_{i\al})\in \mathbb N^{k\x (\mathbb N^n\setminus 0)}\\
          \sum_\al \la_{i\al}=\be_i\\
          \sum_{i\al}\la_{i\al}\al = \ga}}
\frac{\ga!}{\la!} \prod_{\substack{\al\in \mathbb N^n\\ \al>0}} 
\left(\frac1{\al!}\right)^{\sum_i\la_{i\al}} 
\prod_{i,\al>0}(\p^\al g_i(x))^{\la_{i\al}} \\
&=\sum_{\substack{\la=(\la_{i\al})\in \mathbb N^{k\x (\mathbb N^n\setminus 0)}\\
          \sum_{i\al}\la_{i\al}\al = \ga}}
\frac{\ga!}{\la!} \prod_{\substack{\al\in \mathbb N^n\\ \al>0}} 
\left(\frac1{\al!}\right)^{\sum_i\la_{i\al}} 
\left(\p^{\sum_\al \la_{\al}} f\right)(g(x))
\prod_{i,\al>0}(\p^\al g_i(x))^{\la_{i\al}} 
\end{align*}\end{proclaim}
The one dimensional version is due to Fa\`a di Bruno \cite{13}, the only
beatified mathematician. 

\begin{demo}{Proof} We compose the Taylor expansions of 
\begin{align*}
f(g(x)+h):&\quad j^\infty_{g(x)}f(h) 
  = \sum_{m\ge 0} \frac{f^{(m)}(g(x))}{m!} h^m, 
\\
g(x+t):&\quad j^\infty_x g(t) 
  = g(x) + \sum_{n\ge 1} \frac{g^{(n)}(x)}{n!}t^n, 
\\
f(g(x+t)):&\quad j^\infty_x(f\o g)(t) = \sum_{m\ge 0} \frac{f^{(m)}(g(x))}{m!} 
  \left(\sum_{n\ge 1} \frac{g^{(n)}(x)}{n!}t^n \right)^m
\\&\qquad
= \sum_{m\ge 0} \frac{f^{(m)}(g(x))}{m!} 
  \sum_{\al_1,\dots,\al_m>0}\left(\prod_{i=1}^m
  \frac{g^{(\al_i)}(x)}{\al_i!}\right)\;t^{\al_1+\dots+\al_m}. 
\end{align*}
Or we use the multinomial expansion
\begin{displaymath}
\Bigl(\sum_{j=1}^q a_j \Bigr)^m = \sum_{\substack{\la_1,\dots,\la_q\in\mathbb N_{\ge0} \\
\la_1+\dots+\la_q=m}} \frac{m!}{\la_1!\dots\la_q!}a_1^{\la_1}\dots
a_q^{\la_q}
\end{displaymath}
to get 
\begin{displaymath}
j^\infty_x(f\o g)(t) = \sum_{m\ge 0} \frac{f^{(m)}(g(x))}{m!} 
  \sum_{\substack{\la=(\la_n)\in \mathbb N_{\ge 0}^{\mathbb N_{>0}}\\\sum_n\la_n=m}}
  \frac{m!}{\la!}\left(\prod_{n>0}
  \left(\frac{g^{(n)}(x)}{n!}\right)^{\la_n}\right)\;t^{\sum_n \la_n n} 
\end{displaymath}
where $\la!=\la_1!\,\la_2!\,\dots$; most of the $\la_i$ are 0. 
The multidimensional formula just uses more indices.
\qed\end{demo}

\subsection*{\nmb.{6.2} }
The space $\mathcal{S}(\mathbb R)$ of all rapidly decreasing smooth 
functions $f$ for which 
$x\mapsto(1+|x|^2)^k\p_x^n f(x)$ is bounded for all  
$k\in \mathbb N$ and all $n\in \mathbb N_{\ge0}$, with the 
locally convex topology described by these conditions, is a nuclear 
Fr\'echet space. The dual space 
$\mathcal{S}'(\mathbb R)$ is the space of tempered distributions.

$\mathcal{S}(\mathbb R)$ is a commutative algebra under pointwise multiplication
and con\-vo\-lu\-tion $(u*v)(x)=\int u(x-y)v(y)dy$. The Fourier transform 
\begin{displaymath}
\mathcal{F}(u)(\xi)=\hat u(\xi)=\int e^{-ix\xi}u(x)dx,\qquad
\mathcal{F}\i(a)(x)=\frac1{2\pi}\int e^{ix\xi}a(\xi)d\xi,
\end{displaymath}
is an isomorphism of $\mathcal{S}(\mathbb R)$ and also of $L^2(\mathbb R)$ and has the
following further properties:  
\begin{align*}
\widehat {\p_x u}(\xi) &= -i\xi\cdot\hat u(\xi),\quad
\widehat {x\cdot u}(\xi) = -i\p_\xi \hat u(\xi),\quad
\\
\widehat {u(x - a)}(\xi) &= e^{ia\xi}\hat u (\xi), \quad
\widehat {e^{iax}u(x)}(\xi) = e^{ia\xi}\hat u (\xi), \quad
\\
\widehat {u(ax)}(\xi) &= \tfrac1{|a|}\hat u(\tfrac{\xi}{a}),\quad
\widehat {u(-x)}(\xi) =  \hat u(-\xi),
\\
\widehat {u\cdot v} &= \hat u * \hat v, \qquad 
\widehat {u * v} = \hat u \cdot \hat v .
\end{align*}
In particular, for any polynomial $P$ with constant coefficients we have 
\begin{displaymath}
\mathcal{F}(P(-i\p_x)u)(\xi) = P(\xi)\hat u(\xi). 
\end{displaymath}
$\mathcal{S}(\mathbb R)$ satisfies the uniform $\mathcal{V}$-boundedness principle for
every point separating set $\mathcal{V}$ of bounded linear functionals by
\cite{23},~5.24, since it is a Fr\'echet space; 
in particular for the set of all point evaluations
$\{\on{ev}_x:\mathcal{S}(\mathbb R)\to \mathbb R,\quad x\in \mathbb R\}$. 
Thus a linear mapping  
$\ell:E\to\mathcal{S}(\mathbb R)$ is bounded (smooth) if and only if $\on{ev}_x\o f$ 
is bounded for each $x\in \mathbb R$.

\begin{proclaim}{\nmb.{6.3}. Lemma} The space $C^\infty(\mathbb R,\mathcal{S}(\mathbb R))$ of smooth
curves in $\mathcal{S}(\mathbb R)$ consists of all functions 
$f\in C^\infty(\mathbb R^2,\mathbb R)$ satisfying the following 
property:
\begin{enumerate}
\item[$\bullet$]
For all $n,m\in \mathbb N_{\ge0}$ and each $t\in \mathbb R$ the expression
$(1+|x|^2)^k\p_t^{n}\p_x^{m}f(t,x)$  is uniformly bounded in $x$, locally
in $t$.
\end{enumerate} 
\end{proclaim}

\begin{demo}{Proof}
We use \nmb!{A.3} for the set $\{\on{ev}_x: x\in \mathbb R\}$ of point
evaluations in $\mathcal{S}'(\mathbb R)$. Note that $\mathcal{S}(\mathbb R)$ is reflexive. 
Here $c^k(t)=\p_t^kf(t,\quad)$.
\qed\end{demo}

\subsection*{\nmb.{6.4}. Diffeomorphisms which decrease rapidly to the identity
}
Any orientation preserving diffeomorphism $\mathbb R\to \mathbb R$ can be written
as $\on{Id}+f$ for $f$ a smooth function with $f'(x)>-1$ for all 
$x\in \mathbb R$.
Let us denote by $\on{Diff}_{\mathcal{S}}(\mathbb R)_0$ the space of all
diffeomorphisms 
$\on{Id}+f:\mathbb R\to\mathbb R$ (so $f'(x)>-1$ for all $x\in \mathbb R$) for 
$f\in \mathcal{S}(\mathbb R)$.

\begin{proclaim}{Theorem}
$\on{Diff}_{\mathcal{S}}(\mathbb R)_0$ is a regular Lie group.
%and a strong ILH group.
\end{proclaim}

\begin{demo}{Proof}
Let us first check that $\on{Diff}_{\mathcal{S}}(\mathbb R)_0$ is closed under
multiplication.
We have  
\begin{equation}
((\on{Id}+f)\o(\on{Id}+g))(x)=x+g(x)+f(x+g(x)),
\tag{\nmb:{1}}\end{equation}
and $x\mapsto f(x+g(x))$ is
in $\mathcal{S}(\mathbb R)$ by the Fa\`a di Bruno formula \nmb!{6.1} 
and the following estimate:
\begin{equation}
f^{(m)}(x+g(x))=O\Bigl(\frac{1}{(1+|x+g(x)|^2)^{k}}\Bigr)
  =O\Bigl(\frac{1}{(1+|x|^2)^{k}}\Bigr)
\tag{\nmb:{2}}\end{equation}
which holds since $g(x)\to 0$ for $|x|\to \infty$ and thus
\begin{displaymath}
\frac{1+|x|^2}{1+|x+g(x)|^2}\quad\text{  is globally bounded.}
\end{displaymath}

Let us check next that multiplication is smooth.
Suppose that the curves $t\mapsto \on{Id}+f(t,\quad),\on{Id}+g(t,\quad)$
are in
$C^\infty(\mathbb R,\on{Diff}_{\mathcal{S}}(\mathbb R)_0)$ which means
that the functions $f,g\in C^\infty(\mathbb R^2,\mathbb R)$ satisfy the conditions of
lemma \nmb!{6.2}. Then
\begin{displaymath}
(1+|x|^2)^k\p_t^n\p_x^mf(t,x+g(t,x))
\end{displaymath}
is bounded in $x\in \mathbb R$, locally in $t$, 
by the 2-dimensional Fa\'a di Bruno formula
\nmb!{6.1} and the more elaborate version of estimate {\nmb|{2}}
\begin{equation}
(\p^{(n,m)}f)(t,x+g(t,x))=O\Bigl(\frac{1}{(1+|x+g(t,x)|^2)^{k}}\Bigr)
  =O\Bigl(\frac{1}{(1+|x|^2)^{k}}\Bigr)
\tag{\nmb:{3}}\end{equation}
which follows from \nmb!{6.3} for $f$ and $g$.
Thus the multiplication respects smooth curves and is
smooth. 

To check that the inverse $(\on{Id}+g)\i$ is again an element in 
$\on{Diff}_{\mathcal{S}}(\mathbb R)_0$ for $g\in\mathcal{S}(\mathbb R)$, we write
$(\on{Id}+g)\i=\on{Id}+f$ and we have to check that $f\in \mathcal{S}(\mathbb R)$.  
\begin{align}
(\on{Id}+f)\o(\on{Id}+g)=\on{Id} &\implies x+g(x)+f(x+g(x))=x
\notag\\&
\implies x\mapsto f(x+g(x))=-g(x)\text{  is in }\mathcal{S}(\mathbb R).
\tag{\nmb:{4}}\end{align}
Now consider 
\begin{align*}
\p_x(f(x+g(x))) &= f'(x+g(x))(1+g'(x))
\\
\p_x^2(f(x+g(x))) &= f''(x+g(x))(1+g'(x))^2 + f'(x+g(x))g''(x)
\\
\p_x^3(f(x+g(x))) &= f^{(3)}(x+g(x))(1+g'(x))^3 +
\tag{\nmb:{5}}\\&\quad        
  +3f''(x+g(x))(1+g'(x))g''(x) 
  + f'(x+g(x))g^{(3)}(x)
\\
\p_x^m(f(x+g(x))) &= f^{(m)}(x+g(x))(1+g'(x))^m +
\\&\quad
+\sum_{k=1}^{m-1}f^{(m-k)}(x+g(x))a_{mk}(x),
%\\&
%= f^{(n)}(x+g(x))+
%\sum_{k=1}^{n}f^{(n-k)}(x+g(x))a_{nk}(x),
\end{align*}
where $a_{nk}\in \mathcal{S}(\mathbb R)$ for $n\ge k\ge 1$.
%Thus 
%\begin{displaymath}
%%\begin{pmatrix}  \p_x(f(x+g(x))) \\ \p_x^2(f(x+g(x))) \\ 
%          \p_x^3(f(x+g(x))) \\ \dots \end{pmatrix} =
%\begin{pmatrix} 1+a_{11}(x) & 0           & 0           & \dots \\
%         a_{21}(x)      & 1+a_{22}(x) & 0           & \dots    \\
%         a_{31}(x)      & a_{32}(x)   & 1+a_{33}(x) & \dots    \\
%                                \dots & & & \end{pmatrix}
%\begin{pmatrix} f'(x+g(x)) \\ f''(x+g(x)) \\ f^{(3)}(x+g(x)) \\ dots \end{pmatrix}
%\end{displaymath}
We have $1+g'(x)\ge \ep>0$ thus $\frac1{1+g'(x)}$ is bounded and its
derivative is in $\mathcal{S}(\mathbb R)$. Hence we can conclude 
that $(1+|x|^2)^k f^{(n)}(x+g(x))$ is bounded for each $k$. Since 
$(1+|x+g(x)|^2)^k=O(1+|x|^2)$ we conclude that
$(1+|x+g(x)|^2)^kf^{(n)}(x+g(x))$ is bounded for all $k$ and $n$. Inserting 
$y=x+g(x)$ it follows that $f\in\mathcal{S}(\mathbb R)$. 
Thus inversion maps $\on{Diff}_{\mathcal{S}}(\mathbb R)$ into itself. 

Let us check
that inversion is also smooth. So we assume that $g(t,x)$ is a smooth curve in
$\mathcal{S}(\mathbb R)$, satisfies \nmb!{6.3}, and we have to check that then
$f$ does the same. 
Retracing our considerations we see from {\nmb|{4}} that 
$f(t,x+g(t,x))=-g(t,x)$ satisfies
\nmb!{6.3} as a function of $t,x$, and we claim that $f$ then does the
same. Applying $\p_t^n$ to the equations in {\nmb|{5}} we get 
\begin{align*}
\p_t^n\p_x^m(f(t,x+g(t,x))) &= (\p^{(n,m)}f)(t,x+g(t,x))(1+\p_xg(t,x))^m 
+\\&\quad 
+\sum_{\substack{k_1\le n\\k_2\le m+n}}(\p^{(k_1,k_2)}f)(t,x+g(t,x))a_{k_1,k_2}(t,x),
\end{align*}
where $a_{k_1,k_2}(t,x)=O(\frac1{(1+|x|^2)^k})$ uniformly in $x$ and
locally in $t$. 
Again $1+\p_xg(t,x)\ge\ep>0$, locally in $t$ and uniformly in $x$, thus the
function $\frac1{1+\p_xg(t,x)}$ is bounded with any derivative in $\mathcal{}
S(\mathbb R)$ with respect to $x$. Thus we can conclude $f$ satisfies
\nmb!{6.3}. So the inversion is smooth and $\on{Diff}_S(\mathbb R)$ is a Lie
group.

We claim that $\on{Diff}_S(\mathbb R)$ is also a regular Lie group. So let
$t\mapsto X(t,\quad)$ be a smooth curve in the Lie algebra 
$\mathcal{S}(\mathbb R)\p_x$,i.e., $X$ satisfies \nmb!{6.3}. 
The evolution of this time dependent vector field is 
the function given by the ODE
\begin{align}
&\on{Evol}(X)(t,x) = x+f(t,x),
\notag\\
&\begin{cases} \p_t (x+f(t,x)) =f_t(t,x)= X(t,x+f(t,x)),\\
f(0,x)=0. \end{cases}
\tag{\nmb:{6}}\end{align}
We have to show that $f$ satisfies \nmb!{6.3}. 
For $0\le t\le C$ we consider  
\begin{equation}
|f(t,x)|\le \int_0^t |f_t(s,x)|ds=\int_0^t|X(s,x+f(s,x))|\,ds.
\tag{\nmb:{7}}\end{equation}
Since $X(t,x)$ is uniformly bounded in $x$, locally in $t$, the same is
true for $f(t,x)$ by {\nmb|{7}}. But then we may insert 
$X(s,x+f(s,x))=O(\frac1{(1+|x+f(s,x)|^2)^k})=O(\frac1{(1+|x|^2)^k})$ into 
{\nmb|{7}} and can conclude that $f(t,x)=O(\frac1{(1+|x|^2)^k})$
globally in $x$, locally in $t$, for each $k$.
For $\p_t^n\p_x^m f(t,x)$ we differentiate equation {\nmb|{6}} and
arrive at a system of ODE's with functions in $\mathcal{S}(\mathbb R)$ which
we can estimate in the same way.
\qed\end{demo}

\subsection*{\nmb.{6.5}. Sobolev spaces and $HC^n$-spaces}
The differential operator 
\begin{displaymath}
A_k=P_k(-i\p_x)=\sum_{i=0}^k(-1)^i\p_x^{2i},\qquad
P(\xi)=\sum_{i=0}^k\xi^{2i},
\end{displaymath}
will play an important role later on. 
We consider the \idx{\it Sobolev spaces}, namely the Hilbert spaces
\begin{displaymath}
H^n(\mathbb R) =\{f\in \mathcal{S}'(\mathbb R): f, f', f^{(2)},\dots f^{(n)}\in
L^2(\mathbb R)\}.
\end{displaymath}
In terms of the Fourier transform $\hat f$ we have, by the properties
listed in \nmb!{6.2}:
\begin{align*}
f\in H^n &\iff (1+|\xi|)^n\hat f(\xi)\in L^2 \iff (1+|\xi|^2)^{n/2}\hat f(\xi))\in L^2
\\&
\iff (1+|\xi|)^{n-2k}P_k(\xi)\hat f(\xi)\in L^2  \iff A_k(f)\in H^{n-2k}.
\end{align*}
We shall use the norm 
$$
\|f\|_{H^n} := \|\hat f(\xi)(1+|\xi|)^n\|_{L^2}
$$
on $H^n(\mathbb R)$.
Moreover, for $0<\al\le 1$ we consider the Banach space 
\begin{displaymath}
C^{0,\al}_b(\mathbb R)=\{f\in C^0(\mathbb R):
  \sup_{x\in \mathbb R}|f(x)| 
  + \sup_{x\ne y\in \mathbb R}\frac{|f(x)-f(y)|}{|x-y|^\al}<\infty\}
\end{displaymath}
of bounded H\"older continuous functions on $\mathbb R$, and the Banach spaces
\begin{displaymath}
C^{n,\al}_b(\mathbb R) = \{f\in C^n(\mathbb R): 
  f, f', \dots, f^{(n-1)}\text{ bounded, and }
  f^{(n)}\in C^{0,\al}_b(\mathbb R) \}.
\end{displaymath}
Finally we shall consider the space
$$
HC^n(\mathbb R) = H^n(\mathbb R)\cap C^n_b(\mathbb R),\quad \|f\|_{HC^c} =
\|f\|_{H^n} + \|f\|_{C^n_b}.
$$

\begin{proclaim}{\nmb.{6.6}. Lemma}
Consider the differential operator $A_k=\sum_{i=0}^k(-1)^i\p_x^{2i}$.
\begin{enumerate}
\item[(\nmb:{1})] 
   $A_k:\mathcal{S}(\mathbb R)\to \mathcal{S}(\mathbb R)$ is a
   linear isomorphism of the Fr\'echet space of rapidly decreasing smooth
   functions.
\item[(\nmb:{2})] $A_k:H^{n+2k}(S^1)\to H^n(S^1)$ is a
   linear isomorphism of Hilbert spaces for each $n\in \mathbb Z$, where 
   $H^n(S^1)=\{f\in L^2(S^1): A_n(f)\in L^2(S^1)\}$. Note that
   $H^n(S^1)\subseteq C^k(S^1)$ if $n>k+1/2$ (Sobolev inequality).
\item[(\nmb:{3})] 
   $A_k:C^\infty(S^1)\to C^\infty(S^1)$ is a
   linear isomorphism.
\item[(\nmb:{4})] 
   $A_k:HC^{n+2k}(\mathbb R)\to HC^n(\mathbb R)$ is a linear isomorphism of
   Banach spaces for each $n\ge 0$.
\end{enumerate}
\end{proclaim}

\begin{demo}{Proof} Without loss we may consider complex-valued functions. 

\nmb|{1}
Let $\mathcal{F}:C^\infty(S^1)\to s(\mathbb Z)$ be the Fourier transform which is
an isomorphism on the
space of rapidly decreasing sequences. 
Since $\mathcal{F}(f_{xx})(n)=-(2\pi n)^2\mathcal{F}(f)(n)$ we have
$\mathcal{F}\o A_k\o \mathcal{F}\i: 
(c_n)\mapsto ((1+(2\pi n)^2+\dots+ (2\pi n)^{2k})\, c_n)$
which is a linear bibounded isomorphism.

\nmb|{2}
This is obvious from the definition.

\nmb|{3} can be proved similarly to \nmb|{1}, 
using that the Fourier series expansion is an isomorphism between 
$C^\infty(S^1)$ and the space $\mathcal s$ of rapidly decreasing sequences.

\nmb|{4} follows from \nmb|{2}.
%and the fact that 
%$L^2(\mathbb R)\cap C^0_b(\mathbb R)\subset C^0_0(\mathbb R):= 
%\{f\in C^0(\mathbb R): \lim_{x\to -\infty}f(x)=0\}$.
%Introducing spaces $C^n_{-\infty}(\mathbb R)$ of $C^n$-functions where each derivative
%vanishes at $-\infty$ we have 
%$HC^n(\mathbb R)=H^n(\mathbb R)\cap C^n_{-\infty}(\mathbb R)$.
%But obviously $A_k:C^{n+k}_{-\infty}(\mathbb R)\to C^n_{-\infty}(\mathbb R)$   is a
%linear isomorphism since it is linear ODO.
\qed\end{demo}

\begin{proclaim}{\nmb.{6.7}. Sobolev inequality} 
We have bounded linear embeddings ($0<\al\le 1$):
\begin{align*}
&H^n(\mathbb R) \subset C^k_b(\mathbb R) \text{  if } n>k+\tfrac12,
\\
&H^n(\mathbb R) \subset C^{k,\al}_b(\mathbb R) \text{  if }
n>k+\tfrac12+\al.
\end{align*}
\end{proclaim}

\begin{demo}{Proof} 
Since $\p_x^{k}:H^n(\mathbb R)\to H^{n-k}(\mathbb R)$ is bounded we 
may assume that $k=0$. So let $n>\frac12$. Then we use the Cauchy-Schwartz
inequality:
\begin{align*}
2\pi|u(x)|&=\left|\int e^{ix\xi}\hat u(\xi)\,d\xi\right| 
  \le \int |\hat u(\xi)|\,d\xi
= \int |\hat u(\xi)|(1+|\xi|)^{n}\frac1{(1+|\xi|)^n}\,d\xi
\\&
\le\left(\int |\hat u(\xi)|^2(1+|\xi|)^{2n}\,d\xi\right)^{\frac12}
\left(\int\frac1{(1+|\xi|)^{2n}}\,d\xi\right)^{\frac12}
=C\|u\|_{H^n}
\end{align*}
where 
$$
C= \left(\int\frac1{(1+|\xi|)^{2n}}\,d\xi\right)^{\frac12} < \infty
$$
depends only on $n>\frac12$.
For the second assertion we use 
$x>y$ and 
\begin{align*}
e^{ix\xi}-e^{iy\xi}&=(x-y)\int_0^1 i\xi e^{i(y+t(x-y))\xi}dt,
\\
\left|e^{ix\xi}-e^{iy\xi}\right|&\le |x-y|.|\xi|
\end{align*}
to obtain
\begin{align*}
&2\pi\left|\frac{u(x)-u(y)}{(x-y)^\al}\right|
\le \int \left|\frac{e^{ix\xi}-e^{iy\xi}}{x-y}\right|^\al.
  \left|e^{ix\xi}-e^{iy\xi}\right|^{1-\al} |\hat u(\xi)|\,d\xi
\\&
\le 2\int |\hat u(\xi)|(1+|\xi|)^{n}\frac{|\xi|^\al}{(1+|\xi|)^n}\,d\xi
\\&
\le 2\left(\int |\hat u(\xi)|^2(1+|\xi|)^{2n}d\xi\right)^{\frac12}
\left(\int\frac{|\xi|^{2\al}}{(1+|\xi|)^{2n}}\,d\xi\right)^{\frac12}
=C_1\|u\|_{H^n}
\end{align*}
where $C_1$ depends only on $n-\al>\frac12$.
\qed\end{demo}

\begin{proclaim}{\nmb.{6.8}. Banach algebra property}
If $n>\frac12$ then pointwise multiplication 
$\mathcal S(\mathbb R)\x \mathcal S(\mathbb R)\to\mathcal S(\mathbb R)$ 
extends to a bounded bilinear mapping 
$H^n(\mathbb R)\x H^n(\mathbb R)\to H^n(\mathbb R)$.

For $n\ge 0$ multiplication $HC^n(\mathbb R)\x HC^n(\mathbb R)\to HC^n(\mathbb R)$
is bounded bilinear.
\end{proclaim}

See \cite{Eich3} for the most general version of this on open Riemannian
manifolds with bounded geometry. 

\begin{demo}{Proof}
For $f,g\in H^n(\mathbb R)$ we have to show that for $0\le k\le n$ we have
$$
(f.g)^{(k)} = \sum_{l=0}^k\binom{k}{l}f^{(l)}.g^{(k-l)}\quad\in L^2(\mathbb R)  
$$
with norm bounded by a constant times $\|f\|_{H^n}.\|g\|_{H^n}$. If $l<n$
then $f^{(l)}\in C^0_b(\mathbb R)$ by the Sobolev inequality 
and $g^{(k-l)}\in H^l\subset L^2$ so the
product is in $L^2$ with the required bound on the norm. If $l=0$ we
exchange $f$ and $g$.

In the case of $HC^n$, the $L^2$-norm of each product in the sum is bounded
by the sup-norm of the first factor times the $L^2$-norm of the second one. 
And the sup-norm is clearly submultiplicative. 
\qed\end{demo}

\begin{proclaim}{\nmb.{6.9}. Differentiability of composition}
If $n\ge 0$ then composition
$\mathcal S(\mathbb R)\x \mathcal S(\mathbb R)\to\mathcal S(\mathbb R)$ 
extends to a weakly $C^k$-mapping
$HC^{n+k}(\mathbb R)\x (\on{Id}_{\mathbb R}+HC^n(\mathbb R))\to HC^n(\mathbb R)$. 
\end{proclaim}

A mapping $f:E\to F$ is weakly $C^1$ for Banach spaces $E,F$ if $df: E\x
E\to F$ exists and is continuous. We call it strongly $C^1$ if $df:E\to
L(E,F)$ is continuous for the operator norm on the image space. 
Similarly for $C^k$. 
Since I could not find a convincing proof of this result for the spaces
$H^n$ under the assumtion $n>\frac12$, 
I decided to use the spaces $HC^n(\mathbb R)$. This also inproves on the
degree $n$ which we need. 

\begin{demo}{Proof}
We consider the Taylor expansion 
\begin{align*}
f(x+g(x)) &= \sum_{p=0}^k \frac1{p!}f^{(p)}(x).g(x)^p +
\\&\quad
+\int_0^{1}\frac{(1-t)^{k-1}}{(k-1)!}\big(f^{(k)}(x+tg(x))-f^{(k)}(x)\big)\,dt\,.g(x)^k
\end{align*}
For fixed $f$ this is weakly $C^k$ in $g$ by invoking the Banach
algebra property and by estimating the integral in the remainder term. 
We have to show that the integrand is continuous at $(f^{(k)},g=0)$ as a
mapping $H^n\x H^n\to H^n$.
The integral from 0 to 1 does not disturb this so we disregard it.  
By \nmb!{6.1} we have
\begin{multline*}
\p_x^p(f^{(k)}(x+g(x))-f^{(k)}(x))=
\\ 
= p!\sum_{m=0}^p\frac{f^{(k+m)}(x+g(x))}{m!}
\sum_{\substack{\al_1,\dots,\al_m>0\\ \al_1+\dots+\al_m =p}} 
\frac{\p_x^{\al_1}(x+g(x))}{\al_1!}\dots
  \frac{\p_x^{\al_m}(x+g(x))}{\al_m!}
\end{multline*}
The most dangerous term is the one for $p=n$. As soon as a derivative 
of $g$ of order $\ge 2$ is present, this is easily estimated. 
The most difficult term is 
$$
f^{(k+n)}(x+g(x)) - f^{(k+n)}(x)
$$
which should go to 0 in $L^2\cap C^0_b$ for fixed $f$ and for $g\to 0$ in
$HC^n$.
$f^{(k)}$ is continuous and in $L^2$. Off some big compact intervall it has
small $H^n$-norm and small sup-norm (the latter by the lemma of Riemann-Lebesque). 
On this compact
intervall $f^{(k)}$ is uniformly continuous and if we choose $\|g\|_{C^n}$ small
enough, $f^{(k)}(x+t g(x))-f^{(k)}(x)$ is uniformly small there, thus small
in the sup-norm, and also small in $L^2$ (which involves the length of the
compact intervall -- but we can still choose $g$ smaller).
\qed\end{demo}                  

The last result cannot be improved to strongly $C^k$ since we have:

\begin{proclaim}{\nmb.{6.10}. Attention}
Composition 
$HC^{n}(\mathbb R)\x(\on{Id}_{\mathbb R}+ HC^n(\mathbb R))\to HC^n(\mathbb R)$ 
is only continuous and not
Lipschitz in the first variable.
\end{proclaim}

\begin{demo}{Proof}
To see this, consider $(f,t)\mapsto f(\quad -t.g)$ for a given bump
function $g$ which equals 1 on a large intervall. For each $t>0$ we
consider a bump function $f$ with support in $(-\frac t2,\frac t2)$ 
with $\|f\|_{L^2}=1$. Then we have $\|f- f(\quad-t)\|_{L^2}=\sqrt{2}$ by
Pythagoras, and consequently 
$\|f- f(\quad-t.g)\|_{HC^n}\ge \|f-f(\quad-t)\|_{L^2}=\sqrt{2}$.
\qed\end{demo}

\subsection*{\nmb.{6.11}. The topological group $\on{Diff}(\mathbb R)$}
For $n\ge1$ we consider $f:\mathbb R\to \mathbb R$ of the form 
$f(x)=x+g(x)$ for $g\in HC^n$. 
Then $f$ is a $C^n$-diffeomorphism iff $g'(x)>-1$ for all $x$. 
%Consider
%$$
%y= f(x)= x+g(x) \quad\implies\quad  f\i(y) = x = y - g(f\i(x)) = y + h(y),
%\quad\text{  where }h(y)= -g(f\i(y)).
%$$
The inverse is also of the form $f\i(y)=y+h(y)$ for $h\in HC^n(\mathbb R)$ 
iff $g'(x)\ge -1+\ep$ for a constant $\ep$. Indeed, $h(y)=-g(f\i(y))$.
Let us call $\on{DiffHC}^n(\mathbb R)$ the group of all these diffeomorphsms.

\begin{proclaim}{Lemma}
Inversion $\on{DiffHC}^{n+k}(\mathbb R)\to \on{DiffHC}^n(\mathbb R)$ is weakly
$C^k$. 
\end{proclaim}

\begin{demo}{Proof}
As we saw above, $\on{DiffHC}^{n+k}(\mathbb R)$ is stable under inversion.
$(f,g)\mapsto f\o g$ is a weak $C^k$ submersion by \nmb!{6.9}.
So we can use the implicit function theorem for the equation $f\o
f\i=\on{Id}$.
\qed\end{demo}

\begin{proclaim}{\nmb.{6.12} Proposition}
For $n\ge 1$ and $a\in HC^n(\mathbb R)$, the mapping 
$HC^n(\mathbb R)\x \on{DiffHC}^n(\mathbb R)\to HC^{n-1}(\mathbb R)$ given by 
$(f,g)\mapsto (a\p_x(f\o g\i))\o g$
is continuous and Lipschitz in $f$.

For $n>k+\frac12$ and for each linear differential operator $D$ of order
$k$, the mapping
$HC^n(\mathbb R)\x \on{DiffHC}^n(\mathbb R)\to HC^{n-k}(\mathbb R)$ given by 
$(f,g)\mapsto (D(f\o g\i))\o g$
is continuous and Lipschitz in $f$.
\end{proclaim}

Here $\on{Diff}(\mathbb R)=\{\on{Id}_{\mathbb R}+h: \|h'\|_{C^0_b}>-1 \}$.

\begin{demo}{Proof}
We have 
$$
(a\p_x(f\o g\i))\o g = \Big(a.(f_x\o g\i)\frac1{g_x\o g\i}\Big)\o g
=(a\o g).f_x.\frac 1{g_x}
$$
which is Lipschitz by the results above.
\qed\end{demo}

\begin{proclaim}{\nmb.{6.13} Proposition}
For the operator $A_k=\sum_{i=0}^k(-1)^i\p_x^{2i}$ and for $n\ge 2k$,
the mapping \newline
$(f,g)\mapsto (A_k\i(f\o g\i))\o g$ is Lipschitz 
$HC^n(\mathbb R)\x \on{DiffHC}^n(\mathbb R)\to HC^{n+2k}(\mathbb R)$.
\end{proclaim}

\begin{demo}{Proof}
The inverse of $A_k$ is given by the pseudo differential operator
$$
(A_k\i f)(x) 
= \int_{\mathbb R^2} e^{i(x-y)\xi}f(y) \frac 1{1+\xi^2+\quad+\xi^{2n}}d\xi\,dy
$$
Thus the mapping is given by 
\begin{align*}
(A_k\i(f\o g\i))(g(x)) &
= \int_{\mathbb R^2} e^{i(g(x)-y)\xi}f(g\i(y))
  \frac 1{1+\xi^2+\quad+\xi^{2n}}d\xi\,dy
\\&
= \int_{\mathbb R^2} e^{i(g(x)-g(z))\xi}f(z) 
  \frac{g'(z)}{1+\xi^2+\quad+\xi^{2n}}d\xi\,dz
\end{align*}
which is a genuine Fourier integral operator.
By the foregoing results this is visibly locally Lipschitz.
\qed\end{demo}

\section*{\nmb0{7}. The diffeomorphism group of $S^1$ or $\mathbb R$, 
and Burgers' hierarchy
} 

\subsection*{\nmb.{7.1}. Burgers' equation and its curvature } 
We consider the Lie groups $\on{Diff}_{\mathcal S}(\mathbb R)$ and  
$\on{Diff}(S^1)$ with Lie algebras $\X_{\mathcal S}(\mathbb R)$ and $\X(S^1)$ where the  
Lie bracket $[X,Y]=X'Y-XY'$ is the negative of the usual one. For the  
$L^2$-inner product $\ga(X,Y)=\langle X,Y\rangle_0=\int X(x)Y(x)\,dx$ 
integration by parts gives  
\begin{align*}
\langle [X,Y],Z\rangle_0 &= \int_\mathbb R (X'YZ-XY'Z)dx  
\\&
= \int_\mathbb R (2X'YZ+XYZ')dx = \langle Y,\ad(X)^\top Z\rangle,  
\end{align*}
which in turn gives rise to  
\begin{align} 
\ad(X)^\top Z &= 2X'Z+XZ', 
\tag{\nmb:{1}}\\
\al(X)Z=\ad(Z)^\top X &= 2Z'X+ZX' ,
\tag{\nmb:{2}}\\
(\ad(X)^\top+\ad(X)) Z &= 3X'Z ,
\tag{\nmb:{3}}\\
(\ad(X)^\top-\ad(X)) Z &= X'Z + 2XZ' = \al(X)Z . 
\tag{\nmb:{4}}  
\end{align}
Equation {\nmb|{4}} states that $-\frac12\al(X)$ 
is the skew-symmetrization of  
$\ad(X)$ with respect to the inner product  
$\langle\quad,\quad\rangle_0$. From the 
theory of symmetric spaces one  
then expects that $-\frac12\al$ is a 
Lie algebra homomorphism and  
indeed one can check that  
\begin{displaymath}
-\tfrac12\al([X,Y]) = \left[-\tfrac12\al(X),
                       -\tfrac12\al(Y)\right] 
\end{displaymath}
holds for any vector fields $X, Y$. From 
{\nmb|{1}} we get the geodesic equation, whose second part is 
Burgers' equation \cite{9}:  
\begin{displaymath}
\begin{cases} 
&g_t(t,x) = u(t,g(t,x))
\\
&u_t = -\ad(u)^\top u = -3u_xu
\end{cases} 
\tag{\nmb:{5}}
\end{displaymath}
Using the above relations and the general curvature formula  
\nmb!{3.4.2}, we get
\begin{align} 
\mathcal{R}(X,Y)Z &= -X''YZ +XY''Z -2X'YZ' +2XY'Z'  
\notag\\&
= -2[X,Y]Z' - [X,Y]'Z 
= -\al([X,Y])Z. 
\tag{\nmb:{6}}\end{align}
Sectional curvature is non-negative and unbounded:
\begin{align}
-G^0_a(R(X,Y)X,Y) &= \langle\al([X,Y])(X),Y\rangle 
= \langle\ad(X)^\top([X,Y]),Y\rangle  
\notag\\&
= \langle[X,Y],[X,Y]\rangle = \|[X,Y]\|^2, 
\notag\\
k(X\wedge Y)
&= - \frac{G^0_a(R(X,Y)X,Y)}{\|X\|^2\|Y\|^2-G^0_a(X,Y)^2}
\notag\\
&= \frac{\|[X,Y]\|^2}{\|X\|^2\|Y\|^2-\langle X,Y\rangle^2}\ge 0.
\tag{\nmb:{7}}\end{align}
Let us check invariance of the momentum mapping $\bar J$ from \nmb!{4.3}:
\begin{align}
\ga(\bar J(g,X),Y) &= \ga(\Ad(g)^\top X,Y) = \ga(X,\Ad(g)Y) 
=\int X((g'Y)\o g\i)dx
\notag\\& 
= \int X(g'\o g\i)(Y\o g\i)dx 
= \on{sign}(g')\int (X\o g)(g')^2 Y dx 
\notag\\& 
=\on{sign}(g')\ga((g')^2(X\o g),Y)
\notag\\
\bar J(g,X) &= \on{sign}(g_x).(g_x)^2(X\o g).    
\tag{\nmb:{8}}\end{align}
Along a geodesic $t\mapsto g(t,\quad)$, according to {\nmb|{5}} and
\nmb!{4.3}, the momentum 
\begin{equation}
\bar J(g,u=g_t\o g\i) = g_x^2g_t \quad \text{  is constant.}
\tag{\nmb:{9}}\end{equation}
This is what we found in \nmb!{1.3} by chance.
 
\subsection*{\nmb.{7.2}. Jacobi fields for Burgers' equation} 
A Jacobi field $y$ along a geodesic $g$ with velocity field $u$  
is a solution of the  
partial differential equation \nmb!{3.5.1}, which in our case  
becomes: 
\begin{align} 
y_{tt}&= [\ad(y)^\top+\ad(y),\ad(u)^\top]u  
     - \ad(u)^\top y_t -\al(u)y_t + \ad(u)y_t 
\tag{\nmb:{1}}\\&
= - 3u^2y_{xx} -4 uy_{tx} -2u_xy_t
\notag\\
u_t &= -3 u_xu. 
\notag\end{align}
If the geodesic equation has smooth solutions locally in time 
 it is to be expected that 
the space of all Jacobi fields exists and is isomorphic to the space 
of all initial data $(y(0),y_t(0))\in C^\infty(S^1,\mathbb R)^2$ or 
$C^\infty_c(\mathbb R,\mathbb R)^2$,   
respectively. 
The weak symplectic structure on it is given by \nmb!{3.7}: 
\begin{align} 
\om(y,z) &= \langle y, z_t-\tfrac12u_xz+2uz_x \rangle 
     - \langle y_t-\tfrac12u_xy+2uy_x, z \rangle 
\notag\\&
= \int_{S^1\text{ \!\!or }\mathbb R} (yz_t-y_tz+2u(yz_x-y_xz))\,dx. 
\tag{\nmb:{2}}\end{align}

\subsection*{\nmb.{7.3}. The Sobolev $H^k$-metric on $\on{Diff}(S^1)$ and
$\on{Diff}(\mathbb R)$ }
On the Lie algebras $\X_c(\mathbb R)$ and $\X(S^1)$ with 
Lie bracket $[X,Y]=X'Y-XY'$ we consider the $H^k$-inner product  
\begin{align}
\ga(X,Y)=\langle X,Y\rangle_k&=\sum_{i=0}^k\int (\p_x^iX)(\p_x^iY)\,dx
=\int A_k(X)(Y)\,dx 
\notag\\&
=\int X A_k(Y)\,dx,   
\quad\text{  where }\quad
A_k=\sum_{i=0}^k(-1)^i\p_x^{2i}
\tag{\nmb:{1}}\end{align}
is a linear isomorphism  
$\X_c(\mathbb R)\to \X_c(\mathbb R)$ or $\X(S^1)\to \X(S^1)$ whose inverse is a
pseudo differential operator.
$A_k$ is also a bounded linear isomorphism between the Sobolev spaces
$H^{l+2k}(S^1)\to H^l(S^1)$, see lemma \nmb!{6.5}.
On the real line we have to consider functions with fixed
support in some compact set $[-K,K]\subset \mathbb R$.

Integration by parts gives  
\begin{align*}
\langle [X,Y],Z\rangle_k &= \int_\mathbb R (X'Y-XY')A_k(Z)dx  
     = \int_\mathbb R (2X'YA_k(Z)+XYA_k(Z'))dx
\\& 
     = \int_\mathbb R YA_kA_k\i\bigl(2X'A_k(Z)+X A_k(Z')\bigr)dx 
     = \langle Y,\ad(X)^{\top,k},Z\rangle_k,  
\end{align*}
which in turn gives rise to  
\begin{align*} 
\ad(X)^{\top,k} Z &= A_k\i\bigl(2X'A_k(Z)+XA_k(Z')\bigr), 
\\
\al_k(X)Z &= \ad(Z)^{\top,k}(X)= A_k\i\bigl(2Z'A_k(X)+ZA_k(X')\bigr)
\tag{\nmb:{2}}
\end{align*}
Thus the geodesic equation is 
\begin{equation}
\begin{cases}
&g_t(t,x) = u(t,g(t,x)) 
\\
&{\begin{aligned}
u_t&=- \ad(u)^{\top,k} u = -A_k\i\bigl(2u_x A_k(u)+u A_k(u_x)\bigr)
\\&
= -A_k\i\Bigl(2u_x\sum_{i=0}^k(-1)^i\p_x^{2i}u
              +u\sum_{i=0}^k(-1)^i\p_x^{2i+1}u\bigr).
\end{aligned}}\end{cases}
\tag{\nmb:{3}}
\end{equation}
For $k=0$ the second part is Burgers' equation, and for $k=1$ it becomes
\begin{align*}
&u_t - u_{txx}=-3uu_x+2u_xu_{xx}+uu_{xxx}
\tag{\nmb:{4}}\\
\iff & u_t+uu_x +(1-\p_x^2)\i(u^2+\tfrac12 u_x^2)_x = 0
\end{align*}
which is the dispersionfree version of the {\it Camassa-Holm equation}, 
see \cite{10}, \cite{33}, \cite{24}. We met it already in \nmb!{5.6}, and
will meet the full equation in \nmb!{8.7}.
Let us check the invariant momentum mapping from \nmb!{4.3.2}:
\begin{align}
\ga(\bar J(g&,X),Y) = \langle \Ad(g)^\top X,Y\rangle_k 
= \langle X,\Ad(g)Y\rangle_k 
\notag\\& 
= \int A_k(X)(g'\o g\i)(Y\o g\i)dx 
= \on{sign}(g')\int (A_k(X)\o g)(g')^2 Y dx 
\notag\\& 
=\on{sign}(g')\Bigl\langle A_k\i\Bigl((g')^2(A_k(X)\o g)\Bigr),Y\Bigr\rangle_k
\notag\\
\bar J(g,X) &= \on{sign}(g_x).A_k\i\Bigl((g_x)^2(A_k(X)\o g)\Bigr).
\tag{\nmb:{5}}\end{align}
Along a geodesic $t\mapsto g(t,\quad)$, 
by {\nmb|{3}} and \nmb!{4.3}, the expressions 
\begin{equation*}
\on{sign}(g_x)\bar J(g,u=g_t\o g\i) = A_k\i\Bigl((g_x)^2(A_k(u)\o g)\Bigr)
\tag{\nmb:{6}}\end{equation*}
and thus also $(g_x)^2(A_k(u)\o g)$ are constant in $t$.

\begin{proclaim}{\nmb.{7.4}. Theorem}  
Let $k\ge 1$. 
There exists a $HC^{2k+1}$-open neighborhood $V$ of $(\on{Id},0)$ in
$\on{Diff}(S^1)\x \X(S^1)$ such that for each 
$(g_0,u_0)\in V$ there exists a unique $C^3$ geodesic
$g\in C^3((-2,2),\on{Diff}(S^1))$ for the right invariant $H^k$
Riemann metric, starting at $g(0)=g_0$ in the direction $g_t(0)=u_0\o g_0\in
T_{g_0}\on{Diff}(S^1)$. Moreover, the solution depends $C^1$ on the initial
data $(g_0,u_0)\in V$.

The same result holds if we replace $\on{Diff}(S^1)$ by
$\on{Diff}_{\mathcal S(\mathbb R)}$ and $\X(S^1)$ by $\X_{\mathcal S}(\mathbb
R)=\mathcal S(\mathbb R)\p_x$.
\end{proclaim}

This result is stated in \cite{11}, and also this proof follows essentially
\cite{11}. But there is a mistake in \cite{11},~p~795, where the authors
assume that composition and inversion on $H^n(S^1)$ are smooth. This is wrong. 
One needs to use \nmb!{6.12} and \nmb!{6.13}.
The mistake was corrected in \cite{12}, for the more general 
case of the Virasoro group.

In the following proof,
$\on{Diff}$, $\X$, $\on{DiffHC}^n$, $HC^n$ should stand for either 
$\on{Diff}(S^1)$, $\X(S^1)$, $\on{DiffHC}^n(S^1)$, $HC^n(S^1)$ or for 
$\on{Diff}_{\mathcal S}(\mathbb R)$, $\X_{\mathcal S}(\mathbb R)$,  
$\on{DiffHC}^n(\mathbb R)$, $HC^n(\mathbb R)$, respectively.

\begin{demo}{Proof}
For $u\in HC^n$, $n\ge 2k+1$, we have
\begin{align*}
A_k(uu_x)&=\sum_{i=0}^k(-1)^i\p_x^{2i}(uu_x)
=\sum_{i=0}^k(-1)^i\sum_{j=0}^{2i}\tbinom{2i}{j}(\p_x^{j}u)(\p_x^{2i-j+1}u)
\\&
=uA_k(u_x)+\sum_{i=0}^k(-1)^i\sum_{j=1}^{2i}\tbinom{2i}{j}(\p_x^{j}u)(\p_x^{2i-j+1}u)
\\&
=:u\,A_k(u_x)+B_k(u),
\end{align*}
where $B_k:HC^n\to HC^{n-2k}$ is a bounded quadratic operator.
Recall that we have to solve  
\begin{align*}
u_t&=- \ad(u)^{\top,k} u = -A_k\i\bigl(2u_x A_k(u)+u A_k(u_x)\bigr)
\\&
= -A_k\i\bigl(2u_x A_k(u)+ A_k(uu_x)-B_k(u)\bigr)
\\&
= -uu_x -A_k\i\bigl(2u_x A_k(u)-B_k(u)\bigr)
\\&
=: -uu_x +A_k\i C_k(u),
\end{align*}
where $C_k:HC^n\to HC^{n-2k}$ is a bounded quadratic operator, and
where $u=g_t\o g\i\in \X$.
Note that  
\begin{align*}
C_k(u) &= -2u_x A_k(u)+B_k(u)
\\&
= -2u_x A_k(u)
  +\sum_{i=0}^k(-1)^i\sum_{j=1}^{2i}\tbinom{2i}{j}(\p_x^{j}u)(\p_x^{2i-j+1}u).
\end{align*}
We put
\begin{align}
&\begin{cases}
&g_t=:v= u\o g
\\
&{\begin{aligned}
v_t&= u_t\o g + (u_x\o g)g_t = u_t\o g + (uu_x)\o g
= A_k\i C_k(u)\o g  
\\&
= A_k\i C_k(v\o g\i)\o g  
=: \on{pr_2}(D_k\o E_k)(g,v), \quad \text{  where }
\end{aligned}}
\end{cases}
\tag{\nmb:{1}}\\
&E_k(g,v)=(g,C_k(v\o g\i)\o g), \qquad D_k(g,v)= (g,A_k\i(v\o g\i)\o g).
\notag\end{align}
Now consider the topological group and Banach manifold 
$\on{DiffHC}^n$ described in \nmb!{6.11}.

\noindent\therosteritem{\nmb:{2}}
{\it Claim.} The mapping 
$D_k:\on{DiffHC}^n\x HC^{n-2k}\to \on{DiffHC}^n\x HC^n$ 
is strongly $C^1$.

First we check that all directional derivatives exist 
and are in the right spaces.   

For $w\in HC^n$ we have
\begin{align*}
&\p_s|_0 (u\o (g+sw)) = (u_x\o g)w
\\&
\p_s|_0 (g+sw)\i = -\frac{w\o g\i}{g_x\o g\i}
\\&
\p_s|_0 \on{pr}_2 D_k(g+sw,v) =
\\&
= \p_s|_0 A_k\i(v\o g\i)\o(g+sw)
+ \p_s|_0 (A_k\i(v\o (g+sw)\i))\o g
\\&
= ((\p_xA_k\i(v\o g\i))\o g)\,w
-  (A_k\i((v_x\o g\i)\tfrac{w\o g\i}{g_x\o g\i}))\o g
\\&
= (A_k\i(v\o g\i)_x.(w\o g\i))\o g
-  (A_k\i((v\o g\i)_x (w\o g\i)))\o g.
\end{align*}
Therefore, 
\begin{align*}
&A_k((\p_s|_0 \on{pr}_2 D_k(g+sw,v))\o g\i) =
\\&
= A_k(A_k\i(v\o g\i)_x.(w\o g\i))
-  (v\o g\i)_x (w\o g\i)
\\&
= (v\o g\i)_x.(w\o g\i) 
+\sum_{i=0}^k \sum_{j=0}^{2i-1}\binom{2i}{j}
\p_x^{j+1}A_k\i(v\o g\i).\p_x^{2k-j}(w\o g\i)
\\&\quad
-(v\o g\i)_x (w\o g\i) 
\in HC^{n-2k}.
\end{align*}
By \nmb!{6.12} and \nmb!{6.13} this is locally Lipschitz jointly in 
$v,g,w$. Moreover we have 
$\p_s|_0 \on{pr}_2 D_k(g+sw,v)\in HC^n$, and   
$D_k$ is linear in $v$.
Thus $D_k$ is strongly $C^1$.

\noindent\therosteritem{\nmb:{3}}
{\it Claim.}  The mapping 
$E_k:\on{DiffHC}^n\x HC^{n}\to \on{DiffHC}^n\x HC^{n-2k}$ is
strongly $C^1$.
This can be proved similarly, again using \nmb!{6.12} and \nmb!{6.13}.

By the two claims equation {\nmb|{1}} can be viewed as the flow equation of
a $C^1$-vector field on the Hilbert manifold 
$\on{DiffHC}^n\x HC^{n}$. 
Here an existence and uniqueness theorem holds.
Since $v=0$ is a
stationary point, there exist an open neighborhood $W_n$ of $(\on{Id},0)$
in $\on{DiffHC}^n\x HC^{n}$ such that for each initial point
$(g_0,v_0)\in W_n$ equation {\nmb|{1}} has a unique solution 
$\on{Fl}^n_t(g_0,v_0)=(g(t),v(t))$
defined and $C^2$ in $t\in (-2,2)$. Note that $v(t)=g_t(t)$, thus $g(t)$ is
even $C^3$ in $t$. Moreover, the solution depends $C^1$ on the initial
data.

We start with the neighborhood 
\begin{displaymath}
W_{2k+1}\subset \on{DiffHC}^{2k+1}\x HC^{2k+1} 
  \supset \on{DiffHC}^n\x HC^{n}\quad \text{  for }n\ge 2k+1 
\end{displaymath}
and consider the neighborhood $V_n:=W_{2k+1}\cap\on{DiffHC}^n\x HC^{n}$ 
of $(\on{Id},0)$

\noindent\therosteritem{\nmb:{4}}
{\it Claim.} For any initial point $(g_0,v_0)\in V_n$ the unique solution 
$\on{Fl}^n_t(g_0,v_0)=(g(t),v(t))$ exists, is $C^2$ in $t\in (-2,2)$,  
and depends $C^1$ on the
initial point in $V_n$.

We use induction on $n\ge 2k+1$. For $n=2k+1$ the claim holds since
$V_{2k+1}=W_{2k+1}$. Let $(g_0,v_0)\in V_{2k+2}$ and let 
$\on{Fl}^{2k+2}_t(g_0,v_0)=(\tilde g(t),\tilde v(t))$ be maximally defined for
$t\in (t_1,t_2)\ni 0$.
Suppose for contradiction that $t_2<2$.
Since $(g_0,v_0)\in V_{2k+2}\subset V_{2k+1}$ the curve   
$\on{Fl}^{2k+2}_t(g_0,v_0)=(\tilde g(t),\tilde v(t))$ solves \nmb|{1} also in 
$\on{DiffHC}^{2k+1}\x HC^{2k+1}$, thus 
$\on{Fl}^{2k+2}_t(g_0,v_0)=(\tilde g(t),\tilde v(t))=(g(t),v(t))
:=\on{Fl}^{2k+1}_t(g_0,v_0)$ for $t\in (t_1,t_2)\cap(-2,2)$.
By \nmb!{7.3.6}, the expression 
\begin{displaymath}
\tilde J(t)=\tilde J(g,v,t)=g_x(t)^2A_k(u(t))\o g(t)
  =g_x(t)^2A_k(v(t)\o g(t)\i)\o g(t) \tag{\nmb:{5}}
\end{displaymath}
is constant in $t\in (-2,2)$.
Actually, since we used $C^\infty$-theory for deriving this, one
should check it again by differentiating.
Since $u=g_t\o g\i$ we get the following (the exact formulas can be
computed with the help of Fa\`a di Bruno's formula \nmb!{6.1}.
\begin{align*}
u_x&=(g_{tx}\o g\i)(g\i)_x= \frac{g_{tx}}{g_x}\o g\i
\\
\p_x^2u&= (\frac{\p_x^2g_{t}}{g_x^2}-g_{tx}\frac{\p_x^2g}{g_x^3})\o g\i
\\
\p_x(g\i)&= \frac1{g_x}\o g\i
\\
\p_x^2(g\i)\o g&= -\frac{\p_x^2g}{g_x^3} 
\\
\p_x^{2k}(g\i)\o g &= -\frac{\p_x^{2k}g}{g_x^{2k+1}} 
+ \text{ lower order terms in }g
\\
(\p_x^{2k}u)\o g &= \frac{\p_x^{2k}g_t}{g_x^{2k}} -
g_{tx}\frac{\p_x^{2k}g}{g_x^{2k+1}}
+ \text{ lower order terms in }g, g_t=v.
\end{align*}
Thus
\begin{displaymath}
(-1)^k g_x^{2k-1}\tilde J(t)=g_x\p_x^{2k}g_t - g_{tx}\p_x^{2k}g 
+ \text{ lower order terms in }g, g_t=v.
\end{displaymath}
Hence for each $t\in (-2,2)$: 
\begin{align*}
g_x\p_x^{2k}g_t - g_{tx}\p_x^{2k}g 
&= (-1)^kg_x^2\left(g_x^{2k-3}\tilde J(t) + P_k(g,v) \right),\text{  where}
\\
P_k(g,v)&=\frac{Q_k(g,\p_xg,\dots,\p_x^{2k-1}g,v,\p_xv,\dots,\p_x^{2k-1}v)}{g_x^2}
\end{align*}
for a polynomial $Q_k$. Since $\tilde J(t)=\tilde J(0)$ we obtain that
\begin{displaymath}
\left(\frac{\p_x^{2k}g(t)}{g_x(t)} \right)_t 
= (-1)^k\left(g_x^{2k-3}(t)\tilde J(0)+ P_k(g(t),v(t)) \right)
\text{  for all }t\in (-2,2).
\end{displaymath}
This implies
\begin{displaymath}
\frac{\p_x^{2k}g(t)}{g_x(t)} = \frac{\p_x^{2k}g(0)}{g_x(0)} +  
(-1)^k\int_0^t\left(g_x^{2k-3}(s)\tilde J(0)+ P_k(g(s),v(s)) \right)\,ds.
\end{displaymath}
For $t\in(t_1,t_2)$ we have 
\begin{align*}
\p_x^{2k}\tilde g(t) &= \frac{\p_x^{2k}g_0}{\p_xg_0}g_x(t) +  
\tag{\nmb:{6}}\\&\quad
+(-1)^kg_x(t)\int_0^t\left(g_x^{2k-3}(s)\tilde J(0)+ P_k(g(s),v(s)) \right)\,ds.
\end{align*}
Since $(g_0,v_0)\in V_{2k+2}$ we have 
$\tilde J(0)=\tilde J(g_0,v_0,0)\in HC^2$ by {\nmb|{5}}. 
Since $k\ge 1$, by {\nmb|{6}} we see that $\p_x^{2k}\tilde g(t)\in HC^2$.
Moreover, since $t_2<2$, the limit $\lim_{t\to t_2-}\p_x^{2k}\tilde g(t)$
exists in $HC^2$, so $\lim_{t\to t_2-}\tilde g(t)$ exists in
$HC^{2k+2}$.  As this limit equals $g(t_2)$, we conclude that $g(t_2)\in
\on{DiffHC}^{2k+2}$. 
Now $\tilde v=\tilde g_t$; so we may differentiate both sides of {\nmb|{6}}
in $t$ and obtain similarly that $\lim_{t\to t_2-}\tilde v(t)$ exists in 
$HC^{2k+2}$ and equals $v(t_2)$. But then we can prolong the flow
line $(\tilde g,\tilde v)$ in $\on{DiffHC}^{2k+2}\x HC^{2k+2}$ 
beyond $t_2$, so $(t_1,t_2)$ was not maximal.  

By the same method we can iterate the induction.
\qed\end{demo}

\section*{\nmb0{8}. The Virasoro-Bott group and the Korteweg-de  
Vries hierarchy } 
 
\subsection*{\nmb.{8.1}. The Virasoro-Bott group  } 
Let $\on{Diff}$ denote any of the groups $\on{DiffHC}^+(S^1)$,
$\on{Diff}(\mathbb R)_0$ (diffeomorphisms with compact support), or
$\on{Diff}_{\mathcal{S}}(\mathbb R)$ of section \nmb!{6}.
For $\ph\in\on{Diff}$ let $\ph':S^1\text{  or }\mathbb R\to \mathbb R^+$ be the mapping 
given by $T_x\ph\cdot\p_x=\ph'(x)\p_x$. 
Then  
\begin{gather*} 
c:\on{Diff}\x\on{Diff}\to \mathbb R
\\
c(\ph,\ps):=\frac12\int_{S^1}\log(\ph\o\ps)'d\log\ps'  
     = \frac12\int_{S^1}\log(\ph'\o\ps)d\log\ps' 
\end{gather*}satisfies $c(\ph,\ph\i)=0$, $c(\on{Id},\ps)=0$,
$c(\ph,\on{Id})=0$, and is a
smooth group cocycle, i.e., 
\begin{displaymath}
c(\ph_2,\ph_3) - c(\ph_1\o\ph_2,\ph_3)+c(\ph_1,\ph_2\o\ph_3) -
c(\ph_1,\ph_2) =0,
\end{displaymath}
called the Bott cocycle.

\begin{demo}{Proof} 
Let us check first:
\begin{align*}
\int \log(\ph\o\ps)'d\log\ps' &=
\int \log((\ph'\o\ps)\ps')d\log\ps' =
\\&
=\int \log(\ph'\o\ps)d\log\ps' +
\int \log(\ps')d\log\ps',
\\
\int \log(\ps')d\log\ps'&=
\tfrac12
\int d\,\log(\ps')^2 = 0.
\\
2c(\on{Id},\ps) &=
\int \log(1)d\log\ps' = 0.
\\
2c(\ph,\on{Id}) &=
\int \log(\ph')d\log(1) = 0.
\\
2c(\ph\i,\ph) &= \int \log((\ph\i\o\ph)')d\log\ph'
= \int \log(1)d\log\ph' = 0.
\\
c(\ph,\ph\i) &= 0.
\end{align*}
For the cocycle condition we add the following terms:
\begin{align*}
2c(\ph_2,\ph_3) &= \int \log(\ph_2'\o\ph_3)d\log\ph_3'
\\
- 2c(\ph_1&\o\ph_2,\ph_3) = -\int \log((\ph_1\o\ph_2)'\o\ph_3)d\log\ph_3'
\\&
=- \int \log((\ph_1'\o\ph_2\o\ph_3)(\ph_2'\o\ph_3))d\log\ph_3'
\\&
= -\int \log(\ph_1'\o\ph_2\o\ph_3)d\log\ph_3'
  -\int \log(\ph_2'\o\ph_3)d\log\ph_3'
\\
2c(\ph_1,&\ph_2\o\ph_3) =
\int \log(\ph_1'\o\ph_2\o\ph_3)d\log(\ph_2\o\ph_3)'
\\&
=\int \log(\ph_1'\o\ph_2\o\ph_3)d\log((\ph_2'\o\ph_3)\ph_3')
\\&
=\int \log(\ph_1'\o\ph_2\o\ph_3)d\log(\ph_2'\o\ph_3)
+\int \log(\ph_1'\o\ph_2\o\ph_3)d\log\ph_3'
\\&
=\int \log(\ph_1'\o\ph_2)d\log\ph_2'
+\int \log(\ph_1'\o\ph_2\o\ph_3)d\log\ph_3'
\\
-2c(\ph_1,\ph_2) &= -\int \log(\ph_1'\o\ph_2)d\log\ph_2' \qquad\qed
\end{align*}\end{demo}

The corresponding central extension group
$S^1\x_c\on{DiffHC}^+(S^1)$, called the periodic Virasoro-Bott group, is
a trivial $S^1$-bundle  
$S^1\x\on{DiffHC}^+(S^1)$ that becomes a regular Lie  
group relative to the operations  
\begin{displaymath}
\binom{\ph}{\al}\binom{\ps}{\be}  
     =\binom{\ph\o\ps}{\al\be\,e^{2\pi ic(\ph,\ps)}},\quad 
\binom{\ph}{\al}\i=\binom{\ph\i}{\al\i}
\end{displaymath}
for $\ph, \ps \in \on{DiffHC}^+(S^1)$ and $\al,\be \in S^1$. 
Likewise we have the central extension group with compact supports 
$\mathbb R\x_c\on{Diff}(\mathbb R)_0$ with group operations
\begin{displaymath}
\binom{\ph}{\al}\binom{\ps}{\be}  
     =\binom{\ph\o\ps}{\al+\be+ c(\ph,\ps)},\quad 
\binom{\ph}{\al}\i=\binom{\ph\i}{-\al}\quad 
\end{displaymath}
for $\ph, \ps \in \on{DiffHC}^+(\mathbb R)$ and $\al,\be \in \mathbb R$. 
Finally there is the central extension of the rapidly decreasing
Virasoro-Bott group $\mathbb R\x_c\on{Diff}_{\mathcal{S}}^+(\mathbb R)$ 
which is given by the same formulas.

\subsection*{\nmb.{8.2}. The Virasoro Lie algebra }
Let us compute the Lie algebra of the two versions of the the Virasoro-Bott
group. Consider $\mathbb R\x_c\on{Diff}$, where again 
$\on{Diff}$ denotes any one of the groups $\on{DiffHC}^+(S^1)$,
$\on{Diff}(\mathbb R)_0$, or
$\on{Diff}_{\mathcal{S}}(\mathbb R)$.
So let $\ph,\ps:\mathbb R\to \on{Diff}$ 
with $\ph(0)=\ps(0)=\on{Id}$ and
$\ph_t(0)=X$, $\ps_t(0)=Y\in X_c(\mathbb R)$, $\X(S^1)$, or $\mathcal{S}(\mathbb R)\p_x$.
For completeness' sake we also consider $\al,\be:\mathbb R\to \mathbb R$ with 
$\al(0)=0$, $\be(0)=0$.
Then we compute:
\begin{align}
&\on{Ad}\binom{\ph(t)}{\al(t)}\,\binom{Y}{\be'(0)} 
= \p_s|_0 \binom{\ph(t)}{\al(t)}\binom{\ps(s)}{\be(s)}\binom{\ph(t)\i}{-\al(t)}
\notag\\&
= \p_s|_0
\binom{\ph(t)\o\ps(s)\o\ph(t)\i}{\al(t)+\be(s)+c(\ph(t),\ps(s))-\al(t)+
c(\ph(t)\o\ps(s),\ph(t)\i)}
\notag\\&
=  \binom{\ph(t)_*Y=\on{Ad}(\ph(t))Y}
{\be_t(0)+\p_s|_0 c(\ph(t),\ps(s))+ \p_s|_0 c(\ph(t)\o\ps(s),\ph(t)\i)}
\tag{\nmb:{1}}\\
&\Big[\binom{X}{\al_t(0)},\binom{Y}{\be_t(0)} \Big] =
\notag\\&
= \p_t|0 \binom{(\on{Fl}^X_t)_*Y=\on{Ad}(\ph(t))Y}
{\be_t(0)+\p_s|_0 c(\ph(t),\ps(s))+ \p_s|_0 c(\ph(t)\o\ps(s),\ph(t)\i)}
\notag\\&
=\binom{-[X,Y]}
{\p_t|_0\p_s|_0 c(\ph(t),\ps(s))+ \p_t|_0\p_s|_0 c(\ph(t)\o\ps(s),\ph(t)\i)}
\tag{\nmb:{2}}\end{align}
Now we differentiate the Bott cocycle, where sometimes $f'=\p_xf$:
\begin{align*}
2\p_s|_0 c(\ph(t),\ps(s)) &= \p_s|_0 \int\log(\ph(t)'\o\ps(s))\,d\log(\ps(s)')
\\&
=\int\frac{(\ph(t)''\o\ps(0))Y}{\ph(t)'\o\ps(0)}\,
  d\log(\undersetbrace =1\to{\ps(0)'})
  +\int \log(\ph(t)')\,dY'
\\&
=\int\log(\ph(t)')Y''\,dx
\\
2\p_t|_0\p_s|_0 c(\ph(t),\ps(s)) 
&=\p_t|_0\int\log(\ph(t)')Y''\,dx 
=\int\frac{X'Y''}{\ph(0)'}dx = \int X'Y''dx.
\end{align*}
For the second term we first check: 
\begin{align*}
(\ph\i)_x &= \frac1{\ph_x\o\ph\i},\qquad
(\ph\i)_{xx}=-\frac{\ph_{xx}\o\ph\i}{(\ph_x\o\ph\i)^3},
\\
\ph\i(x)&=y, \qquad \frac1{\ph_x\o\ph\i}dx = dy
\\
d\log((\ph\i)_x) &= - \frac{\ph''\o\ph\i}{(\ph'\o\ph\i)^2}dx 
=  - \frac{\ph''}{\ph'}dy
\end{align*}
and continue to compute
\begin{align*}
&2\p_s|_0 c(\ph(t)\o\ps(s),\ph(t)\i) 
= \p_s|_0 \int\log((\ph(t)\o\ps(s))_x\o\ph(t)\i)\,d\log(\ph(t)\i_x)
\\&
=\! \int\!\frac{(\ph(t)''\!\o\ph(t)\i)(Y\!\o\ph(t)\i)
  +(\ph(t)'\!\o\ph(t)\i)(Y'\!\o\ph(t)\i)}{(\ph(t)'\o\ph(t)\i)(\ps(0)'\o
  \ph(t)\i)}\,
  d\log(\ph(t)\i_x)
\\&
= -\int\frac{(\ph(t)'')^2\,Y    +\ph(t)'\ph(t)''\,Y'}{(\ph(t)')^2}\,dy
\\
&2\p_t|_0\p_s|_0 c(\ph(t)\o\ps(s),\ph(t)\i) 
= -\p_t|_0\int\frac{(\ph(t)'')^2\,Y     +\ph(t)'\ph(t)''\,Y'}{(\ph(t)')^2}\,dy
\\&
= -\int\frac{0  + 0 +\ph(0)'X''\,Y' - 0}{(\ph(0)'=1)^4}\,dy
\\&
= -\int X''Y' \,dy =\int X'Y''\,dx.
\end{align*}
Finally we get from {\nmb|{2}}: 
\begin{equation}
\left[\binom{X}{a},\binom{Y}{b} \right]
=\binom{-[X,Y]}{\om(X,Y)} = \binom{X'Y-XY'}{\om(X,Y)}
\tag{\nmb:{3}}\end{equation}
where
\begin{displaymath}
\om(X,Y)=\om(X)Y=\int X'dY'=\int X'Y''dx =  
\tfrac12\int \det\begin{pmatrix} X'& Y'\\ X''&Y''\end{pmatrix}\,dx, 
\end{displaymath}
is the \idx{\it Gelfand-Fuchs Lie algebra cocycle}  
$\om:\g\x \g\to \mathbb R$, which is a bounded skew-symmetric bilinear mapping
satisfying the cocycle condition
\begin{displaymath}
\om([X,Y],Z)+\om([Y,Z],X)+\om([Z,X],Y)=0.
\end{displaymath}
It is a generator of the 1-dimensional bounded Chevalley cohomology  
$H^2(\g,\mathbb R)$ for any of the Lie algebras 
$\g=\X(S^1)$, $\X_c(\mathbb R)$, or $\mathcal{S}(\mathbb R)\p_x$.
The Lie algebra of the Virasoro-Bott Lie group is thus the central extension 
$\mathbb R\x_\om \g$ of $\g$ induced by this cocycle.
We have $H^2(\X_c(M),\mathbb R)=0$ for each
finite dimensional manifold of dimension $\ge 2$ (see \cite{15}),  
which blocks the way to  
find a higher dimensional analog of the Korteweg -- de Vries  
equation in a way similar to that sketched below. 

For further use we also note the expression for the adjoint
action on the Virasoro-Bott groups which we computed along the way. For the
integral in the central term in {\nmb|{1}} we have: 
\begin{multline*}
\frac12\int\Bigl(\log(\ph')Y''-\frac{(\ph'')^2Y+\ph'\ph''Y'}{(\ph')^2}\Bigr)\,dx
=\frac{1}2\int\Bigl(-2\frac{\ph''}{\ph'}Y'-\Bigl(\frac{\ph''}{\ph'}\Bigr)^2Y\Bigr)\,dx
=\\=
\int\Bigl(\Bigl(\frac{\ph''}{\ph'}\Bigr)'
  -\frac{1}2\Bigl(\frac{\ph''}{\ph'}\Bigr)^2\Bigr)Y\,dx
=\int S(\ph)Y\,dx,
\end{multline*}
where a new character appears on stage, the \idx{\it Schwartzian
derivative}:
\begin{align}
S(\ph) &
= \Bigl(\frac{\ph''}{\ph'}\Bigr)' -\frac{1}2\Bigl(\frac{\ph''}{\ph'}\Bigr)^2
= \frac{\ph'''}{\ph'} -\frac{3}2\Bigl(\frac{\ph''}{\ph'}\Bigr)^2
= \log(\ph')'' -\frac12(\log(\ph')')^2
\tag{\nmb:{4}}
\end{align}
which measures the deviation of $\ph$ from being a Moebius transformation:
\begin{displaymath}
S(\ph)=0 \iff \ph(x)=\frac{ax+b}{cx+d}\text{  for }
\begin{pmatrix} a & b \\ c & d\end{pmatrix} \in SL(2,\mathbb R).
\end{displaymath}
Indeed, $S(\ph)=0$ if and only if $g=\log(\ph')'=\frac{\ph''}{\ph'}$
satisfies the differential equation $g'= g^2/2$, so that
$\frac{2\,dg}{g^2}=dx$ or $\frac{-2}{g}=x+\frac dc$
which means
$\log(\ph')'(x)=g(x)=\frac{-2}{x+d/c}$ or again
$\log(\ph'(x))=\int\frac{-2dx}{x+d/c}=-2\log(x+d/c)-2\log(c)=\log(\frac 1{(cx+d)^{2}})$.
Therefore, $\ph'(x)= \frac 1{(cx+d)^{2}}=\p_x\frac{ax+b}{cx+d}$.

For completeness' sake, let us note here the Schwartzian derivative of a
composition and an inverse (which follow since the adjoint action
{\nmb|{5}} below is an action):
\begin{displaymath}
S(\ph\o \ps) = (S(\ph)\o\ps)(\ps')^2 + S(\ps),\quad
S(\ph\i)=-\frac{S(\ph)}{(\ph')^2}\o\ph\i
\end{displaymath}
So finally,  the adjoint action is given by: 
\begin{equation}
\on{Ad}\binom{\ph}{\al}\,\binom{Y}{b} 
=  \binom{\on{Ad}(\ph)Y=\ph_*Y=T\ph\o Y\o \ph\i}
{b+\int S(\ph)Y\,dx}
\tag{\nmb:{5}}\end{equation}

\subsection*{\nmb.{8.3}. $H^0$-Geodesics on the Virasoro-Bott groups
} 
We shall use the $L^2$-inner product on $\mathbb R\x_\om \g$, where 
$\g=\X(S^1), \X_c(\mathbb R), \mathcal{S}(\mathbb R)\p_x$: 
\begin{equation}
\left\langle \binom{X}{a},\binom{Y}{b}\right\rangle_0 :=  
\int XY\,dx + ab. 
\tag{\nmb:{1}}\end{equation}
Integrating by parts we get 
\begin{align*} 
\left\langle \ad\binom{X}{a}\binom{Y}{b},\binom{Z}{c}\right\rangle_0  
&= \left\langle \binom{X'Y-XY'}{\om(X,Y)}, 
     \binom{Z}{c}\right\rangle_0
\\&
= \int (X'YZ-XY'Z+cX'Y'')\,dx 
\\&
= \int (2X'Z+XZ'+cX''')Y\,dx
\\&
= \left\langle \binom{Y}{b}, 
     \ad\binom{X}{a}^\top\binom{Z}{c}\right\rangle_0, 
          \quad {\text {where}}
\\
\ad\binom{X}{a}^\top\binom{Z}{c}&=\binom{2X'Z+XZ'+cX'''}{0}. 
\end{align*}Using matrix notation we get therefore (where $\p :=
\p_x$) 
\begin{align*} 
\ad\binom{X}{a} &= \begin{pmatrix} X' - X\p & 0 \\  
                                   \om(X) & 0 \end{pmatrix} 
\\
\ad\binom{X}{a}^\top &= \begin{pmatrix} 2X'+X\p & X'''\\ 
                                      0 & 0 \end{pmatrix}  
\\
\al\binom{X}{a} &=\ad\binom{\quad}{\quad}^\top\binom{X}{a}  
     = \begin{pmatrix} X'+2X\p +a\p^3 & 0 \\ 
                                      0 & 0 \end{pmatrix}  
\\
\ad{\binom{X}{a}}^\top + \ad\binom{X}{a} 
     &= \begin{pmatrix} 3X' & X''' \\ 
                \om(X) & 0 \end{pmatrix}  
\\
\ad\binom{X}{a}^\top - \ad\binom{X}{a} 
     &= \begin{pmatrix} X'+2X\p & X''' \\ 
                -\om(X) & 0 \end{pmatrix}.  
\\
\end{align*}Formula \nmb!{3.2.2} gives the $H^0$ geodesic equation on the  
Virasoro-Bott group:   
\begin{align}
\binom{u_t}{a_t}&=-\ad\binom{u}{a}^\top\binom{u}{a} 
     =\binom{-3u_xu-au_{xxx}}{0}\quad\text{  where } 
\tag{\nmb:{2}}\\
\binom{u(t)}{a(t)}&=\p_s\binom{\ph(s)}{\al(s)}.\binom{\ph(t)\i}{-\al(t)}\Bigr|_{s=t}
\notag\\&
=\p_s\binom{\ph(s)\o\ph(t)\i}{\al(s)-\al(t)+c(\ph(s),\ph(t)\i)}\Bigr|_{s=t}
\notag\\&
=\binom{\ph_t\o\ph\i}{\al_t-\int \frac{\ph_{tx}\ph_{xx}}{2\ph_x^2}dx}
\notag\end{align}
since we have
\begin{align*}
2\p_s &c(\ph(s),\ph(t)\i)|_{s=t} = \p_s\int\log(\ph(s)'\o
  \ph(t)\i)\,d\log((\ph(t)\i)')|_{s=t}
\\&
= \int \frac{\ph_{t}(t)'\o \ph(t)\i}{\ph(t)'\o
\ph(t)\i}\left(-\frac{\ph(t)''\o\ph(t)\i}{(\ph(t)'\o\ph(t)\i)^2} \right)\,dx
\quad\text{  by }\nmb!{8.2}
\\&
= -\int \frac{\ph_{t}'\ph''}{(\ph')^2}dy
= -\int \frac{\ph_{tx}\ph_{xx}}{\ph_x^2}dx.
\end{align*}
Thus $a$ is a constant in time and the geodesic equation  
is hence the \idx{\it Korteweg-de Vries equation} 
\begin{equation}
u_t+3u_xu+au_{xxx}=0. 
\tag{\nmb:{3}}\end{equation}
with its natural companions
\begin{displaymath}
\ph_t=u\o\ph,\qquad \al_t = a+\int \frac{\ph_{tx}\ph_{xx}}{2\ph_x^2}dx. 
\end{displaymath}
It is the periodic equation, if we work on $S^1$.

The derivation above is direct and
does not use the Euler-Poincar\'e equations; for a derivation of the
Korteweg-de Vries   equation from this point of view 
see \cite{25}, section 13.8.

Let us compute the invariant momentum mapping from \nmb!{4.3.2}. First we
need the transpose of the adjoint action \nmb!{8.2.5}:
\begin{align*}
&\left\langle \on{Ad}\binom{\ph}{\al}^\top\,\binom{Y}{b},\binom{Z}{c} \right\rangle_0
=\left\langle \binom{Y}{b},\on{Ad}\binom{\ph}{\al}\,\binom{Z}{c} \right\rangle_0
\\&\qquad
=\left\langle \binom{Y}{b},\binom{\ph_*Z}{c+\int S(\ph)Z\,dx}\right\rangle_0
\\&\qquad
=\int Y((\ph'\o\ph\i)(Z\o\ph\i)\,dx +bc+\int bS(\ph)Z\,dx
\\&\qquad
=\int ((Y\o\ph)(\ph')^2 +bS(\ph))Z\,dx + bc
\\
&\on{Ad}\binom{\ph}{\al}^\top\,\binom{Y}{b} =
\binom{(Y\o\ph)(\ph')^2 +bS(\ph)}{b}.
\end{align*}
Thus the invariant momentum mapping \nmb!{4.3.2}
turns out as 
\begin{equation}
\bar J\left(\binom{\ph}{\al},\binom{Y}{b} \right) =
\on{Ad}\binom{\ph}{\al}^\top\binom{Y}{b} 
\binom{(Y\o\ph)(\ph')^2 +bS(\ph)}{b}.
\tag{\nmb:{4}}\end{equation}
Along a geodesic $t\mapsto g(t,\quad)=\binom{\ph(t,\quad)}{\al(t)}$, 
according to {\nmb|{3}} and
\nmb!{4.3}, the momentum 
\begin{equation}
\bar J\left(\binom{\ph}{\al},\binom{u=\ph_t\o \ph\i}{a}\right) 
= \binom{(u\o \ph)\ph_x^2 +aS(\ph)}{a} = \binom{\ph_t\ph_x^2+aS(\ph)}{a}
\tag{\nmb:{5}}\end{equation}
is constant in $t$.

\subsection*{\nmb.{8.4}. The curvature } 
The computation of the curvature at the identity element has been done
independently by \cite{31} and Misiolek \cite{32}. Here
we proceed with a completely general computation that takes advantage
of the formalism introduced so far. Inserting the matrices of
differential-   and integral   operators
$\ad\tbinom{X}{a}^\top$,
$\al\tbinom{X}{a}$, and  
$\ad\tbinom{X}{a}$ etc.\ 
given above into formula \nmb!{3.4.2} and recalling that the matrix  
is applied to vectors of the form $\binom{Z}{c}$, where $c$ is a  
constant, we see that  
$4\mathcal{R}\left(\tbinom{X_1}{a_1},\tbinom{X_2}{a_2}\right)$ is the following  
$2\x 2$-matrix whose entries are differential- and integral  
operators:  
\begin{displaymath}
\begin{pmatrix} 
{\begin{aligned} 
&4(X_1X_2''-X_1''X_2) +2(a_1X_2^{(4)}-a_2X_1^{(4)})
\\
+(&8(X_1X_2'-X_1'X_2)+10(a_1X_2'''-a_2X_1'''))\p
\\
+&18(a_1X_2''-a_2X_1'')\p^2
\\
+&(12(a_1X_2'-a_2X_1')+2\om(X_1,X_2))\p^3
\\
-&X_1'''\om(X_2)+X_2'''\om(X_1) 
\end{aligned}} 
&\qquad& 
{\begin{aligned}  
&2(X_1'''X_2' - X_1'X_2''')
\\
+&2(X_1X_2^{(4)}-X_1^{(4)}X_2)
\\
+&(a_1X_2^{(6)}-a_2X_1^{(6)})\end{aligned}} 
\\
\vphantom{\int_A^A} & & 
\\
{\begin{aligned} 
&\om(X_2)(4X_1'+2X_1\p+a_1\p^3)
\\
-&\om(X_1)(4X_2'+2X_2\p+a_2\p^3) 
\end{aligned}} 
&& 0 \end{pmatrix} 
\end{displaymath}
Therefore, $4\mathcal{R}\left(\tbinom{X_1}{a_1},\tbinom{X_2}{a_2}\right)
\tbinom{X_3}{a_3}$ has the following expression
\begin{displaymath}
\begin{pmatrix}
{\begin{aligned}
&4(X_1X_2'' - X_1''X_2)X_3 + 2(a_1X_2^{(4)} - a_2X_1^{(4)})X_3
\\&
+\left(8(X_1X_2'-X_1'X_2) + 10(a_1X_2''' - a_2X_1''')\right)X_3'
\\&
+18(a_1X_2'' - a_2X_1'')X_3'' + 12(a_1X_2' - a_2X_1')X_3'''
\\&
+2X_3'''\int X_1'X_2''dx - X_1'''\int X_2'X_3''dx
+X_2'''\int X_1'X_3''dx
\\&
+2a_3(X_1'''X_2'\! -\! X_1'X_2''')
+2a_3(X_1X_2^{(4)}\! -\! X_1^{(4)}X_2)
+a_3(a_1X_2^{(6)}\! -\! a_2X_1^{(6)})
\end{aligned}}
\\
\vphantom{\int_A^A}
\\
{\begin{aligned}
&\int  X_3'''(a_1X_2'''-a_2X_1''')dx
\\&
+\int 2X_3'(X_1X_2'''-X_1'''X_2 - 2X_1'X_2'' + 2X_1''X_2')dx
\end{aligned}}
\end{pmatrix} 
\end{displaymath}
which coincides with formula (2.3) in Misiolek \cite{32}.
This in turn leads to the following expression for the sectional
curvature  
\begin{align*} 
&\left\langle4\mathcal{R}\left(\tbinom{X_1}{a_1},\tbinom{X_2}{a_2}\right) 
     \tbinom{X_1}{a_1},\tbinom{X_2}{a_2}\right\rangle_0 = 
\\
=\int \Bigl( 
&4(X_1X_2''-X_1''X_2)X_1X_2 +8(X_1X_2'-X_1'X_2)X_1'X_2
\\&
+2(a_1X_2^{(4)}-a_2X_1^{(4)})X_1X_2+10(a_1X_2'''-a_2X_1''')X_1'X_2
\\&
+18(a_1X_2''-a_2X_1'')X_1''X_2
\\&
+12(a_1X_2'-a_2X_1')X_1'''X_2+2\om(X_1,X_2)X_1'''X_2
\\&
-X_1'''\om(X_2,X_1)X_2+X_2'''\om(X_1,X_1)X_2
\\&
+2(X_1'''X_2' - X_1'X_2''')a_1X_2
\\&
+2(X_1X_2^{(4)}-X_1^{(4)}X_2)a_1X_2
\\&
+(a_1X_2^{(6)}-a_2X_1^{(6)})a_1X_2 
\\&
+(4X_1'X_1X_2'''+2X_1X_1'X_2'''+a_1X_1'''X_2'''
\\&
\qquad-4X_2'X_1X_1'''-2X_2X_1'X_1'''-a_2X_1'''X_1''')a_2    \Bigr)\;dx
\\
%\allowdisplaybreak 
%=\int \Bigl( 
%&-4[X_1,X_2]'X_1X_2 -8[X_1,X_2]X_1'X_2\\ 
%&+2(a_1X_2^{(4)}-a_2X_1^{(4)})X_1X_2+10(a_1X_2'''-a_2X_1''')X_1'X_2\\ 
%&+18(a_1X_2''-a_2X_1'')X_1''X_2\\ 
%&+12(a_1X_2'-a_2X_1')X_1'''X_2+2\om(X_1,X_2)X_1'''X_2\\ 
%&-X_1'''\om(X_2,X_1)X_2+X_2'''\om(X_1,X_1)X_2\\ 
%&+2(X_1'''X_2' - X_1'X_2''')(a_1X_2 - a_2X_1)\\ 
%&+2(X_1X_2^{(4)}-X_1^{(4)}X_2)(a_1X_2 -a_2X_1)\\  
%&+(a_1X_2^{(6)}-a_2X_1^{(6)})a_1X_2 \\ 
%&+(a_1X_1'''X_2'''-a_2X_1'''X_1''')a_2    \Bigr)\;dx\\ 
%\allowdisplaybreak 
%=\int \Bigl( 
%&4[X_1,X_2](X_1X_2' - X_1'X_2)\\ 
%&+4(a_1X_2-a_2X_1)(X_1X_2^{(4)}-X_1^{(4)}X_2) 
%     -4(a_1X_2-a_2X_1)(X_1'X_2'''-X_1'''X_2')\\ 
%&-(X_2''')^2a_1^2 +2X_1'''X_2'''a_1a_2 -(X_1''')^2a_2^2 \Bigr)\;dx\\ 
%+3\om(&X_1,X_2)^2 \\ 
%\allowdisplaybreak 
=\int \Bigl( 
&-4[X_1,X_2]^2 
+4(a_1X_2-a_2X_1)(X_1X_2^{(4)}-X_1'X_2'''+X_1'''X_2'-X_1^{(4)}X_2)
\\&
-(X_2''')^2a_1^2 +2X_1'''X_2'''a_1a_2 -(X_1''')^2a_2^2 \Bigr)\;dx
\\
+3\om(&X_1,X_2)^2 .
\\
\end{align*}
This formula shows that the sign of the sectional curvature is not
constant. Indeed, choosing $h_1(x) = \sin x, h_2(x) = \cos x$ we get 
$-\pi(8 + a_1^2 + a_2^2 - 3\pi)$ which can be positive and negative
by choosing the constants $a_1, a_2$ judiciously.

\subsection*{\nmb.{8.5}. Jacobi fields } 
A Jacobi field $y=\binom{y}{b}$ along a geodesic with velocity field  
$\binom{u}{a}$ is a solution of the partial differential  
equation \nmb!{3.5.1} which in our case looks as follows.  
%Here $u'$ or $y^{(4)}$ denote derivatives with respect to $x$. 
\begin{align*} 
\binom{y_{tt}}{b_{tt}} 
&= \left[\ad\binom{y}{b}^\top+\ad\binom{y}{b}, 
                     \ad\binom{u}{a}^\top\right]\binom{u}{a} 
\\&
\quad- \ad\binom{u}{a}^\top \binom{y_t}{b_t}  
     -\al\binom{u}{a}\binom{y_t}{b_t} + \ad\binom{u}{a}\binom{y_t}{b_t}
\\&
= \left[\begin{pmatrix} 3y_x & y_{xxx} \\
                  \om(y) & 0 \end{pmatrix}, 
           \begin{pmatrix} 2u_x+u\p_x & u_{xxx}\\ 
                     0 & 0 \end{pmatrix}      \right]\binom{u}{a} 
\\&
\quad+ \begin{pmatrix} -2u_x-4u\p_x-a\p_x^3 & -u_{xxx}\\ 
                   \om(u) & 0 \end{pmatrix} 
\binom{y_t}{b_t}, \\ 
\end{align*}which leads to  
\begin{align} 
y_{tt}&=-u(4y_{tx}+3uy_{xx}+ay_{xxxx}) -u_x(2y_t+2ay_{xxx})
\tag{\nmb:{1}}\\&
\qquad-u_{xxx}(b_t+\om(y,u)-3ay_x)-ay_{txxx}, 
\notag\\
b_{tt}&=\om(u,y_t)+\om(y,3u_xu)+\om(y,au_{xxx}).
\tag{\nmb:{2}}\end{align}
Equation {\nmb|{2}} is equivalent to: 
\begin{equation}
b_{tt} = \int (-y_{txxx}u + y_{xxx}(3u_xu+au_{xxx}))dx. 
\tag{2$'$}\end{equation}
Next, let us show that the integral term in equation
{\nmb|{1}} is constant:
\begin{equation}
b_t+\om(y,u) = b_t +\int  y_{xxx}u\,dx =: B_1. 
\tag{\nmb:{3}}\end{equation}
Indeed its $t$-derivative along the geodesic for $u$ (that is, $u$
satisfies the Korteweg-de Vries equation)
coincides with \thetag{2$'$}:
\begin{displaymath}
b_{tt} + \int  (y_{txxx}u + y_{xxx}u_t)\, dx  
= b_{tt} + \int  (y_{txxx}u + y_{xxx}(-3u_xu-au_{xxx}))\, dx = 0. 
\end{displaymath}
Thus $b(t)$ can be explicitly solved from {\nmb|{3}} as 
\begin{equation}
b(t) = B_0 + B_1t - \int_a^t\!\!\int y_{xxx}u\,dx\,dt.
\tag{\nmb:{4}}\end{equation}
The first component of the Jacobi equation on the Virasoro-Bott group  
is a genuine partial differential equation. Thus the Jacobi equations
are given by the following  system: 
\begin{align} 
y_{tt} &= -u(4y_{tx}+3uy_{xx}+ay_{xxxx}) -u_x(2y_t+2ay_{xxx})
\notag\\&
\qquad-u_{xxx}(B_1-3ay_x)-ay_{txxx}, 
\tag{\nmb:{5}}\\
u_t&= -3u_xu - au_{xxx},
\notag\\
a&= \text{ constant, } 
\notag\end{align}
where $u(t,x),y(t,x)$ are either smooth functions in  
$(t,x)\in I\x S^1$ or in  
$(t,x)\in I\x \mathbb R$, where $I$ is an interval or $\mathbb R$, and  
where in the latter case $u$, $y$, $y_t$  
have compact support with respect to $x$.

Choosing $u = c \in \mathbb R$, a constant, these equations coincide with
(3.1) in Misiolek \cite{32} where it is shown by direct inspection that
there are solutions of this equation which vanish at non-zero values
of $t$, thereby concluding that there are conjugate points along 
geodesics emanating from the identity element of the Virasoro-Bott
group on $S^1$.
 
\subsection*{\nmb.{8.6}. The weak symplectic structure on the space of  
Jacobi fields on the Virasoro Lie algebra } 
Since the Korteweg - de Vries equation has local solutions  
depending smoothly on the initial conditions (and global solutions if  
$a\ne 0$), we expect that 
the space of all Jacobi fields exists and is isomorphic to the space 
of all initial data  
$(\mathbb R\x_\om \X(S^1))\x(\mathbb R\x_\om \X(S^1))$. The weak symplectic  
structure is given in section \nmb!{3.7}: 
\begin{align} 
\om\left(\binom{y}{b},\binom{z}{c}\right)  
&= \left\langle \binom{y}{b},\binom{z_t}{c_t}\right\rangle_0  
    \! - \left\langle \binom{y_t}{b_t},\binom{z}{c}\right\rangle_0 
    \! + \left\langle \left[\binom{u}{a},\binom{y}{b}\right], 
          \binom{z}{c}\right\rangle_0 
\notag\\&
\qquad- \left\langle \binom{y}{b}, 
          \left[\binom{u}{a},\binom{z}{b}\right] \right\rangle_0  
     - \left\langle \left[\binom{y}{b},\binom{z}{c}\right], 
          \binom{u}{a}\right\rangle_0 
\notag\\&
= \int (yz_t-y_tz+2u(yz_x-y_xz))\,dx 
\notag\\&
\qquad+b(c_t+\om(z,u)) -c(b_t+\om(y,u)) -a\om(y,z) 
\notag\\&
= \int (yz_t-y_tz+2u(yz_x-y_xz))\,dx 
\tag{\nmb:{1}}\\&
\qquad+bC_1 -cB_1 -a\int y'z''\,dx,
\end{align}
where the constant $C_1$ relates to $c$ as $B_1$ does to $b$, see 
\nmb!{8.5.3} and \nmb!{8.5.4}. 

\subsection*{\nmb.{8.7}. The geodesics of the $H^k$-metric on the Virasoro group
}
We shall use the $H^k$-inner product on $\mathbb R\x_\om \g$, where $\g$ is
any of the Lie algebras $\X(S^1)$ or 
$\X_{\mathcal S}(\mathbb R)=\mathcal{S}(\mathbb R)\p_x$. The Lie algebra
$\X_c(\mathbb R)$ does not work here any more since
$A_k=\sum_{j=0}^k(-1)^j\p_x^{2j}$ 
is no longer a linear isomorphism here. 
\begin{align}
\left\langle \binom{X}{a},\binom{Y}{b}\right\rangle_k :&=  
\int (XY+X'Y'+\dots+X^{(k)}Y^{(k)})\,dx + ab
\tag{\nmb:{1}}\\&
=\int A_k(X)Y\,dx + ab
=\int XA_k(Y)\,dx + ab,
\notag\\&\text{  where }A_k=\sum_{i=0}^k(-1)^i\p_x^{2i} \text{  as in
\nmb!{7.3.1}.}
\notag\end{align}
Integrating by parts we get 
\begin{align} 
\left\langle \ad\binom{X}{a}\binom{Y}{b},\binom{Z}{c}\right\rangle_k  
&= \left\langle \binom{X'Y-XY'}{\om(X,Y)}, 
     \binom{Z}{c}\right\rangle_k
\notag\\&
= \int (X'YA_k(Z)-XY'A_k(Z)+cX'Y'')\,dx 
\notag\\&
= \int (2X'YA_k(Z)+XYA_k(Z')+cX''')\,dx
\notag\\&
= \int YA_kA_k\i(2X'A_k(Z)+XA_k(Z')+cX''')\,dx
\notag\\&
= \left\langle \binom{Y}{b}, 
     \ad\binom{X}{a}^\top\binom{Z}{c}\right\rangle_0, 
          \quad {\text {where}}
\notag\\
\ad\binom{X}{a}^\top\binom{Z}{c}&=\binom{A_k\i(2X'A_k(Z)+XA_k(Z')+cX''')}{0}. 
\tag{\nmb:{2}}\end{align}
Using matrix notation we get therefore (where $\p :=
\p_x$) 
\begin{align*} 
\ad\binom{X}{a} &= \begin{pmatrix} X' - X\p & 0 \\  
                                   \om(X) & 0 \end{pmatrix} 
\\
\ad\binom{X}{a}^\top &= \begin{pmatrix} A_k\i.(2X'.A_k+XA_k.\p_x) & A_k\i(X''')\\ 
                                      0 & 0 \end{pmatrix}  
\\
\al\binom{X}{a} &=\ad\binom{\quad}{\quad}^\top\binom{X}{a}  
     = \begin{pmatrix} A_k\i.(A_k(X')+2A_k(X)\p_x +a\p^3) & 0 \\ 
                                      0 & 0 \end{pmatrix}.  
\end{align*}
Formula \nmb!{3.2.2} gives the geodesic equation on the  
Virasoro-Bott group:   
\begin{align}
\binom{u_t}{a_t}&=-\ad\binom{u}{a}^\top\binom{u}{a} 
     =\binom{-A_k\i(2u_xA_k(u)+uA_k(u_x)+au_{xxx})}{0},
\tag{\nmb:{3}}\\&
\text{ where }\binom{u(t)}{a(t)}
=\binom{\ph_t\o\ph\i}{\al_t-\int \frac{\ph_{tx}\ph_{xx}}{2\ph_x^2}dx}
\notag\end{align}
as in \nmb!{8.3.2}
Thus $a$ is a constant in time and the geodesic equation  
contains the equation from the Korteweg-de Vries hierarchy:  
\begin{displaymath}
A_k(u_t)=-2u_xA_k(u)-uA_k(u_x)-au_{xxx}
\tag{\nmb:{4}}
\end{displaymath}
For $k=0$ this gives the Korteweg--de~Vries equation.

For $k=1$ we get the equation
\begin{displaymath}
u_t-u_{txx}=-3uu_x+2u_xu_{xx}+uu_{xxx}-au_{xxx},
\end{displaymath}
the {\it Camassa-Holm equation}, \cite{11}, \cite{26}.  
See  \nmb!{7.3.4} for the dispersionfree version. 

Let us compute the invariant momentum mapping from \nmb!{4.3.2}. First we
need the transpose of the adjoint action \nmb!{8.2.5}:
\begin{align*}
&\left\langle \on{Ad}\binom{\ph}{\al}^\top\,\binom{Y}{b},\binom{Z}{c}
\right\rangle_k
=\left\langle \binom{Y}{b},\on{Ad}\binom{\ph}{\al}\,\binom{Z}{c}
\right\rangle_k
\\&\qquad
=\left\langle \binom{Y}{b},\binom{\ph_*Z}{c+\int
S(\ph)Z\,dx}\right\rangle_k
\\&\qquad
=\int A_k(Y)(\ph_*Z)\,dx +bc+\int bS(\ph)Z\,dx
\\&\qquad
=\int A_k(Y)((\ph'Z)\o\ph\i)\,dx +bc+\int bS(\ph)Z\,dx
\\&\qquad
=\int (A_k(Y)\o\ph)(\ph')^2 Z\,dy +bc+\int bS(\ph)Z\,dx
\\&\qquad
=\int ((A_k(Y)\o\ph)(\ph')^2 +bS(\ph))Z\,dx + bc
\\&\qquad
=\int A_k A_k\i((A_k(Y)\o\ph)(\ph')^2 +bS(\ph))Z\,dx + bc
\\&\qquad
= \left\langle \binom{A_k\i((A_k(Y)\o\ph)(\ph')^2 +bS(\ph))}{b},
\binom{Z}{c}\right\rangle.
\\
&\on{Ad}\binom{\ph}{\al}^\top\,\binom{Y}{b} =
\binom{A_k\i((A_k(Y)\o\ph)(\ph')^2 +bS(\ph))}{b}
\end{align*}
Thus the invariant momentum mapping \nmb!{4.3.2}
turns out as 
\begin{equation}
\bar J\left(\binom{\ph}{\al},\binom{Y}{b} \right) =
\on{Ad}\binom{\ph}{\al}^\top\binom{Y}{b} 
=\binom{A_k\i((A_k(Y)\o\ph)(\ph')^2 +bS(\ph))}{b}.
\tag{\nmb:{5}}\end{equation}
Along a geodesic $t\mapsto g(t,\quad)=\binom{\ph(t,\quad)}{\al(t)}$, 
according to {\nmb|{4}} and
\nmb!{4.3}, the momentum  
\begin{align}
\bar J\left(\binom{\ph}{\al},\binom{u=\ph_t\o \ph\i}{a}\right) 
&= \binom{A_k\i\bigl(A_k(u)\o \ph)\ph_x^2 +aS(\ph)\bigr)}{a}
\notag\\&
= \binom{A_k\i\bigl((A_k(\ph_t\o\ph\i)\o \ph)\ph_x^2 +aS(\ph)\bigr)}{a}
\tag{\nmb:{6}}\end{align}
is constant in $t$, and thus also 
\begin{displaymath}
\tilde J(a,\ph) := \bigl(A_k(\ph_t\o\ph\i)\o \ph\bigr)\ph_x^2 +aS(\ph)
\tag{\nmb:{7}}
\end{displaymath}
is constant in $t$. 

\begin{proclaim}{\nmb.{8.8}.  Theorem} 
{\rm \cite{12}} Let $k\ge 2$. 
There exists a $HC^{2k+1}$-open neighborhood $V$ of $(\on{Id},0)$ in the
space 
$(S^1\x_c\on{Diff}(S^1))\x(\mathbb R\x_\om \X(S^1))$ such that for each 
$(g_0,\al,u_0,a)\in V$ there exists a unique $C^3$ geodesic
$g\in C^3((-2,2),S^1\x_c\on{Diff}(S^1))$ for the right invariant $H^k$
Riemann metric, starting at $g(0)=g_0$ in the direction $g_t(0)=u_0\o g_0\in
T_{g_0}\on{Diff}(S^1)$. Moreover, the solution depends $C^1$ on the initial
data $(g_0,u_0)\in V$.

The same result holds if we replace $S^1\x_c\on{Diff}(S^1)$ by $\mathbb R\x_c
\on{Diff}_{\mathcal S}(\mathbb R)$ and $\X(S^1)$ by $\mathcal S(\mathbb
R)\p_x=\X_{\mathcal S}(\mathbb R)$. 
\end{proclaim}

In the following proof $\on{Diff}$, $\X$, $\on{DiffHC}^n$, $HC^n$ will mean either 
$\on{Diff}(S^1)$, $\X(S^1)$, $\on{DiffHC}^n(S^1)$, $HC^n(S^1)$, or 
$\on{Diff}_{\mathcal S}(\mathbb R)$, $\X_{\mathcal S}(\mathbb R)$, 
$\on{DiffHC}^n(\mathbb R)$, $HC^n(\mathbb R)$, respectively.

\begin{demo}{Proof}
For $u\in HC^n$, $n\ge 2k+1$, we have as in the proof of \nmb!{7.4}
\begin{displaymath}
%\begin{matrix}
A_k(uu_x)
%&=\sum_{i=0}^k(-1)^i\p_x^{2i}(uu_x)
%=\sum_{i=0}^k(-1)^i\sum_{j=0}^{2i}\tbinom{2i}{j}(\p_x^{j}u)(\p_x^{2i-j+1}u)
%\\&
=uA_k(u_x)+\sum_{i=0}^k(-1)^i\sum_{j=1}^{2i}\tbinom{2i}{j}(\p_x^{j}u)(\p_x^{2i-j+1}u)
%\\&
=:u\,A_k(u_x)+B_k(u),
%\end{matrix}
\end{displaymath}
where $B_k:HC^n\to HC^{n-2k}$ is a bounded quadratic operator.
Recall from \nmb!{8.7.4} that we have to solve (where $a$ is a real constant) 
\begin{align*}
u_t&= -A_k\i\bigl(2u_x A_k(u)+u A_k(u_x)+au_{xxx}\bigr)
\\&
= -A_k\i\bigl(2u_x A_k(u)+ A_k(uu_x)-B_k(u)+au_{xxx}\bigr)
\\&
= -uu_x -A_k\i\bigl(2u_x A_k(u)-B_k(u)+au_{xxx}\bigr)
\\&
=: -uu_x +A_k\i C_k(u,a),
\end{align*}
where $u=g_t\o g\i\in \X$, and
where $C_k:HC^n\to HC^{n-2k}$ is a bounded polynomial operator, 
given by 
\begin{align*}
C_k(a,u) &= -2u_x A_k(u)+B_k(u)-au_{xxx}
\\&
= -2u_x A_k(u)
  +\sum_{i=0}^k(-1)^i\sum_{j=1}^{2i}\tbinom{2i}{j}(\p_x^{j}u)(\p_x^{2i-j+1}u)
  -au_{xxx}.
\end{align*}
Note that here we need $2k\ge 3$. 
In \cite{3m} this result was obtained for $k\ge 3/2$.
We put
\begin{align}
&\begin{cases}
&g_t=:v= u\o g
\\
&{\begin{aligned}
v_t&= u_t\o g + (u_x\o g)g_t = u_t\o g + (uu_x)\o g
= A_k\i C_k(a,u)\o g  
\\&
= A_k\i C_k(a,v\o g\i)\o g  
=: \on{pr_2}(D_k\o E_k)(g,v), \quad \text{  where }
\end{aligned}}
\end{cases}
\tag{\nmb:{1}}\\
E_k(a,g,v)&=(g,C_k(a,v\o g\i)\o g), \qquad D_k(g,v)= (g,A_k\i(v\o g\i)\o g).
\notag\end{align}
Now consider the topological group and Banach manifold 
$\on{DiffHC}^n$.

{\it Claim.} The mapping 
$D_k:\on{DiffHC}^n\x HC^{n-2k}\to \on{DiffHC}^n\x HC^n$ is
strongly $C^1$.

Let us assume that we have $C^1$-curves $s\mapsto g(s)\in \on{DiffHC}^n$
and $s\mapsto v(s)\in HC^{n-2k}$.
Then we have:
\begin{align*}
&\p_s \on{pr}_2 D_k(a,g(s),v(s)) = \p_s A_k\i(v\o g\i)\o g
\\&
=A_k\i(v_s\o g\i)\o g 
+ A_k\i\left((v_x\o g\i)(-\frac{g_s\o g\i}{g_x\o g\i})\right)\o g 
\\&\qquad
+(A_k\i(v\o g\i)_x\o g)g_s
\\
&A_k\Bigl(\Bigl(\p_s \on{pr}_2 D_k(a,g(s),v(s))\Bigr)\o g\i\Bigr) =
\\&
= v_s\o g\i - (v\o g\i)_x(g_s\o g\i) + A_k(A_k\i(v\o g\i)_x (g_s\o g\i))
\\&
= v_s\o g\i - (v\o g\i)_x(g_s\o g\i) + (v\o g\i)_x(g_s\o g\i) +
\\&\quad
+\sum_{i=0}^k \sum_{j=0}^{2i-1}\binom{2i}{j}(\p_x^{j+1}(A_k\i(v\o
g\i))\p_x^{2i-j} (g_s\o g\i) \in HC^{n-2k}
\\
&\p_s \on{pr}_2 D_k(a,g(s),v(s))
= A_k\i(v_s\o g\i)\o g
\\&\quad
+\sum_{i=0}^k \sum_{j=0}^{2i-1}\binom{2i}{j}
A_k\i\Big(\big(\p_x^{j+1}(A_k\i(v\o g\i)\big)\p_x^{2i-j} (g_s\o g\i)\Big)\o g
\end{align*}
and by \nmb!{6.12} and \nmb!{6.13} we can conclude that this is continuous in
$a,g,g_s,v,v_s$ jointly and Lipschitz in $g_s$ and $v_s$. 
Thus $D_k$ is strongly $C^1$.

%For $w\in HC^n$ we have
%\begin{displaymath}
%%\begin{matrix}
%&\p_s|_0 (u\o (g+sw)) = (u_x\o g)w
%\\&
%\p_s|_0 (g+sw)\i = -\frac{w\o g\i}{g_x\o g\i}
%\\&
%\p_s|_0 \on{pr}_2 D_k(g+sw,v) =
%\\&
%= \p_s|_0 A_k\i(v\o g\i)\o(g+sw)
%+ \p_s|_0 (A_k\i(v\o (g+sw)\i))\o g
%\\&
%= ((\p_xA_k\i(v\o g\i))\o g)\,w
%-  (A_k\i((v_x\o g\i)\tfrac{w\o g\i}{g_x\o g\i}))\o g
%\\&
%= (A_k\i(v\o g\i)_x.(w\o g\i))\o g
%-  (A_k\i((v\o g\i)_x (w\o g\i)))\o g
%\end{matrix}
%\end{displaymath}
%%Therefore 
%\begin{displaymath}
%%\begin{matrix}
%&A_k((\p_s|_0 \on{pr}_2 D_k(g+sw,v))\o g\i) =
%\\&
%= A_k(A_k\i(v\o g\i)_x.(w\o g\i))
%-  (v\o g\i)_x (w\o g\i)
%\\&
%= %(v\o g\i)_x.(w\o g\i) +
%\sum_{i=0}^{k-1}\tbinom{2k}{i}\p_x^{2i+1}A_k\i(v\o g\i).\p_x^{2k-i}(w\o g\i)
%%-(v\o g\i)_x (w\o g\i) 
%\in HC^{n-2k}
%\end{matrix}
%\end{displaymath}
%%and so $\p_s|_0 \on{pr}_2 D_k(g+sw,v)\in HC^n$. The rest is easier.

{\it Claim.} The mapping 
$E_k:\on{DiffHC}^n\x HC^{n}\to \on{DiffHC}^n\x HC^{n-2k}$ is
strongly $C^1$.

This can be proved in a similar way as the last claim.
%We make the same assumption and compute
%\begin{align*}
%&\p_s \on{pr}_2 E_k(a, g(s), v(s)) = \p_s C_k(a,v\o g\i)\o g
%\\&
%=\p_s \Bigl(-2((v\o g\i)_x\o g) (A_k(v\o g\i)\o g)
%\\&\qquad\qquad
%  +\sum_{i=0}^k(-1)^i\sum_{j=1}^{2i}\tbinom{2i}{j}((\p_x^{j}(v\o g\i))\o g)
%  ((\p_x^{2i-j+1}(v\o g\i))\o g)
%\\&\qquad\qquad
%  -a(v\o g\i)_{xxx}\o g
%\Bigr)
%\end{align*}
%which by \nmb!{6.12} is 

By the two claims equation \thetag 1 can be viewed as the flow equation of
a $C^1$-vector field on the Hilbert manifold 
$\on{DiffHC}^n\x HC^{n}$. 
Here an existence and uniqueness theorem holds.
Since $v=0$ is a
stationary point, there exists an open neighborhood $W_n$ of $(\on{Id},0)$
in $\on{DiffHC}^n\x HC^{n}$ such that for each initial point
$(g_0,v_0)\in W_n$ equation {\nmb|{1}} has a unique solution 
$\on{Fl}^n_t(g_0,v_0)=(g(t),v(t))$
defined and $C^2$ in $t\in (-2,2)$. Note that $v(t)=g_t(t)$, thus $g(t)$ is
even $C^3$ in $t$. Moreover, the solution depends $C^1$ on the initial
data.

We start with the neighborhood 
\begin{displaymath}
W_{2k+1}\subset \on{DiffHC}^{2k+1}\x HC^{2k+1} 
  \supset \on{DiffHC}^n\x HC^{n}\quad \text{  for }n\ge 2k+1 
\end{displaymath}
and consider the neighborhood $V_n:=W_{2k+1}\cap\on{DiffHC}^n\x HC^{n}$ 
of $(\on{Id},0)$.

{\it Claim.} For any initial point $(g_0,v_0)\in V_n$ the solution 
$\on{Fl}^n_t(g_0,v_0)=(g(t),v(t))$ exists, is unique, is $C^2$ in $t\in
(-2,2)$, and depends $C^1$ on the
initial point in $V_n$.

We use induction on $n\ge 2k+1$. For $n=2k+1$ the claim holds since
$V_{2k+1}=W_{2k+1}$. Let $(g_0,v_0)\in V_{2k+2}$ and let 
$\on{Fl}^{2k+2}_t(g_0,v_0)=(\tilde g(t),\tilde v(t))$ be maximally defined for
$t\in (t_1,t_2)\ni 0$.
Suppose for contradiction that $t_2<2$.
Since $(g_0,v_0)\in V_{2k+2}\subset V_{2k+1}$ the curve   
$\on{Fl}^{2k+2}_t(g_0,v_0)=(\tilde g(t),\tilde v(t))$ solves {\nmb|{1}} also in 
$\on{DiffHC}^{2k+1}\x HC^{2k+1}$, thus 
$\on{Fl}^{2k+2}_t(g_0,v_0)=(\tilde g(t),\tilde v(t))=(g(t),v(t))
:=\on{Fl}^{2k+1}_t(g_0,v_0)$ for $t\in (t_1,t_2)\cap(-2,2)$.
By \nmb!{7.3.6}, the expression 
\begin{displaymath}
\tilde J(t)=\tilde J(g,v,t)=g_x(t)^2A_k(u(t))\o g(t)
  =g_x(t)^tA_k(v(t)\o g(t))\o g(t) \tag{\nmb:{2}}
\end{displaymath}
is constant in $t\in (-2,2)$.
Actually, since we used $C^\infty$-theory for deriving this, one
should check it again by differentiating.
Since $u=g_t\o g\i$ we get the following (the exact formulas can be
computed with the help of Fa\`a di Bruno's formula \nmb!{6.1}:
\begin{align*}
u_x&=(g_{tx}\o g\i)(g\i)_x= \frac{g_{tx}}{g_x}\o g\i
\\
\p_x^2u&= (\frac{\p_x^2g_{t}}{g_x^2}-g_{tx}\frac{\p_x^2g}{g_x^3})\o g\i
\\
\p_x(g\i)&= \frac1{g_x}\o g\i
\\
\p_x^2(g\i)\o g&= -\frac{\p_x^2g}{g_x^3} 
\\
\p_x^{2k}(g\i)\o g &= -\frac{\p_x^{2k}g}{g_x^{2k+1}} 
+ \text{ lower order terms in }g
\\
(\p_x^{2k}u)\o g &= \frac{\p_x^{2k}g_t}{g_x^{2k}} -
g_{tx}\frac{\p_x^{2k}g}{g_x^{2k+1}}
+ \text{ lower order terms in }g, g_t=v.
\end{align*}
Thus
\begin{displaymath}
(-1)^k g_x^{2k-1}\tilde J(t)=g_x\p_x^{2k}g_t - g_{tx}\p_x^{2k}g 
+ \text{ lower order terms in }g, g_t=v.
\end{displaymath}
Hence for each $t\in (-2,2)$: 
\begin{align*}
g_x\p_x^{2k}g_t - g_{tx}\p_x^{2k}g 
&= (-1)^kg_x^2\left(g_x^{2k-3}\tilde J(t) + P_k(g,v) \right),\text{  where}
\\
P_k(g,v)&=\frac{Q_k(g,\p_xg,\dots,\p_x^{2k-1}g,v,\p_xv,\dots,\p_x^{2k-1}v)}{g_x^2}
\end{align*}
for a polynomial $Q_k$. Since $\tilde J(t)=\tilde J(0)$ we obtain that
\begin{displaymath}
\left(\frac{\p_x^{2k}g(t)}{g_x(t)} \right)_t 
= (-1)^k\left(g_x^{2k-3}(t)\tilde J(0)+ P_k(g(t),v(t)) \right)
\text{  for all }t\in (-2,2).
\end{displaymath}
This implies
\begin{displaymath}
\frac{\p_x^{2k}g(t)}{g_x(t)} = \frac{\p_x^{2k}g(0)}{g_x(0)} +  
(-1)^k\int_0^t\left(g_x^{2k-3}(s)\tilde J(0)+ P_k(g(s),v(s)) \right)\,ds.
\end{displaymath}
For $t\in(t_1,t_2)$ we have  
\begin{align*}
\p_x^{2k}\tilde g(t) = &\frac{\p_x^{2k}g_0}{\p_xg_0}g_x(t) +  
\tag{\nmb:{3}}\\&
+(-1)^kg_x(t)\int_0^t\left(g_x^{2k-3}(s)\tilde J(0)+ P_k(g(s),v(s)) \right)\,ds.
\end{align*}
Since $(g_0,v_0)\in V_{2k+2}$ we have 
$\tilde J(0)=\tilde J(g_0,v_0,0)\in HC^2$ by \thetag2. 
Since $k\ge 1$, by {\nmb|{3}} we see that $\p_x^{2k}\tilde g(t)\in HC^2$.
Moreover, since $t_2<2$, the limit $\lim_{t\to t_2-}\p_x^{2k}\tilde g(t)$
exists in $HC^2$, so $\lim_{t\to t_2-}\tilde g(t)$ exists in
$HC^{2k+2}$.  As this limit equals $g(t_2)$, we conclude that $g(t_2)\in
\on{DiffHC}^{2k+2}$. 
Now $\tilde v=\tilde g_t$; so we may differentiate both sides of {\nmb|{3}}
in $t$ and obtain similarly that $\lim_{t\to t_2-}\tilde v(t)$ exists in 
$HC^{2k+2}$ and equals $v(t_2)$. But then we can prolong the flow
line $(\tilde g,\tilde v)$ in $\on{DiffHC}^{2k+2}\x HC^{2k+2}$ 
beyond $t_2$, so $(t_1,t_2)$ was not maximal.  

By the same method we can iterate the induction.
\qed\end{demo}

\section*{Appendix \nmb0{A}.  Smooth calculus beyond Banach spaces} 

The traditional differential calculus works 
well for finite dimensional vector spaces and for Banach spaces. For 
more general locally convex spaces 
we sketch here the convenient approach as explained in 
\cite{14} and \cite{21}.
The main difficulty is that the composition of 
linear mappings stops to be jointly continuous at the level of Banach 
spaces, for any compatible topology. 
We use the notation of \cite{21} and this is the
main reference for the whole appendix. We list results in the order in
which one can prove them, without proofs for which we refer to \cite{21}. 
This should explain how to use these results.
Later we also explain the fundamentals about regular infinite dimensional
Lie groups. 

\subsection*{\nmb.{A.1}. Convenient vector spaces} Let $E$ be a 
locally convex vector space. 
A curve $c:\mathbb R\to E$ is called 
{\it smooth} or $C^\infty$ if all derivatives exist and are 
continuous - this is a concept without problems. Let 
$C^\infty(\mathbb R,E)$ be the space of smooth functions. It can be 
shown that $C^\infty(\mathbb R,E)$ does not depend on the locally convex 
topology of $E$, but only on its associated bornology (system of bounded 
sets).

$E$ is said to be a {\it convenient 
vector space} if one of the following equivalent
conditions is satisfied (called $c^\infty$-completeness):
\begin{enumerate}
\item For any $c\in C^\infty(\mathbb R,E)$ the (Riemann-) integral 
       $\int_0^1c(t)dt$ exists in $E$.
\item A curve $c:\mathbb R\to E$ is smooth if and only if $\la\o c$ is 
       smooth for all $\la\in E'$, where $E'$ is the dual consisting 
       of all continuous linear functionals on $E$.
\item Any Mackey-Cauchy-sequence (i.\ e.\  $t_{nm}(x_n-x_m)\to 0$  
       for some $t_{nm}\to \infty$ in $\mathbb R$) converges in $E$. 
       This is visibly a weak completeness requirement.
\end{enumerate}
The final topology with respect to all smooth curves is called the 
$c^\infty$-topology on $E$, which then is denoted by $c^\infty E$. 
For Fr\'echet spaces it coincides with 
the given locally convex topology, but on the space $\Cal D$ of test 
functions with compact support on $\mathbb R$ it is strictly finer.

\subsection*{\nmb.{A.2}. Smooth mappings} Let $E$ and $F$ be locally 
convex vector spaces, and let $U\subset E$ be $c^\infty$-open. 
A mapping $f:U\to F$ is called {\it smooth} or 
$C^\infty$, if $f\o c\in C^\infty(\mathbb R,F)$ for all 
$c\in C^\infty(\mathbb R,U)$.
{\it
The main properties of smooth calculus are the following.
\begin{enumerate}
\item For mappings on Fr\'echet spaces this notion of smoothness 
       coincides with all other reasonable definitions. Even on 
       $\mathbb R^2$ this is non-trivial.
\item Multilinear mappings are smooth if and only if they are 
       bounded.
\item If $f:E\supseteq U\to F$ is smooth then the derivative 
       $df:U\x E\to F$ is  
       smooth, and also $df:U\to L(E,F)$ is smooth where $L(E,F)$ 
       denotes the space of all bounded linear mappings with the 
       topology of uniform convergence on bounded subsets.
\item The chain rule holds.
\item The space $C^\infty(U,F)$ is again a convenient vector space 
       where the structure is given by the obvious injection
$$
C^\infty(U,F)\to \prod_{c\in C^\infty(\mathbb R,U)} C^\infty(\mathbb R,F)
\to \prod_{c\in C^\infty(\mathbb R,U), \la\in F'} C^\infty(\mathbb R,\mathbb R).
$$
\item The exponential law holds:
$$
C^\infty(U,C^\infty(V,G)) \cong C^\infty(U\x V, G)
$$
     is a linear diffeomeorphism of convenient vector spaces. Note 
     that this is the main assumption of variational calculus.
\item A linear mapping $f:E\to C^\infty(V,G)$ is smooth (bounded) if 
       and only if $E \East{f}{} C^\infty(V,G) \East{\on{ev}_v}{} G$ is smooth 
       for each $v\in V$. This is called the smooth uniform 
       boundedness theorem and it is quite applicable.
\end{enumerate}
}

\begin{proclaim}{\nmb.{A.3}. Theorem. {\rm \cite{14},~4.1.19}.}
Let $c:\mathbb R\to E$ be a curve in a convenient vector space $E$. Let 
$\mathcal{V}\subset E'$ be a subset of bounded linear functionals such that 
the bornology of $E$ has a basis of $\sigma(E,\mathcal{V})$-closed sets. 
Then the following are equivalent:
\begin{enumerate}
\item $c$ is smooth
\item There exist locally bounded curves $c^{k}:\mathbb R\to E$ such that
      $\ell\o c$ is smooth $\mathbb R\to \mathbb R$ with $(\ell\o c)^{(k)}=\ell\o
      c^{k}$. 
\end{enumerate}
If $E$ is reflexive, then for any point separating subset
$\mathcal{V}\subset E'$ the bornology of $E$ has a basis of 
$\si(E,\mathcal{V})$-closed subsets, by {\rm\cite{14},~4.1.23}.
\end{proclaim}

\subsection*{\nmb.{A.4}. Counterexamples in infinite dimensions against 
common beliefs on ordinary differential equations}
Let $E:=s$ be the Fr\'echet space of rapidly decreasing sequences; 
note that by the theory of Fourier series we have $s=C^\infty(S^1,\mathbb
R)$.
Consider the continuous linear operator $T:E\to E$ given by
$T(x_0,x_1,x_2,\dots):=(0,1^2x_1,2^2x_2,3^2x_3,\dots)$. 
The ordinary linear differential equation $x'(t)=T(x(t))$ 
with constant coefficients has no solution in $s$ for certain initial 
values. By recursion one sees that 
the general solution should be given by 
$$
x_n(t)=\sum_{i=0}^n
     \left(\tfrac{n!}{i!}\right)^2 x_i(0)\frac{t^{n-i}}{(n-i)!}.
$$
If the initial value is a finite sequence, say $x_n(0)=0$ for $n>N$
and $x_N(0)\ne 0$, then 
\begin{align*}
x_n(t)&=\sum_{i=0}^N
     \left(\tfrac{n!}{i!}\right)^2 x_i(0)\frac{t^{n-i}}{(n-i)!}\\
&=\frac{(n!)^2}{(n-N)!}\;t^{n-N}\sum_{i=0}^N
     \left(\tfrac{1}{i!}\right)^2 x_i(0)\tfrac{(n-N)!}{(n-i)!}\;t^{N-i}\\
|x_n(t)| &\ge \frac{(n!)^2}{(n-N)!}\;|t|^{n-N}
     \left(|x_N(0)|\left(\tfrac1{N!}\right)^2 -
          \sum_{i=0}^{N-1}
     \left(\tfrac{1}{i!}\right)^2 
     |x_i(0)|\tfrac{(n-N)!}{(n-i)!}|t|^{N-i}
     \right)\\
&\ge \frac{(n!)^2}{(n-N)!}\;|t|^{n-N}
     \left(|x_N(0)|\left(\tfrac1{N!}\right)^2 -
          \sum_{i=0}^{N-1}
     \left(\tfrac{1}{i!}\right)^2 |x_i(0)||t|^{N-i}
     \right)\\
\end{align*}
where the first factor does not lie in the space $s$ of rapidly 
decreasing sequences and where the second factor is larger than 
$\ep>0$ for $t$ small enough.
So at least for a dense set of initial values this differential 
equation has no local solution. 

This shows also, that the theorem of 
Frobenius is wrong, in the 
following sense: The vector field $x\mapsto T(x)$ generates a 
1-dimensional subbundle $E$ of the tangent bundle on the open subset 
$s\setminus \{0\}$. It is involutive since it is 1-dimensional. But 
through points representing finite sequences there exist no local integral 
submanifolds ($M$ with $TM=E|M$). Namely, if $c$ were a smooth 
nonconstant curve with $c'(t)=f(t).T(c(t))$ for some smooth 
function $f$, then $x(t):=c(h(t))$ would satisfy $x'(t)=T(x(t))$, where 
$h$ is a solution of $h'(t)=1/f(h(t))$. 

As next example consider 
$E:=\mathbb R^\mathbb N$ and the continuous linear operator 
$T:E\to E$ given by
$T(x_0,x_1,\dots):= (x_1,x_2,\dots )$.
The corresponding differential equation has solutions for every 
initial value $x(0)$, since the coordinates must satisfy the recusive 
relations $x_{k+1}(t)=x_k'(t)$ and hence any smooth functions 
$x_0:\mathbb R\to \mathbb R$ gives rise to a
solution $x(t):=(x_0^{(k)}(t))_k$ with initial value
$x(0)=(x_0^{(k)}(0))_k$. So by Borel's theorem there exist solutions 
to this equation for any initial value and the difference of any two 
functions with same initial value
is an arbitray infinite flat function.
Thus the solutions are far from being unique.
Note that $\mathbb R^\mathbb N$ is a topological direct summand in 
$C^\infty(\mathbb R,\mathbb R)$
via the projection $f\mapsto (f(n))_n$, and hence the same situation
occurs in $C^\infty(\mathbb R,\mathbb R)$.

Let now $E:=C^\infty(\mathbb R,\mathbb R)$ and consider the continuous 
linear operator
$T:E\to E$ given by $T(x):=x'$. Let 
$x:\mathbb R\to C^\infty(\mathbb R,\mathbb R)$ be a
solution of the equation $x'(t)=T(x(t))$. In terms of 
$\hat x:\mathbb R^2\to \mathbb R$
this says 
$\frac{\partial}{\partial t}\hat x(t,s)
     =\frac{\partial}{\partial s}\hat x(t,s)$. 
Hence
$r\mapsto \hat x(t-r,s+r)$ has vanishing derivative everywhere and so 
this function is constant, and in particular 
$x(t)(s)=\hat x(t,s)=\hat x(0,s+t)=x(0)(s+t)$.
Thus we have a smooth solution $x$ uniquely determined by the initial
value $x(0)\in C^\infty(\mathbb R,\mathbb R)$ which even describes a flow 
for the vector field $T$ in the sense of \nmb!{A.6} below.
In general this solution is however not real-analytic, since for any
$x(0)\in C^\infty(\mathbb R,\mathbb R)$, which is not real-analytic in a 
neighborhood of
a point $s$ the composite $\on{ev}_s\o x=x(s+\quad)$ is not real-analytic 
around $0$.

\subsection*{\nmb.{A.5}. Manifolds and vector fields}
In the sequel we shall use smooth manifolds $M$ modelled on 
$c^\infty$-open subsets of convenient vector spaces.
Since we shall need it we also 
include some results on vector fields and their flows.

{\it
Consider vector fields $X_i\in C^\infty(TM)$
and $Y_i\in \Ga(TN)$ for $i=1,2$, and a smooth mapping $f:M\to N$.
If $X_i$ and $Y_i$ are $f$-related for $i=1,2$, i.\ e.\  $Tf\o X_i = 
Y_i\o f$, then also
$[X_1,X_2]$ and $[Y_1,Y_2]$ are $f$-related.
}

In particular if $f:M\to N$ is a local
diffeomorphism (so $(T_xf)\i$ makes sense for each $x\in M$),
then for $Y\in \Ga(T N)$ a vector field 
$f^*Y\in \Ga(T M)$ is defined
by $(f^*Y)(x) = (T_xf)\i.Y(f(x))$. The linear mapping
$f^*:\Ga(T N) \to \Ga(T M)$ is then a Lie algebra 
homomorphism.

\subsection*{\nmb.{A.6}. The flow of a vector field }
Let $X\in \Ga(TM)$ be a vector field. A 
\idx{\it local flow} 
$\on{Fl}^X$ for $X$ is a smooth mapping 
$\on{Fl}^X: M\x\mathbb R\supset U \to M$ defined on a $c^\infty$-open 
neighborhood $U$ of $M\x0$ such that
\begin{enumerate}
\item $\frac d{dt} \on{Fl}^X_t(x)=X(\on{Fl}^X_t(x))$.
\item $\on{Fl}^X_0(x)=x$ for all $x\in M$.
\item $U\cap (\{x\}\x \mathbb R)$ is a connected open interval. 
\item $\on{Fl}^X_{t+s} = \on{Fl}^X_t\o \on{Fl}^X_s$ holds in the following sense. 
     If the right hand side exists then also the left hand side 
     exists and we have equality. Moreover:
     If $\on{Fl}^X_s$ exists, then the existence of both 
     sides is equivalent and they are equal.
\end{enumerate}

{\it Let $X\in \Ga(TM)$ be a 
vector field which admits a local flow $\on{Fl}^X_t$. 
Then for each integral curve $c$ of $X$ we have $c(t)=\on{Fl}^X_t(c(0))$, 
thus there exists a unique maximal flow. 
Furthermore, $X$ is $\on{Fl}^X_t$-related to itself, i.\ e.\  
$T(\on{Fl}^X_t)\o X = X\o \on{Fl}^X_t$. 

Let $X\in \Ga(TM)$ and $Y\in \Ga(TN)$ be $f$-related vector fields for a smooth mapping 
$f:M \to N$ which have local flows $\on{Fl}^X$ and $\on{Fl}^Y$. 
Then we have $f\o \on{Fl}^X_t = \on{Fl}^Y_t\o f$, whenever both
sides are defined. 

Moreover, if $f$ is a diffeomorphism we
have $\on{Fl}^{f^*Y}_t = f\i\o \on{Fl}^Y_t\o f$ in the following sense: If 
one side exists then also the other side exists, and they are equal. 
}

For $f=Id_M$ this implies that if there exists a flow then 
there exists a unique maximal flow $\on{Fl}^X_t$.

\subsection*{ \nmb.{A.7}. The Lie derivative }
There are situations where we do not know that 
the flow of $X$ exists but where we will be able to produce the 
following assumption: 
Suppose that 
$\ph:\mathbb R\x M\supset U\to M$ is a smooth mapping such that 
$(t,x)\mapsto (t,\ph(t,x)=\ph_t(x))$ is a diffeomorphism 
$U\to V$, where $U$ and $V$ are open neighborhoods of 
$\{0\}\x M$ in $\mathbb R\x M$, and such that
$\ph_0=\on{Id}_M$ and $\p_t \ph_t=X\in \Ga(TM)$.
Then again 
$\p_t|_0(\ph_t)^*f = \p_t|_0f\o \ph_t = 
df\o X = X(f)$.

{\it In this situation we have for $Y\in \Ga(TM)$, 
and for a $k$-form $\om\in \Om^k(M)$:
\begin{gather*}
\p_t|_0(\ph_t)^*Y = [X,Y],\\
\p_t|_0 (\ph_t)^*\om = \L_X\om.
\end{gather*}
}

\section*{Appendix \nmb0{B}. Regular infinite dimensional Lie groups}

\subsection*{ \nmb.{B.1}. Lie groups } 
A \idx{\it Lie group} $G$ is a
smooth manifold modelled on $c^\infty$-open subsets of a convenient 
vector space, and a group such that the multiplication
$\mu:G\x G\to G$ and the inversion $\nu:G\to G$ are smooth. 
We shall use the following notation: \newline
$\mu:G\x G\to G$, multiplication, $\mu(x,y) = x.y$. \newline
$\mu_a: G\to G$, left translation, $\mu_a(x) = a.x$.\newline
$\mu^a: G\to G$, right translation, $\mu^a(x) = x.a$.\newline
$\nu: G\to G$, inversion, $\nu(x) = x\i$.\newline
$e\in G$, the unit element. 

{\it The tangent mapping 
$T_{(a,b)}\mu:T_aG\x T_bG \to T_{ab}G$ is given by 
$$T_{(a,b)}\mu.(X_a,Y_b) = T_a(\mu^b).X_a + T_b(\mu_a).Y_b$$
and $T_a\nu:T_aG\to T_{a\i}G$ is given by
$$T_a\nu = - T_e(\mu^{a\i}).T_a(\mu_{a\i}) 
= - T_e(\mu_{a\i}).T_a(\mu^{a\i}).$$ 
}

\subsection*{ \nmb.{B.2}. Invariant vector fields and Lie algebras }
Let $G$ be a (real) Lie group. A vector field $\xi$ on $G$ is
called \idx{\it left invariant}, 
if $\mu_a^*\xi = \xi$ for all $a\in
G$, where $\mu_a^*\xi = T(\mu_{a\i})\o\xi\o\mu_a$.
Since we have $\mu_a^*[\xi,\et] =
[\mu_a^*\xi,\mu_a^*\et]$, the space $\X_L(G)$ of all left invariant
vector fields on $G$ is closed under the Lie bracket, so it is a
sub Lie algebra of $\X(G)$. 
Any left invariant vector field
$\xi$ is uniquely determined by $\xi(e)\in T_eG$, since
$\xi(a)=T_e(\mu_a).\xi(e)$. Thus the Lie algebra $\X_L(G)$ of
left invariant vector fields is linearly isomorphic to $T_eG$,
and on $T_eG$ the Lie bracket on $\X_L(G)$ induces a Lie algebra
structure, whose bracket is again denoted by $[\quad,\quad]$.
This Lie algebra will be denoted as usual by $\mathfrak g$,
sometimes by $\operatorname{Lie}(G)$.

We will also give a name to the isomorphism with the space of left 
invariant vector fields: $L:\mathfrak g\to \X_L(G)$, $X\mapsto L_X$, 
where $L_X(a)= T_e\mu_a.X$. Thus $[X,Y] = [L_X,L_Y](e)$.

Similarly a vector field $\eta$ on $G$ is called 
\idx{\it right invariant}, if
$(\mu^a)^*\eta =\eta$ for all $a\in G$. If $\xi$ is left
invariant, then $\nu^*\xi$ is right invariant.
The right invariant vector fields form a sub Lie algebra
$\X_R(G)$ of $\X(G)$, which is again linearly isomorphic to
$T_eG$ and induces the negative of the Lie algebra structure on $T_eG$. 
We will denote by
$R:\mathfrak g = T_eG\to \X_R(G)$ the isomorphism discussed, which
is given by $R_X(a)=T_e(\mu^a).X$.

{\it If $L_X$ is a left invariant vector
field and $R_Y$ is a right invariant vector field, then $[L_X,R_Y]=0$.
So if the flows of $L_X$ and $R_Y$ exist, they commute.

Let $\ph:G\to H$ be a smooth homomorphism of Lie
groups.
Then $\ph' := T_e\ph:\mathfrak g=T_eG \to \mathfrak
h=T_eH$ is a Lie algebra homomorphism.
}

\subsection*{ \nmb.{B.3}.  One parameter subgroups } 
Let $G$ be a Lie
group with Lie algebra $\mathfrak g$. A \idx{\it one parameter subgroup}
of $G$ is a Lie group homomorphism $\al:(\mathbb R,+) \to G$, i.e.
a smooth curve $\al$ in $G$ with $\al(s+t)=\al(s).\al(t)$, and hence 
$\al(0)=e$.  

Note that a smooth mapping 
$\be:(-\ep,\ep)\to G$ satisfying $\be(t)\be(s)=\be(t+s)$ for $|t|$, 
$|s|$, $|t+s|<\ep$ is the restriction of a one parameter subgroup.
Namely, choose $0<t_0<\ep/2$. Any $t\in\mathbb R$ can be uniquely 
written as $t=N.t_0+t'$ for $0\le t'<t_0$ and $N\in\mathbb Z$. Put
$\al(t)=\be(t_0)^N\be(t')$. The required properties are easy to 
check.

{\it Let $\al:\mathbb R\to G$ be a smooth curve with
$\al(0)=e$. Let $X\in \mathfrak g$. 
Then the following assertions are equivalent.
\begin{enumerate}
\item $\al$ is a one parameter subgroup with $X=\p_t \al(t)$.
\item $\al(t)$ is an integral curve of the left invariant vector 
       field $L_X$, and also an integral curve of the right invariant 
       vector field $R_X$. 
\item $\on{Fl}^{L_X}(t,x):= x.\al(t)$ (or
     $\on{Fl}^{L_X}_t=\mu^{\al(t)}$) is the (unique by \nmb!{A.6}) 
       global flow of $L_X$ in the sense of \nmb!{A.6}.
\item $\on{Fl}^{R_X}(t,x):= \al(t).x$ (or
     $\on{Fl}^{R_X}_t=\mu_{\al(t)}$) is the (unique)
       global flow of $R_X$.
\end{enumerate}
Moreover, each of these properties determines $\al$ uniquely.
}

\subsection*{ \nmb.{B.4}. Exponential mapping }
Let $G$ be a Lie group with Lie algebra $\g$. We say that $G$ admits 
an \idx{\it exponential mapping} if there exists a smooth mapping
$\exp:\mathfrak g \to G$ such that $t\mapsto \exp(tX)$ is the (unique by 
\nmb!{B.3}) 
1-parameter subgroup with tangent vector $X$ at 0. Then we have by 
\nmb!{B.3}
\begin{enumerate}
\item $\on{Fl}^{L_X}(t,x) = x.\exp(tX)$.
\item $\on{Fl}^{R_X}(t,x) = \exp(tX).x$.
\item $\exp(0)=e$ and $T_0\exp = Id:T_0\mathfrak g=\mathfrak g \to
        T_eG=\mathfrak g$ since $T_0\exp.X = \p_t|_0 \exp(0+t.X) =
        \p_t|_0 \on{Fl}^{L_X}(t,e) = X$.
\item Let $\ph:G\to H$ be a smooth
     homomorphism between Lie groups admitting exponential 
     mappings. Then the diagram 
$$\xymatrix{ 
\mathfrak g   \ar[r]^{\ph'} \ar[d]_{\exp^G}
     &  \mathfrak h \ar[d]^{\exp^H}\\
G \ar[r]^{\ph} &   H}  
$$
     commutes, since $t\mapsto \ph(\exp^G(tX))$ is a one parameter
     subgroup of $H$ and $\p_t|_0\ph(\exp^GtX) = \ph'(X)$, so
     $\ph(\exp^GtX)=\exp^H(t\ph'(X))$.
\end{enumerate}

We shall strengthen this notion in \nmb!{B.9} below
and call it a `regular Fr\'echet Lie groups'. 

If $G$ admits an exponential mapping,
it follows from \nmb!{B.4}.\therosteritem3 that $\exp$ is a 
diffeomorphism from a neighborhood of 0 in $\g$ onto a neighborhood 
of $e$ in $G$, if a suitable inverse function theorem is applicable. 
This is true for example for smooth Banach Lie groups, also for gauge 
groups, but it is wrong for diffeomorphism groups.

If $E$ is a Banach space, then in the Banach Lie group $GL(E)$ 
of all bounded linear automorphisms of $E$ the exponential 
mapping is given by the von Neumann series 
$\exp(X)=\sum_{i=0}^\infty\frac1{i!}X^i$. 

If $G$ is connected with exponential 
mapping and $U\subset
\mathfrak g$ is open with $0\in U$, then one may ask whether the group 
generated by $\exp(U)$ equals $G$. 
Note that this is a normal subgroup. So if $G$ is simple, the answer 
is yes.
This is true for connected components of diffeomorphism 
groups and many of their important subgroups.

\subsection*{ \nmb.{B.5}. The adjoint representation } 
Let $G$ be a Lie group with Lie algebra $\g$. 
For $a\in G$ we define $\on{conj}_a:G\to G$ by $\on{conj}_a(x)=axa\i$. It
is called the \idx{\it conjugation} or the \idx{\it inner automorphism}
by $a\in G$. This defines a smooth action of $G$ on itself by 
automorphisms.

The adjoint representation $\Ad:G\to GL(\mathfrak g)\subset L(\g,\g)$ 
is given by 
$\Ad(a)=(\on{conj}_a)'=T_e(\on{conj}_a):\mathfrak g\to \mathfrak g$ for $a\in G$. 
By \nmb!{B.2} $\Ad(a)$ is a Lie algebra homomorphism.
By \nmb!{B.1} we have
$\Ad(a)=T_e(\on{conj}_a) = T_a(\mu^{a\i}).T_e(\mu_a) =
     T_{a\i}(\mu_a).T_e(\mu^{a\i})$.

Finally we define the (lower case) \idx{\it adjoint representation} of
the Lie algebra $\mathfrak g$,
$\ad:\mathfrak g\to \mathfrak g\mathfrak l(\mathfrak g):=L(\mathfrak g,\mathfrak g)$,
by $\ad:=\Ad'=T_e\Ad$.

We shall also use the {\it right} \idx{\it Maurer-Cartan form} 
$\ka^r\in\Om^1(G,\g)$, given by 
$\ka^r_g=T_g(\mu^{g\i}):T_gG\to \g$; similarly the 
{\it left Maurer-Cartan form} $\ka^l\in\Om^1(G,\g)$ is 
given by $\ka^l_g=T_g(\mu_{g\i}):T_gG\to \g$.

{\it
\begin{enumerate}
\item $L_X(a) = R_{\Ad(a)X}(a)$ for $X\in \mathfrak g$
        and $a\in G$.
\item $\ad(X)Y = [X,Y]$ for $X,Y\in \mathfrak g$.
\item $d\Ad = (\ad\o\ka^r).\Ad = 
       \Ad.(\ad\o \ka^l):TG\to L(\g,\g)$.
\end{enumerate}
}

\subsection*{ \nmb.{B.6}. Right actions} 
Let $r:M\x G\to M$ be a right action, so
$\check r:G\to \operatorname{Diff}(M)$ is a group anti-homomorphism. 
We will use the following notation:
$r^a:M\to M$ and $r_x:G\to M$, given by $r_x(a)=r^a(x)=r(x,a)=x.a$.

For any $X\in \mathfrak g$ we define the 
\idx{\it fundamental vector field} $\zeta_X=\ze_X^M\in \X(M)$ by 
$\ze_X(x) = T_e(r_x).X = T_{(x,e)}r.(0_x,X)$.

{\it In this situation the following assertions hold:
\begin{enumerate}
\item $\ze:\mathfrak g\to\X(M)$ is a Lie algebra homomorphism.
\item $T_x(r^a).\ze_X(x) = \ze_{Ad(a\i)X}(x.a)$.
\item $0_M\x L_X \in \X(M\x G)$ is $r$-related to 
       $\ze_X\in \X(M)$.
\end{enumerate}
}

\subsection*{ \nmb.{B.7}. The right and left logarithmic derivatives }
Let $M$
be a manifold and let $f:M\to G$ be a smooth mapping into a Lie
group $G$ with Lie algebra $\mathfrak g$. We define the mapping 
$\de^r f:TM\to \mathfrak g$ by the formula
$$
\de^r f(\xi_x) := T_{f(x)}(\mu^{f(x)\i}).T_xf.\xi_x 
\text{ for }\xi_x\in T_xM.
$$
Then $\de^r f$ is a $\mathfrak g$-valued 1-form on $M$, $\de^r f\in
\Om^1(M;\mathfrak g)$. We call $\de^r f$ the
\idx{\it right logarithmic derivative} of $f$, since for $f:\mathbb
R\to(\mathbb R^+,\cdot)$ we have 
$\de^r f(x).1=\frac{f'(x)}{f(x)}=(\log\o f)'(x)$.

Similarly the
\idx{\it left logarithmic derivative} 
$\de^lf\in\Om^1(M,\mathfrak g)$ 
of a smooth mapping $f:M\to G$ is given by 
$$\de^lf.\xi_x= T_{f(x)}(\mu_{f(x)\i}).T_xf.\xi_x.$$

{\it Let $f,g:M\to G$ be smooth. Then the Leibniz rule 
holds:
$$\de^r(f.g)(x) 
     = \de^r f(x) + \Ad(f(x)).\de^r g(x).$$
Moreover, the differential form $\de^rf\in\Om^1(M;\g)$ 
satifies the `left Maurer-Cartan equation' (left because it stems 
from the left action of $G$ on itself) 
\begin{gather*}
d\de^rf(\xi,\et) - 
[\de^rf(\xi),\de^rf(\et)]^\g=0,\\
\text{ or }\quad d\de^rf - \frac12 [\de^rf,\de^rf]^\g_\wedge=0,
\end{gather*}
where $\xi,\et\in T_x M$, and 
where for $\ph\in\Om^p(M;\g),\ps\in\Om^q(M;\g)$ one puts
$$[\ph,\ps]^\g_\wedge(\xi_1,\dots,\xi_{p+q}) 
     := \frac1{p!q!}\sum_\si\on{sign}(\si)
     [\ph(\xi_{\si1},\dots),\ps(\xi_{\si(p+1)},\dots)]^\g.$$
For the left logarithmic derivative 
the corresponding Leibniz rule is uglier, 
and it satisfies the `right Maurer Cartan equation':
\begin{gather*}
\de^l(fg)(x) = \de^lg(x) + 
Ad(g(x)\i)\de^lf(x),\\
d\de^lf + \frac12 
[\de^lf,\de^lf]^\g_\wedge=0.
\end{gather*}
}

For `regular Lie groups' a converse to this statement holds, see
\cite{21},~40.2.
The proof of this result in infinite dimensions uses 
principal bundle geometry 
for the trivial principal bundle 
$\operatorname{pr_1}:M\x G\to M$ with right principal action. Then 
the submanifolds $\{(x,f(x).g):x\in M\}$ for $g\in G$ form a 
foliation of $M\x G$ whose tangent distribution is complementary to the 
vertical bundle $M\x TG \subseteq T(M\x G)$ and is invariant under the 
principal right $G$-action. So it is the horizontal distribution of a 
principal connection on $M\x G\to G$. Thus this principal connection has 
vanishing curvature which translates into the result for the right
logarithmic derivative.

\subsection*{\nmb.{B.8}}
Let $G$ be a Lie group with Lie algebra $\g$. 
For a closed interval $I\subset \mathbb R$  and for 
$X\in C^\infty(I,\g)$ we consider the ordinary differential equation
\begin{equation*}
\begin{cases} g(t_0)&=e \\
       \p_tg(t)&=T_e(\mu^{g(t)})X(t) = R_{X(t)}(g(t)),
          \quad\text{ or }\ka^r(\p_tg(t)) = X(t),
\end{cases}\tag{1}
\end{equation*}
for local smooth curves $g$ in $G$, where $t_0\in I$.

{\it 
\begin{enumerate}
\item[(2)] Local solution curves $g$ of the differential equation 
       \thetag1 are unique.
\item[(3)]  If for fixed $X$ the differential equation \thetag1 has a 
       local solution near each $t_0\in I$, then it has also a global 
       solution $g\in C^\infty(I,G)$.  
\item[(4)]  If for all $X\in C^\infty(I,\g)$ the differential equation 
       \thetag1 has a local solution near one fixed $t_0\in I$, then 
       it has also a global solution $g\in C^\infty(I,G)$ for each 
       $X$. Moreover, if the local solutions near $t_0$ depend 
       smoothly on the vector fields $X$ then so does the global solution. 
\item[(5)]  The curve $t\mapsto g(t)\i$ is the unique local smooth curve $h$ in 
       $G$ which satifies  
$$
\begin{cases} h(t_0)=e \\
\p_th(t) = T_e(\mu_{h(t)})(-X(t)) = L_{-X(t)}(h(t)),
          \quad\text{ or }\ka^l(\p_th(t)) = -X(t).
\end{cases}$$
\end{enumerate}
}

\subsection*{ \nmb.{B.9}. Regular Lie groups }
If for each $X\in C^\infty(\mathbb R,\g)$ there exists 
$g\in C^\infty(\mathbb R,G)$ satisfying 
\begin{equation*}
\begin{cases} g(0)=e, \\
       \p_tg(t)=T_e(\mu^{g(t)})X(t) = R_{X(t)}(g(t)),\\
          \quad\text{ or }\ka^r(\p_tg(t))= \de^rg(\partial_t) = X(t),
\end{cases}\tag{1}
\end{equation*}
then we write 
\begin{gather*}
\on{evol}^r_G(X) = \on{evol}_G(X):=g(1),\\
\on{Evol}^r_G(X)(t) := \on{evol}_G(s\mapsto tX(ts)) = g(t),
\end{gather*} 
and call it the \idx{\it right evolution} 
of the curve $X$ in $G$. By lemma \nmb!{B.8} the 
solution of the differential equation \thetag1 is unique, 
and for global existence it is sufficient that it has a local solution.
Then 
$$
\on{Evol}^r_G: C^\infty(\mathbb R,\g) \to \{g\in C^\infty(\mathbb R,G):g(0)=e\}
$$
is bijective with inverse the right logarithmic derivative $\de^r$.

The Lie group $G$ is called a
\idx{\it regular Lie group} if 
$\on{evol}^r:C^\infty(\mathbb R,\g)\to G$ exists and is smooth.

We also write 
\begin{gather*}
\on{evol}^l_G(X) = \on{evol}_G(X):=h(1),\\
\on{Evol}^l_G(X)(t) := \on{evol}^l_G(s\mapsto tX(ts)) = h(t),
\end{gather*} 
if $h$ is the (unique) solution of 
\begin{equation*}
\begin{cases} h(0)=e \\
     \p_th(t) = T_e(\mu_{h(t)})(X(t)) = L_{X(t)}(h(t)),\\
          \quad\text{ or }\ka^l(\p_th(t))=\de^lh(\partial_t) = X(t).
\end{cases}\tag{2}
\end{equation*}
Clearly $\on{evol}^l:C^\infty(\mathbb R,\g)\to G$ exists and is 
also smooth if $\on{evol}^r$ does, 
since we have $\on{evol}^l(X)=\on{evol}^r(-X)\i$ by lemma \nmb!{B.8}. 

Let us collect some easily seen properties of the evolution mappings.
If $f\in C^\infty(\mathbb R,\mathbb R)$, then we have 
\begin{align*}
\on{Evol}^r(X)(f(t)) &= \on{Evol}^r(f'.(X\o f))(t).\on{Evol}^r(X)(f(0)),\\
\on{Evol}^l(X)(f(t)) &= \on{Evol}^l(X)(f(0)).\on{Evol}^l(f'.(X\o f))(t).
\end{align*}
If $\ph:G\to H$ is a smooth homomorphism between regular Lie groups 
then the diagram
$$\xymatrix{ 
C^\infty(\mathbb R,\mathfrak g)   \ar[r]^{\ph'_*} \ar[d]_{\on{evol}_G}  &  
     C^\infty(\mathbb R,\mathfrak h) \ar[d]^{\on{evol}_H}\\ 
G \ar[r]^{\ph} &   H}  
$$
commutes, since 
$\p_t\ph(g(t))=T\ph.T(\mu^{g(t)}).X(t)=T(\mu^{\ph(g(t))}).\ph'.X(t)$.

Note that each regular Lie group admits an exponential mapping, 
namely the restriction of $\on{evol}^r$ to the 
constant curves $\mathbb R\to \g$. A Lie group is regular if 
and only if its universal covering group is regular.

Up to now the following statement holds:
\begin{enumerate}
\item[]  All known Lie groups are regular.
\end{enumerate}
Any Banach Lie group is regular since we may consider the time 
dependent right invariant vector field $R_{X(t)}$ on $G$ and its 
integral curve $g(t)$ starting at $e$, which exists and depends 
smoothly on (a further parameter in) $X$. In particular 
finite dimensional Lie groups are regular.  

For diffeomorphism groups the evolution operator 
is just integration of time dependent 
vector fields with compact support. 

\subsection*{\nmb.{B.10}. Extensions of Lie groups }
Let $H$ and $K$ be Lie groups.
A Lie group $G$ is called a smooth \idx{\it extension \ign{  of groups}} of 
$H$ with kernel $K$ if we have a short exact sequence of groups
\begin{equation*}
\{e\} \to K \East{i}{} G \East{p}{} H \to \{e\},\tag{1}
\end{equation*}
such that $i$ and $p$ are smooth
and one of the following two equivalent conditions is satisfied:
\begin{enumerate}
\item[2] $p$ admits a local smooth section $s$ near $e$ (equivalently near 
       any point), and $i$ is initial (i.\ e.\,  any $f$ into $K$ is 
       smooth if and only if $i\o f$ is smooth).
\item $i$ admits a local smooth retraction $r$ near $e$ (equivalently 
       near any point), and $p$ is final (i.\ e.\  $f$ from $H$ is smooth 
       if and only if $f\o p$ is smooth).
\end{enumerate}
Of course by $s(p(x))i(r(x))=x$ the two conditions are equivalent, 
and then $G$ is locally diffeomorphic to $K\x H$ via $(r,p)$ with 
local inverse $(i\o\on{pr}_1).(s\o\on{pr}_2)$.

Not every smooth exact sequence of Lie groups admits local sections 
as required in \therosteritem2. Let for example $K$ be a closed 
linear subspace in a convenient vector space $G$ which is not a 
direct summand, and let $H$ be $G/K$. Then the tangent mapping at 0 
of a local smooth splitting would make $K$ a direct summand.

{\it Let $\{e\} \to K \East{i}{} G \East{p}{} H 
\to \{e\}$ be a smooth extension of Lie groups.
Then $G$ is regular if and only if both $K$ and $H$ are regular. 
}

\subsection*{\nmb.{B.11}. Subgroups of regular Lie groups }
Let $G$ and $K$ be Lie groups, let $G$ be regular and let $i:K\to G$ 
be a smooth homomorphism which is initial (see \nmb!{B.10}) with
$T_ei=i':\mathfrak k\to \mathfrak g$ injective. 
We suspect that $K$ is then regular, but we know a proof for this  
only under the following assumption.
\begin{narrower}
There is an open neighborhood $U\subset G$ of $e$ and a smooth 
mapping $p:U\to E$ into a convenient vector space $E$ such that 
$p\i(0)=K\cap U$ and $p$ constant on left cosets $Kg\cap U$.
\end{narrower}

%\newpage
\bibliographystyle{plain}

\end{document}